\documentclass{amsart}
\usepackage[dvips]{graphicx}
\usepackage{amscd}
\usepackage{amsmath}
\usepackage{amsxtra}
\usepackage{amsfonts}
\usepackage{amssymb}
\newcommand{\frogdefs}{\newcommand{\frogleg}{\ }\newcommand{\frogleft}{(}\newcommand{\frogright}{)}}
\newtheorem{theorem}{Theorem}[section]
\newtheorem{corollary}[theorem]{Corollary}
\newtheorem{lemma}[theorem]{Lemma}

\theoremstyle{definition}

\newtheorem{remark}[theorem]{Remark}

\theoremstyle{remark}
\newtheorem*{acknowledgements}{Acknowledgements}

\newlength{\displayboxwidth}
\renewcommand{\theenumi}{\alph{enumi}}

\numberwithin{equation}{section}
\newlength{\qedskip}
\newlength{\qedadjust}

\newbox\frogdown
\newlength\frogdrop
\def\hookdownarrow{\setlength{\unitlength}{0.4pt}\setbox\frogdown=\hbox to 0pt{\hss $\displaystyle \downarrow $\hss }\setlength{\frogdrop}{2.5\ht\frogdown}\raisebox{0pt}[5\unitlength][\frogdrop]
{\begin{picture}(0,5)(10,0)
\put(5,0){\oval(10,10)[t]}
\end{picture}\lower\ht\frogdown\box\frogdown}}
\def\openone
{\mathchoice
{\hbox{\upshape \small1\kern-3.3pt\normalsize1}}
{\hbox{\upshape \small1\kern-3.3pt\normalsize1}}
{\hbox{\upshape \tiny1\kern-2.3pt\SMALL1}}
{\hbox{\upshape \Tiny1\kern-2pt\tiny1}}}
\makeatletter
\newbox\ipbox
\newcommand{\ip}[2]{\left\langle #1\mathrel{\mathchoice
{\setbox\ipbox=\hbox{$\displaystyle \left\langle\mathstrut #1#2\right\rangle$}
\vrule height\ht\ipbox width0.25pt depth\dp\ipbox}
{\setbox\ipbox=\hbox{$\textstyle \left\langle\mathstrut #1#2\right\rangle$}
\vrule height\ht\ipbox width0.25pt depth\dp\ipbox}
{\setbox\ipbox=\hbox{$\scriptstyle \left\langle\mathstrut #1#2\right\rangle$}
\vrule height\ht\ipbox width0.25pt depth\dp\ipbox}
{\setbox\ipbox=\hbox{$\scriptscriptstyle \left\langle\mathstrut #1#2\right\rangle$}
\vrule height\ht\ipbox width0.25pt depth\dp\ipbox}
} #2\right\rangle}
\newcommand{\diracb}[1]{\left\langle #1\mathrel{\mathchoice
{\setbox\ipbox=\hbox{$\displaystyle \left\langle\mathstrut #1\right.$}
\vrule height\ht\ipbox width0.25pt depth\dp\ipbox}
{\setbox\ipbox=\hbox{$\textstyle \left\langle\mathstrut #1\right.$}
\vrule height\ht\ipbox width0.25pt depth\dp\ipbox}
{\setbox\ipbox=\hbox{$\scriptstyle \left\langle\mathstrut #1\right.$}
\vrule height\ht\ipbox width0.25pt depth\dp\ipbox}
{\setbox\ipbox=\hbox{$\scriptscriptstyle \left\langle\mathstrut #1\right.$}
\vrule height\ht\ipbox width0.25pt depth\dp\ipbox}
}\right. }
\newcommand{\dirack}[1]{\left. \mathrel{\mathchoice
{\setbox\ipbox=\hbox{$\displaystyle \left.\mathstrut #1\right\rangle$}
\vrule height\ht\ipbox width0.25pt depth\dp\ipbox}
{\setbox\ipbox=\hbox{$\textstyle \left.\mathstrut #1\right\rangle$}
\vrule height\ht\ipbox width0.25pt depth\dp\ipbox}
{\setbox\ipbox=\hbox{$\scriptstyle \left.\mathstrut #1\right\rangle$}
\vrule height\ht\ipbox width0.25pt depth\dp\ipbox}
{\setbox\ipbox=\hbox{$\scriptscriptstyle \left.\mathstrut #1\right\rangle$}
\vrule height\ht\ipbox width0.25pt depth\dp\ipbox}
} #1\right\rangle}
\newcommand{\rip}[2]{\left( #1\mathrel{\mathchoice
{\setbox\ipbox=\hbox{$\displaystyle \left(\mathstrut #1#2\right)$}
\vrule height\ht\ipbox width0.25pt depth\dp\ipbox}
{\setbox\ipbox=\hbox{$\textstyle \left(\mathstrut #1#2\right)$}
\vrule height\ht\ipbox width0.25pt depth\dp\ipbox}
{\setbox\ipbox=\hbox{$\scriptstyle \left(\mathstrut #1#2\right)$}
\vrule height\ht\ipbox width0.25pt depth\dp\ipbox}
{\setbox\ipbox=\hbox{$\scriptscriptstyle \left(\mathstrut #1#2\right)$}
\vrule height\ht\ipbox width0.25pt depth\dp\ipbox}
} #2\right)}
\let\subsubsubsectionname\@empty
\newcounter{subsubsubsection}[subsubsection]

\def\l@subsubsubsection{\@tocline{4}{0pt}{1pc}{9pc}{}}

\def\subsubsubsection{\@startsection{subsubsubsection}{4}%
  \z@{.5\linespacing\@plus.7\linespacing}{-.5em}%
  {\normalfont\itshape}}
\AtBeginDocument{%
  \@for\@tempa:=-1,0,1,2,3,4\do{%
    \@ifundefined{r@tocindent\@tempa}{%
      \@xp\gdef\csname r@tocindent\@tempa\endcsname{0pt}}{}%
  }%
}
\def\@writetocindents{%
  \begingroup
  \@for\@tempa:=-1,0,1,2,3,4\do{%
    \immediate\write\@auxout{%
      \string\newlabel{tocindent\@tempa}{%
        \csname r@tocindent\@tempa\endcsname}}%
  }%
  \endgroup}
\makeatother
\def\LaTeXparent#1{}%
\def\ChildStyles#1{}%

\frogdefs
\input cyracc.def

\begin{document}
\title[Cascade algorithm at irregular scaling functions]{Convergence of the cascade algorithm at irregular scaling functions}
\author{Ola~Bratteli}
\address{Department of Mathematics\\
University of Oslo\\
PB 1053 -- Blindern\\
N-0316 Oslo\\
Norway}
\email{bratteli@math.uio.no}
\author{Palle~E.~T.~ Jorgensen}
\address{Department of Mathematics\\
The University of Iowa\\
14 MacLean Hall\\
Iowa City, IA 52242-1419\\
U.S.A.}
\email{jorgen@math.uiowa.edu}
\thanks{Research supported by the University of Oslo and the U.S. National Science Foundation.}
\subjclass{Primary 46L60, 47D25, 42A16, 43A65; Secondary 46L45, 42A65, 41A15}
\keywords{Wavelet, cascade algorithm, refinement operator, representation, orthogonal
expansion, quadrature mirror filter, isometry in Hilbert space}

\begin{abstract}
The spectral properties of the Ruelle transfer operator which arises from a
given polynomial wavelet filter are related to the convergence question for
the cascade algorithm for approximation of the corresponding wavelet scaling function.
\end{abstract}\maketitle
\tableofcontents

\setlength{\displayboxwidth}{\textwidth}\addtolength{\displayboxwidth
}{-2\leftmargini}
\renewcommand{\frogleg}{\\}
\renewcommand{\frogleft}{}
\renewcommand{\frogright}{}

\section{\label{Intro}Introduction}

Two operators from wavelet theory are studied: the refinement operator (alias
the cascade approximation operator) $M$
(see (\ref{eq9})),
and the transfer operator $R$
(see (\ref{eq25})).
In the
case of a compactly supported scaling function $\varphi$, we then consider the
approximation problem
\begin{equation}
\lim_{n\rightarrow\infty}M^{n}\psi^{\left(  0\right)  }
=\varphi
\text{\qquad in }%
L^{2}\left(  \mathbb{R}\right)  , \label{eqIntro.Approx}%
\end{equation}
where $\psi^{\left(  0\right)  }\in L^{2}\left(  \mathbb{R}\right)  $ is
given. A result is proved which relates the spectrum of $R$ to the question of
when $\varphi=\lim_{n\rightarrow\infty}M^{n}\psi^{\left(  0\right)  }$. For
those vectors $\psi^{\left(  0\right)  }$ where the approximation holds, the
rapidity of the approximation is related to the spectral data for $R$.

\section{\label{General}General Theory}

It is well known that compactly supported scaling functions $\varphi$ of a
multiresolution analysis satisfy the functional equation%
\begin{equation}
\varphi\left(  x\right)  =\sqrt{2}\sum_{k=0}^{N}a_{k}\varphi\left(
2x-k\right)  \label{eq1}%
\end{equation}
(at least after an integer translation of $\varphi$) \cite{Dau92,CoRy95}. The
standard requirement that $\left\{  \varphi\left(  \,\cdot\,-k\right)  \mid
k\in\mathbb{Z}\right\}  $ forms an orthonormal set of functions in
$L^{2}\left(  \mathbb{R}\right)  $ implies the conditions
\begin{equation}
\sum_{k\in\mathbb{Z}}\bar{a}_{k}a_{k+2l}=\delta_{l} \label{eq2}%
\end{equation}
for all $l\in\mathbb{Z}$, and the second standard requirement that
$\hat{\varphi}\left(  0\right)  =1$ translates into the condition
\begin{equation}
\sum_{k\in\mathbb{Z}}a_{k}=\sqrt{2}. \label{eq3}%
\end{equation}
If%
\begin{equation}
m_{0}\left(  z\right)  =\sum_{k\in\mathbb{Z}}a_{k}z^{k} \label{eq4}%
\end{equation}
for $z=e^{-it}\in\mathbb{T}$, condition (\ref{eq2}) is equivalent to%
\begin{equation}
\left|  m_{0}\left(  z\right)  \right|  ^{2}+\left|  m_{0}\left(  -z\right)
\right|  ^{2}=2, \label{eq5}%
\end{equation}
and (\ref{eq3}) is equivalent to%
\begin{equation}
m_{0}\left(  1\right)  =\sqrt{2}. \label{eq6}%
\end{equation}
However, in general, orthogonality of
$\left\{ \varphi \left( x-k\right) \mid k\in \mathbb{Z}\right\}$
is a condition which is more restrictive than either one of
the two equivalent conditions (\ref{eq2}) or (\ref{eq5}).
The Fourier transform of (\ref{eq1}) is%
\begin{equation}
\hat{\varphi}\left(  t\right)  =\frac{1}{\sqrt{2}}m_{0}\left(  \frac{t}%
{2}\right)  \hat{\varphi}\left(  \frac{t}{2}\right)  . \label{eq7}%
\end{equation}
We apologize to engineers for using $t$ to denote frequency. Since $\varphi$
has compact support (and then (\ref{eq1}) implies that the support is in
$\left[  0,N\right]  $), $\hat{\varphi}$ is continuous at $0$ and an iteration
of (\ref{eq7}) gives
\begin{equation}
\hat{\varphi}\left(  t\right)  =\prod_{k=1}^{\infty}\left(  \frac{m_{0}\left(
t2^{-k}\right)  }{\sqrt{2}}\right)  . \label{eq8}%
\end{equation}
Since $m_{0}$ is a polynomial, this expansion converges uniformly on compacts.

Now, let $\psi^{\left(  0\right)  }$ be any bounded function of compact
support such that $\widehat{\psi^{\left(  0\right)  }}\left(  0\right)  =1$,
and define by iteration%
\begin{align}
\psi^{\left(  n+1\right)  }\left(  x\right)   &  =\left(  M\psi^{\left(
n\right)  }\right)  \left(  x\right) \label{eq9}\\
&  =\sqrt{2}\sum_{k=0}^{N}a_{k}\psi^{\left(  n\right)  }\left(  2x-k\right)
.\nonumber
\end{align}
Then%
\begin{equation}
\widehat{\psi^{\left(  n+1\right)  }}\left(  t\right)  =\frac{1}{\sqrt{2}%
}m_{0}\left(  \frac{t}{2}\right)  \widehat{\psi^{\left(  n\right)  }}\left(
\frac{t}{2}\right)  , \label{eq10}%
\end{equation}
and hence%
\begin{equation}
\widehat{\psi^{\left(  n\right)  }}\left(  t\right)  =\prod_{k=1}^{n}\left(
\frac{m_{0}\left(  t2^{-k}\right)  }{\sqrt{2}}\right)  \widehat{\psi^{\left(
0\right)  }}\left(  t2^{-n}\right)  . \label{eq11}%
\end{equation}
As a result, $\widehat{\psi^{\left(  n\right)  }}\underset{n\rightarrow\infty
}{\longrightarrow}\hat{\varphi}$, uniformly on compacts, and thus
$\psi^{\left(  n\right)  }\underset{n\rightarrow\infty}{\longrightarrow
}\varphi$ in the distribution sense. In short, if $\left\{  a_{k}\mid
k=0,\dots,N\right\}  $ is a finite set of coefficients satisfying (\ref{eq3})
alone, the refinement equation (\ref{eq1}) possesses a distribution solution
$\varphi$ with $\hat{\varphi}\left(  0\right)  =1$ and compact support in
$\left[  0,N\right]  $. This solution is defined by (\ref{eq8}), and is given
as the distribution limit%
\begin{equation}
\varphi=\lim_{n\rightarrow\infty}M^{n}\psi^{\left(  0\right)  }, \label{eq12}%
\end{equation}
where $\psi^{\left(  0\right)  }$ is any integrable function of compact
support such that $\widehat{\psi^{\left(  0\right)  }}\left(  0\right)  =1$.
When this is used to depict the graph of $\varphi$, it is common to take
$\psi^{\left(  0\right)  }$ to be the Haar function%
\begin{equation}
\psi^{\left(  0\right)  }\left(  x\right)  =%
\begin{cases}
1 & \text{if }0\leq x\leq1,  \\
0 & \text{otherwise;}
\end{cases}%
\label{eq13}%
\end{equation}
see \cite{Coh92}. In this paper, we discuss a variety of choices. The question
of when the convergence in (\ref{eq12}) is stronger than distribution
convergence, and the related question of what regularity properties the limit
$\varphi$ has, have received much attention in the literature; see for example
\cite{Str96}, \cite{Vil94}, \cite{CoDa96}, \cite{Pol92}, and references cited
therein. In all these references, it is also assumed that the condition
(\ref{eq2}), alias (\ref{eq5}), holds. If%
\begin{equation}
\hat{\varphi}_{n}\left(  t\right)  =\prod_{k=1}^{n}\left(  \frac{m_{0}\left(
t2^{-k}\right)  }{\sqrt{2}}\right)  \label{eq14}%
\end{equation}
is the $n$'th partial product in (\ref{eq8}) one then shows,
using (\ref{eq5}),
that%
\begin{align}
\int_{-2^{n}\pi}^{2^{n}\pi}\left|  \hat{\varphi}_{n}\left(  t\right)  \right|
^{2}\,dt  &  =\int_{0}^{2^{n+1}\pi}\left|  \hat{\varphi}_{n}\left(  t\right)
\right|  ^{2}\,dt\label{eq15}\\
&  =\int_{0}^{2^{n}\pi}\left|  \hat{\varphi}_{n-1}\left(  t\right)  \right|
^{2}\frac{1}{2}\left(  \left|  m_{0}\left(  2^{-n}t\right)  \right|
^{2}+\left|  m_{0}\left(  2^{-n}t+\pi\right)  \right|  ^{2}\right)
\,dt\nonumber\\
&  =\int_{0}^{2^{n}\pi}\left|  \hat{\varphi}_{n-1}\left(  t\right)  \right|
^{2}\,dt=\dots=\int_{-\pi}^{\pi}\left|  \hat{\varphi}_{0}\left(  t\right)
\right|  ^{2}\,dt=2\pi.\nonumber
\end{align}
Since the convergence $\hat{\varphi}_{n}\underset{n\rightarrow\infty
}{\longrightarrow}\hat{\varphi}$ is uniform on compacts, it follows that%
\begin{equation}
\left\|  \hat{\varphi}\right\|  _{2}^{2}\leq2\pi. \label{eq16}%
\end{equation}
Hence $\varphi\in L^{2}\left(  \mathbb{R}\right)  $, and $\left\|
\varphi\right\|  _{2}\leq1$. Since $\hat{\varphi}_{n}\left(  \,\cdot\,\right)
\chi_{\left[  -\pi,\pi\right]  }\left(  \,\cdot\,2^{-n}\right)  $ converges,
as $n\rightarrow\infty$, uniformly on compacts, and has constant $L^{2}$-norm
equal to $\sqrt{2\pi}$ by (\ref{eq15}), it follows that this sequence
converges weakly to $\hat{\varphi}$ in $L^{2}\left(  \mathbb{R}\right)  $.
Consequently it converges in $L^{2}$-norm to $\hat{\varphi}$ if and only if%
\begin{equation}
\left\|  \hat{\varphi}\right\|  _{2}^{2}=2\pi, \label{eq17}%
\end{equation}
i.e., if and only if%
\begin{equation}
\left\|  \varphi\right\|  _{2}=1. \label{eq18}%
\end{equation}

It is well known that this is equivalent to several other more practical
conditions, like (\ref{eq19})--(\ref{eq23}) below:
\begin{align}
&  \begin{minipage}[t]{\displayboxwidth}\raggedright$\left\{ \varphi
\left( \,\cdot\,-k\right) \right\} _{k\in\mathbb{Z}}%
$ is an orthonormal set in $L^{2}\left( \mathbb{R}\right) $.\end{minipage}%
\label{eq19}\\
&  \begin{minipage}[t]{\displayboxwidth}\raggedright$\sum_{k\in\mathbb{Z}%
}\left| \hat{\varphi}\left( t+2\pi k\right) \right| ^{2}=1$ for all $t\in
\mathbb{R}$. \end{minipage}\label{eq20}\\
&  \begin{minipage}[t]{\displayboxwidth}\raggedright
The only trigonometric polynomials $\xi$ satisfying \[\xi\left( z\right
) =\frac{1}{2}\sum_{w^{2}=z}\left| m_{0}\left( w\right) \right| ^{2}\xi
\left( w\right) \] are the constants.\end{minipage}\label{eq21}\\
&  \begin{minipage}[t]{\displayboxwidth}\raggedright
There is no nontrivial cycle in \[\left\{ z\in\mathbb{T}\mid\left| m_{0}%
\left( z\right) \right| =\smash{\sqrt{2}}\right\} =\left\{ z\in\mathbb{T}%
\mid\left| m_{0}\left( -z\right) \right| =0\right
\} \] for the doubling map $z\mapsto z^{2}$. \end{minipage} \label{eq22}%
\end{align}
See for example \cite[Theorem 3.3.6]{Hor95} or \cite[Chapter 2]{CoRy95} for
details. We will especially need (\ref{eq21}) in the sequel. By the previous
remarks we also note that this is equivalent to:
\begin{equation}
\begin{minipage}[t]{\displayboxwidth}\raggedright
The cascade algorithm, with \[\psi^{\left( 0\right) }\left( x\right
) =\check{\chi}_{\left[ -\pi,\pi\right] }\left( x\right) =\frac{1}{2\pi}%
\int_{-\pi}^{\pi}e^{itx}\,dt =\frac{1}{\pi x}\sin\left( \pi x\right
) ,\] converges in $L^{2}$-norm to $\varphi$; \end{minipage} \label{eq23}%
\end{equation}
and then $\left\{  \varphi\left(  \,\cdot\,-k\right)  \right\}  $ is an
orthonormal set.

In contrast to the distribution convergence, this latter convergence depends
very sensitively on the choice of initial function $\psi^{\left(  0\right)  }%
$. If for example the $\psi^{\left(  0\right)  }$ starting vector above is
replaced by $\frac{1}{N\pi}\sin\left(  N\pi x\right)  $, then $\left\|
\psi^{\left(  n\right)  }\right\|  _{2}^{2}=N$ for all $n$, so $\psi^{\left(
n\right)  }$ cannot converge to $\varphi$ in norm (although it does so
weakly). The question is then which initial functions $\psi^{\left(  0\right)
}$ can be used. One approach, developed by Strang \cite{Str96}, establishes
$L^{2}$-convergence under general circumstances if $\psi^{\left(  0\right)  }$
is chosen such that%
\begin{equation}
\left\{  \psi^{\left(  0\right)  }\left(  \,\cdot\,-k\right)  \right\}
_{k\in\mathbb{Z}} \label{eq24}%
\end{equation}
is an orthonormal set. Thus the standard choice $\psi^{\left(  0\right)
}=\chi_{\left[  0,1\right]  }$ is included. To describe these ``general
circumstances'' we introduce the Ruelle operator%
\begin{equation}
\left(  R\xi\right)  \left(  z\right)  =\frac{1}{2}\sum_{w^{2}=z}\left|
m_{0}\left(  w\right)  \right|  ^{2}\xi\left(  w\right)  . \label{eq25}%
\end{equation}
We may view $R$ as an operator on any of the spaces%
\begin{equation}
\mathbb{C}\left[  z,z^{-1}\right]  \subset C\left(  \mathbb{T}\right)  \subset
L^{\infty}\left(  \mathbb{T}\right)  \subset L^{2}\left(  \mathbb{T}\right)  .
\label{eq26}%
\end{equation}
and in particular it is clear from (\ref{eq25}), (\ref{eq4}) that $R$ maps any
of these spaces into itself. Since $R\left(  \mathbb{C}\left[  z,z^{-1}%
\right]  \right)  \subset\mathbb{C}\left[  z,z^{-1}\right]  $, the invariance
of $C\left(  \mathbb{T}\right)  $ follows from the Stone--Weierstrass theorem.
This also follows directly from (\ref{eq25}), using the continuity of $m_{0}$.

If $P\left[  n,m\right]  $, $n\leq m$, is the subspace of $\mathbb{C}\left[
z,z^{-1}\right]  $ consisting of trigonometric polynomials of the form
$\sum_{k=n}^{m}b_{k}z^{k}$, we note that
\begin{equation}
R\left(  P\left[  n,m\right]  \right)  \subset P\left[  -\left(  \frac{N-n}%
{2}\right)  ,\left(  \frac{m+N}{2}\right)  \right]  \label{eq27}%
\end{equation}
where $\left[  x\right]  $ is the largest integer ${}\leq x$. Thus any
$P\left[  n,m\right]  $ will ultimately be mapped into $P\left[  -N,N\right]
$ by repeated applications of $R$, so all the spaces%
\begin{equation}
P\left[  -N,N\right]  \subset\mathbb{C}\left[  z,z^{-1}\right]  \subset
C\left(  \mathbb{T}\right)  \subset L^{\infty}\left(  \mathbb{T}\right)
\subset L^{2}\left(  \mathbb{T}\right)  \label{eq28}%
\end{equation}
are invariant under $R$.

Let us study the spectral properties and the norm of $R$ on the various
spaces. Note that $\xi\left(  z\right)  =1$ (the constant function $\openone$)
is an eigenvector of $R$ with eigenvalue $1$ in all the subspaces; and, by
(\ref{eq21}), it is the unique eigenvector up to a scalar with eigenvalue $1$
in $\mathbb{C}\left[  z,z^{-1}\right]  $ if and only if $\left\|
\varphi\right\|  _{2}=1$, i.e., if and only if the cascade (\ref{eq23})
converges in $L^{2}$-norm. In any case, (\ref{eq25}) and (\ref{eq5})
immediately imply that%
\begin{equation}
\left\|  R\right\|  _{\infty\rightarrow\infty}=1, \label{eq29}%
\end{equation}
and hence the spectral radius of $R$, as an operator on the four left-hand
subspaces in (\ref{eq28}), is%
\begin{equation}
\rho_{\infty}\left(  R\right)  =1. \label{eq30}%
\end{equation}

Let us note \emph{en passant} that the behaviour of $R$ as a operator on
$L^{2}\left(  \mathbb{T}\right)  $ is different. {}From (\ref{eq25}), we get
for the $L^{2}\left(  \mathbb{T}\right)  $-adjoint operator:%
\begin{equation}
\left(  R^{\ast}\xi\right)  \left(  z\right)  =\left|  m_{0}\left(  z\right)
\right|  ^{2}\xi\left(  z^{2}\right)  , \label{eq31}%
\end{equation}
and hence%
\begin{align}
\left(  RR^{\ast}\xi\right)  \left(  z\right)   &  =\frac{1}{2}\sum_{w^{2}%
=z}\left|  m_{0}\left(  w\right)  \right|  ^{2}\left(  R^{\ast}\xi\right)
\left(  w\right) \label{eq32}\\
&  =\frac{1}{2}\sum_{w^{2}=z}\left|  m_{0}\left(  w\right)  \right|  ^{4}%
\xi\left(  z\right)  .\nonumber
\end{align}
Thus $RR^{\ast}$ is a multiplication operator, and it follows from (\ref{eq5})
and (\ref{eq6}) that%
\begin{equation}
\left\|  R\right\|  =\left\|  RR^{\ast}\right\|  ^{\frac{1}{2}}=\left(
\frac{1}{2}\left(  \sqrt{2}\right)  ^{4}\right)  ^{\frac{1}{2}}=\sqrt{2},
\label{eq33}%
\end{equation}
so $R$ is not contractive, nor normal, on $L^{2}\left(  \mathbb{T}\right)  $.

What is the spectral radius? {}From (\ref{eq25}) and (\ref{eq31}), we verify
the identities%
\begin{equation}
\left(  R^{n}\xi\right)  \left(  z\right)  =\frac{1}{2^{n\mathstrut}}%
\sum_{w^{2^{n}}=z}\prod_{k=0}^{n-1}\left|  m_{0}\left(  w^{2^{k}}\right)
\right|  ^{2}\xi\left(  w\right)  \label{eq34}%
\end{equation}
and%
\begin{equation}
\left(  R^{\ast\,n}\xi\right)  \left(  z\right)  =\prod_{k=0}^{n-1}\left|
m_{0}\left(  z^{2^{k}}\right)  \right|  ^{2}\xi\left(  z^{2^{n}}\right)  ,
\label{eq35}%
\end{equation}
and hence%
\begin{equation}
\left(  R^{n}R^{\ast\,n}\xi\right)  \left(  z\right)  =\frac{1}{2^{n\mathstrut
}}\sum_{w^{2^{n}}=z}\prod_{k=0}^{n-1}\left|  m_{0}\left(  w^{2^{k}}\right)
\right|  ^{4}\xi\left(  z\right)  \label{eq36}%
\end{equation}
holds. Thus $R^{n}R^{\ast\,n}$ is the operator of multiplication by the
function%
\begin{equation}
\frac{1}{2^{n\mathstrut}}\sum_{w^{2^{n}}=z}\prod_{k=0}^{n-1}\left|
m_{0}\left(  w^{2^{k}}\right)  \right|  ^{4}=p_{n}\left(  z\right)  ,
\label{eq37}%
\end{equation}
and it follows that the spectral radius of $R$ is%
\begin{equation}
\rho_{2}\left(  R\right)  =\lim_{n\rightarrow\infty}\left\|  p_{n}\right\|
_{\infty}^{\frac{1}{2n}}. \label{eq38}%
\end{equation}
Note that
\begin{equation}
p_{n}\left(  z\right)  =\frac{1}{2}\sum_{w^{2}=z}\left|  m_{0}\left(
w\right)  \right|  ^{4}p_{n-1}\left(  w\right)  \label{eq39}%
\end{equation}
and
\begin{equation}
p_{0}\left(  z\right)  =1. \label{eq40}%
\end{equation}
But, by (\ref{eq5}) and (\ref{eq6}), we have%
\begin{equation}
1\leq\frac{1}{2}\sum_{w^{2}=z}\left|  m_{0}\left(  w\right)  \right|  ^{4}%
\leq2, \label{eq41}%
\end{equation}
with both extremities attained. By (\ref{eq39}) and (\ref{eq41}), and
induction, we further have%
\begin{equation}
1\leq p_{n}\left(  z\right)  \leq2^{n}, \label{eq42}%
\end{equation}
and hence, by (\ref{eq38}),%
\begin{equation}
1\leq\rho_{2}\left(  R\right)  \leq\sqrt{2}. \label{eq43}%
\end{equation}

After this diversion into Hilbert space, let us return to the cascade
algorithm. Strang's criterion is based on the following basic lemma: If
$\psi_{1}$, $\psi_{2}$ are $L^{2}$-functions of bounded support, define their
relative polynomial by%
\begin{equation}
p\left(  \psi_{1},\psi_{2}\right)  \left(  z\right)  =\sum_{k}z^{k}%
\int_{\mathbb{R}}\overline{\psi_{1}\left(  x-k\right)  }\psi_{2}\left(
x\right)  \,dx.\label{eq44}%
\end{equation}
Note that%
\[
p\left(  L^{2}\left(  \left[  0,N\right]  \right)  ,L^{2}\left(  \left[
0,N\right]  \right)  \right)  \subset P\left[  -N,N\right]  ,
\]
i.e., $p$ maps $L^{2}\left(  \left[  0,N\right]  \right)  \times L^{2}\left(
\left[  0,N\right]  \right)  $ into a $2N+1$-dimensional space of
trigonometric polynomials.

\begin{lemma}
\label{LemGeneral.1}If $M$ is the cascade transform defined by
\textup{(\ref{eq9}),} then%
\begin{equation}
p\left(  M\psi_{1},M\psi_{2}\right)  =R\left(  p\left(  \psi_{1},\psi
_{2}\right)  \right)  \label{eq45}%
\end{equation}
for any pair $\psi_{1}$, $\psi_{2}$ of $L^{2}$-functions of compact support,
where $R$ is the Ruelle operator \textup{(\ref{eq25}) }on $\mathbb{C}\left[
z,z^{-1}\right]  $.
\end{lemma}

\begin{proof}
If $\xi\left(  z\right)  =\sum_{k}x_{k}z^{k}\in\mathbb{C}\left[
z,z^{-1}\right]  $, its Fourier transform is $\left\langle x_{k}\right\rangle
_{k\in\mathbb{Z}}$, and the Ruelle operator (\ref{eq25}) transforms into
\begin{equation}
\left(  \hat{R}x\right)  _{k}=\sum_{j}\sum_{l}\bar{a}_{j}a_{l}x_{j-l+2k}%
.\label{eq46}%
\end{equation}
This follows from the next computation on%
\begin{equation}
\xi\left(  z\right)  =\sum_{k}x_{k}z^{k}.\label{eq47}%
\end{equation}
We have%
\begin{align}
\left(  R\xi\right)  \left(  z\right)   &  =\frac{1}{2}\sum_{w^{2}=z}%
\overline{m_{0}\left(  w\right)  }m_{0}\left(  w\right)  \xi\left(  w\right)
\label{eq48}\\
&  =\frac{1}{2}\sum_{w^{2}=z}\sum_{jl}\bar{a}_{j}a_{l}w^{-j}w^{l}\sum_{k}%
x_{k}w^{k}\nonumber\\
&  =\frac{1}{2}\sum_{w^{2}=z}\sum_{jlk}\bar{a}_{j}a_{l}x_{k}w^{-j+l+k}%
\nonumber\\
&  =\frac{1}{2}\sum_{w^{2}=z}\sum_{jln}\bar{a}_{j}a_{l}x_{j-l+n}%
w^{n}\nonumber\\
&  =\frac{1}{2}\sum_{n}\sum_{w^{2}=z}w^{n}\sum_{jl}\bar{a}_{j}a_{l}%
x_{j-l+n}.\nonumber
\end{align}
The terms with odd $n$ disappear in the last sum, so%
\begin{equation}
\left(  R\xi\right)  \left(  z\right)  =\sum_{k}z^{k}\sum_{jl}\bar{a}_{j}%
a_{l}x_{j-l+2k},\label{eq49}%
\end{equation}
and (\ref{eq46}) follows.

If $\psi_{1}$, $\psi_{2}$ are $L^{2}$-functions with compact support, define%
\begin{align}
b_{k}  &  =\int_{\mathbb{R}}\overline{\psi_{1}\left(  x-k\right)  }\psi
_{2}\left(  x\right)  \,dx,\label{eq50}\\
c_{k}  &  =\int_{\mathbb{R}}\overline{\left(  M\psi_{1}\right)  \left(
x-k\right)  }\left(  M\psi_{2}\right)  \left(  x\right)  \,dx. \label{eq51}%
\end{align}
To prove (\ref{eq45}), we have to show%
\begin{equation}
c_{k}=\left(  \hat{R}b\right)  _{k}. \label{eq52}%
\end{equation}
We have%
\begin{align}
c_{k}  &  =\int_{\mathbb{R}}\overline{\left(  M\psi_{1}\right)  \left(
x-k\right)  }\left(  M\psi_{2}\right)  \left(  x\right)  \,dx\label{eq53}\\
&  =2\int_{\mathbb{R}}\sum_{mn}\overline{a_{m}\psi_{1}\left(  2x-2k-m\right)
}a_{n}\psi_{2}\left(  2x-n\right)  \,dx\nonumber\\
&  =\int_{\mathbb{R}}\sum_{mn}\bar{a}_{m}a_{n}\overline{\psi_{1}\left(
x-2k-m\right)  }\psi_{2}\left(  x-n\right)  \,dx\nonumber\\
&  =\int_{\mathbb{R}}\sum_{mn}\bar{a}_{m}a_{n}\overline{\psi_{1}\left(
x+n-m-2k\right)  }\psi_{2}\left(  x\right)  \,dx\nonumber\\
&  =\sum_{mn}\bar{a}_{m}a_{n}b_{m-n+2k}=\left(  \hat{R}b\right)
_{k}.\nonumber
\end{align}
This shows (\ref{eq52}), and Lemma \ref{LemGeneral.1} is proved.
\end{proof}

We will list a few more preliminaries on the relative polynomials before we
prove our theorem on spectrum and cascade approximation. As in Lemma
\ref{LemGeneral.1}, we consider functions $\psi$, $\psi_{1}$, $\psi_{2}$ in
$L^{2}\left(  \mathbb{R}\right)  $ of compact support. Usually the support of
the functions in question will be assumed contained in a fixed interval
$\left[  0,N\right]  \subset\mathbb{R}$. Consider $\xi\left(  z\right)
=\sum_{k}\xi_{k}z^{k}\in\mathbb{C}\left[  z,z^{-1}\right]  \subset C\left(
\mathbb{T}\right)  $, $\xi_{k}=\tilde{\xi}\left(  k\right)  $ denoting the
Fourier coefficients. Defining%
\begin{equation}
\left(  \xi\ast\psi\right)  \left(  x\right)  :=\sum_{k\in\mathbb{Z}}\xi
_{k}\psi\left(  x-k\right)  , \label{eqOla2.54}%
\end{equation}
we note that%
\begin{equation}
\left(  \xi\ast\psi\right)  \sphat\left(  t\right)  =\xi\left(  e^{-it}%
\right)  \hat{\psi}\left(  t\right)  , \label{eqOla2.55}%
\end{equation}
where ${}\sphat\;$ as usual refers to the $\mathbb{R}$-Fourier transform. The
sesquilinear operator $p\left(  \psi_{1},\psi_{2}\right)  $ in (\ref{eq44})
will be viewed as a quadratic form%
\begin{equation}
p\colon L^{2}\left(  \mathbb{R}\right)  _{c}\times L^{2}\left(  \mathbb{R}%
\right)  _{c}\longrightarrow\mathbb{C}\left[  z,z^{-1}\right]  \subset
C\left(  \mathbb{T}\right)  \subset L^{2}\left(  \mathbb{T}\right)  ,
\label{eqOla2.56}%
\end{equation}
with the subscript $c$ standing for compact support.

Positivity of this $C\left(  \mathbb{T}\right)  $-valued form $p$ follows from
the following identity:%
\begin{align}
p\left(  \psi,\psi\right)  \left(  z\right)   &  =\sum_{k\in\mathbb{Z}}%
z^{k}\int_{\mathbb{R}}\overline{\psi\left(  x-k\right)  }\psi\left(  x\right)
\,dx\label{eqOla2.57}\\
&  =\sum_{k}z^{k}\frac{1}{2\pi}\int_{0}^{2\pi}e_{k}\left(  t\right)
\operatorname*{PER}\left(  \left|  \hat{\psi}\right|  ^{2}\right)  \left(
e^{it}\right)  \,dt\nonumber\\
&  =\operatorname*{PER}\left(  \left|  \hat{\psi}\right|  ^{2}\right)  \left(
z\right)  ,\nonumber
\end{align}
with the convention $z=e^{-it}$, $e_{k}\left(  t\right)  =z^{k}=e^{-ikt}$, and

\[
\operatorname*{PER}\left(  f\right)  \left(  e^{it}\right)  =\sum
_{n\in\mathbb{Z}}f\left(  t+2\pi n\right)  .
\]

Similarly we derive the formula%
\begin{equation}
p\left(  \psi_{1},\psi_{2}\right)  \left(  z\right)  =\operatorname*{PER}%
\left(  \Bar{\Hat{\psi}}_{1}\hat{\psi}_{2}\right)  \left(  z\right)
.\label{eqpPER}%
\end{equation}
{}From this it is immediate that the relative polynomial $p$ has the
sesquilinearity property%
\begin{equation}
p\left(  \psi_{1},\psi_{2}\right)  \left(  z\right)  =\overline{p\left(
\psi_{2},\psi_{1}\right)  \left(  z\right)  },\label{eqOla2.59}%
\end{equation}
and that it is densely defined on $L^{2}\left(  \mathbb{R}\right)  \times
L^{2}\left(  \mathbb{R}\right)  $. We next show that when $\psi_{1}\in
L^{2}\left(  \mathbb{R}\right)  _{c}$ is given, then $p\left(  \psi
_{1},\,\cdot\,\right)  $ extends boundedly to $L^{2}\left(  \mathbb{R}\right)
$ in the second variable, and that it becomes a module mapping in that
variable relative to $\mathbb{C}\left[  z,z^{-1}\right]  \subset L^{\infty
}\left(  \mathbb{T}\right)  $. This is actually slightly more than what we
need (see the proof of (\ref{ThmGeneral.2(2)}) $\Rightarrow$
(\ref{ThmGeneral.2(1)}) in Theorem \ref{ThmGeneral.2}). On the other hand, one
may use the Zak transform \cite{Dau92} to prove a stronger continuity
estimate: $p$ may be extended to a sesquilinear function $L^{2}\left(
\mathbb{R}\right)  \times L^{2}\left(  \mathbb{R}\right)  \rightarrow
L^{1}\left(  \mathbb{T}\right)  $ such that%
\[
\left\|  p\left(  \psi_{1},\psi_{2}\right)  \right\|  _{1}\leq\left\|
\psi_{1}\right\|  _{2}\cdot\left\|  \psi_{2}\right\|  _{2}.
\]
We turn to the details.

\begin{lemma}
\label{LemPropQuad}The relative polynomial 
$p$ has the following properties:

\begin{enumerate}
\item \label{LemPropQuad(1)}Norm estimate:%
\begin{align}
\left\|  p\left(  \psi_{1},\psi_{2}\right)  \right\|  _{L^{2}\left(
\mathbb{T}\right)  }  &  =\frac{1}{2\pi}\int_{0}^{2\pi}\left|
\operatorname*{PER}\left(  \Bar{\Hat{\psi}}_{1}\hat{\psi}_{2}\right)  \left(
t\right)  \right|  ^{2}\,dt\label{eqOla2.60}\\
&  \leq\left\|  \operatorname*{PER}\left|  \hat{\psi}_{1}\right|
^{2}\right\|  _{\infty}\cdot\left\|  \psi_{2}\right\|  _{L^{2}\left(
\mathbb{R}\right)  }^{2}.\nonumber
\end{align}

\item \label{LemPropQuad(2)}Module property:%
\begin{equation}
p\left(  \psi_{1},\xi\ast\psi_{2}\right)  \left(  z\right)  =\xi\left(
z\right)  p\left(  \psi_{1},\psi_{2}\right)  \left(  z\right)  .
\label{eqOla2.61}%
\end{equation}

\item \label{LemPropQuad(3)}Ruelle operator covariance:
\begin{equation}
M_{m_{0}}\left(  \xi\ast\psi\right)  =M_{m_{0}\xi}\left(  \psi\right)
\label{eqOla2.62}%
\end{equation}
and thus%
\begin{equation}
R\left(  p\left(  \psi_{1},\xi\ast\psi_{2}\right)  \right)  \left(  z\right)
=p\left(  M_{m_{0}}\psi_{1},M_{\xi m_{0}}\psi_{2}\right)  \left(  z\right)  ,
\label{eqOla2.63}%
\end{equation}
where $\xi m_{0}$ is the pointwise product.

\item \label{LemPropQuad(4)}Cascade smoothing:%
\begin{equation}
\left(  M^{n}\left(  \xi\ast\psi\right)  \right)  \sphat\left(  t\right)
=\xi\left(  e^{-\frac{it}{2^{n\mathstrut}}}\right)  \widehat{\left(  M^{n}%
\psi\right)  }\left(  t\right)  , \label{eqOla2.64}%
\end{equation}
where $M=M_{m_{0}}$ is the cascade operator.
\end{enumerate}
\end{lemma}

\begin{proof}
\emph{Ad} (\ref{LemPropQuad(1)}): By Parseval's identity in $L^{2}\left(
\mathbb{T}\right)  $,%
\begin{align}
\left\|  p\left(  \psi_{1},\psi_{2}\right)  \right\|  _{L^{2}\left(
\mathbb{T}\right)  }^{2}  &  =\frac{1}{2\pi}\int_{0}^{2\pi}\left|  p\left(
\psi_{1},\psi_{2}\right)  \left(  e^{-it}\right)  \right|  ^{2}%
\,dt\label{eqOla2.65}\\
&  =\sum_{k\in\mathbb{Z}}\left|  \int_{\mathbb{R}}\overline{\psi_{1}\left(
x-k\right)  }\psi_{2}\left(  x\right)  \,dx\right|  ^{2}\nonumber\\
&  =\sum_{k\in\mathbb{Z}}\left|  \frac{1}{2\pi}\int_{\mathbb{R}}e_{k}\left(
t\right)  \Bar{\Hat{\psi}}_{1}\left(  t\right)  \hat{\psi}_{2}\left(
t\right)  \,dt\right|  ^{2}\nonumber\\
&  =\sum_{k\in\mathbb{Z}}\left|  \frac{1}{2\pi}\int_{0}^{2\pi}e_{k}\left(
t\right)  \operatorname*{PER}\left(  \Bar{\Hat{\psi}}_{1}\hat{\psi}%
_{2}\right)  \left(  t\right)  \,dt\right|  ^{2}\nonumber\\
&  =\frac{1}{2\pi}\int_{0}^{2\pi}\left|  \operatorname*{PER}\left(  \Bar
{\Hat{\psi}}_{1}\hat{\psi}_{2}\right)  \right|  ^{2}\,dt\nonumber\\
&  \leq\frac{1}{2\pi}\int_{0}^{2\pi}\left(  \operatorname*{PER}\left|
\hat{\psi}_{1}\right|  ^{2}\cdot\operatorname*{PER}\left|  \hat{\psi}%
_{2}\right|  ^{2}\right)  \,dt\nonumber\\
&  \leq\left\|  \operatorname*{PER}\left|  \hat{\psi}_{1}\right|
^{2}\right\|  _{\infty}\cdot\frac{1}{2\pi}\int_{0}^{2\pi}\operatorname*{PER}%
\left|  \hat{\psi}_{2}\right|  ^{2}\,dt\nonumber\\
&  =\left\|  \operatorname*{PER}\left|  \hat{\psi}_{1}\right|  ^{2}\right\|
_{\infty}\cdot\frac{1}{2\pi}\int_{\mathbb{R}}\left|  \hat{\psi}_{2}\right|
^{2}\,dt\nonumber\\
&  =\left\|  \operatorname*{PER}\left|  \hat{\psi}_{1}\right|  ^{2}\right\|
_{\infty}\cdot\int_{\mathbb{R}}\left|  \psi_{2}\left(  x\right)  \right|
^{2}\;dx.\nonumber
\end{align}

\emph{Ad} (\ref{LemPropQuad(2)}): We have%
\begin{align}
p\left(  \psi_{1},\xi\ast\psi_{2}\right)  \left(  z\right)   &  =\sum
_{k\in\mathbb{Z}}z^{k}\int_{\mathbb{R}}\overline{\psi_{1}\left(  x-k\right)
}\left(  \xi\ast\psi_{2}\right)  \left(  x\right)  \,dx\label{eqOla2.66}\\
&  =\sum_{k}\sum_{l}\xi_{l}z^{k}\int_{\mathbb{R}}\overline{\psi_{1}\left(
x-k\right)  }\psi_{2}\left(  x-l\right)  \,dx\nonumber\\
&  =\sum_{k}\sum_{l}\xi_{l}z^{l}z^{k-l}\int_{\mathbb{R}}\overline{\psi
_{1}\left(  x+l-k\right)  }\psi_{2}\left(  x\right)  \,dx\nonumber\\
&  =\xi\left(  z\right)  p\left(  \psi_{1},\psi_{2}\right)  \left(  z\right)
.\nonumber
\end{align}

\emph{Ad} (\ref{LemPropQuad(3)}): Since%
\begin{equation}
R\left(  p\left(  \psi_{1},\xi\ast\psi_{2}\right)  \right)  =p\left(
M\psi_{1},M\left(  \xi\ast\psi_{2}\right)  \right)  \label{eqOla2.67}%
\end{equation}
by Lemma \ref{LemGeneral.1}, we need only calculate $M\left(  \xi\ast\psi
_{2}\right)  $, where $M$ is the cascade operator $M=M_{m_{0}}$ given by the
low-pass filter $m_{0}$.%
\begin{align}
M\left(  \xi\ast\psi_{2}\right)  \left(  x\right)   &  =\sqrt{2}\sum_{k}%
a_{k}\left(  \xi\ast\psi_{2}\right)  \left(  2x-k\right) \label{eqOla2.68}\\
&  =\sqrt{2}\sum_{k}a_{k}\sum_{l}\xi_{l}\psi_{2}\left(  2x-k-l\right)
\nonumber\\
&  =\sqrt{2}\sum_{n}\sum_{k}a_{k}\xi_{n-k}\psi_{2}\left(  2x-n\right)
\nonumber\\
&  =\sum_{n}\left(  m_{0}\xi\right)  \sptilde\left(  n\right)  \,\psi
_{2}\left(  2x-n\right) \nonumber\\
&  =M_{m_{0}\xi}\psi_{2}\left(  x\right)  ,\nonumber
\end{align}
where $M_{m_{0}\xi}$ is the cascade operator corresponding to the product
filter $\left(  m_{0}\xi\right)  \left(  z\right)  =m_{0}\left(  z\right)
\xi\left(  z\right)  $.

\emph{Ad} (\ref{LemPropQuad(4)}): From the definition of the cascade operator
$M$, we get $\left(  M\psi\right)  \sphat\left(  t\right)  =\frac{1}{\sqrt{2}%
}m_{0}\left(  \frac{t}{2}\right)  \hat{\psi}\left(  \frac{t}{2}\right)  $. Now
apply this to $\xi\ast\psi$, and use $\left(  \xi\ast\psi\right)
\sphat\left(  t\right)  =\xi\left(  e^{-it}\right)  \hat{\psi}\left(
t\right)  $. Iteration yields%
\begin{align}
\left(  M^{n}\left(  \xi\ast\psi\right)  \right)  \sphat\left(  t\right)   &
=\prod_{k=1}^{n}\frac{1}{\sqrt{2}}m_{0}\left(  \frac{t}{2^{k\mathstrut}%
}\right)  \left(  \xi\ast\psi\right)  \sphat\left(  \frac{t}{2^{n\mathstrut}%
}\right) \label{eqOla2.69}\\
&  =\xi\left(  e^{-\frac{it}{2^{n\mathstrut}}}\right)  \prod_{k=1}^{n}\frac
{1}{\sqrt{2}}m_{0}\left(  \frac{t}{2^{k\mathstrut}}\right)  \hat{\psi}\left(
\frac{t}{2^{n\mathstrut}}\right) \nonumber\\
&  =\xi\left(  e^{-\frac{it}{2^{n\mathstrut}}}\right)  \left(  M^{n}%
\psi\right)  \sphat\left(  t\right)  .%
\settowidth{\qedskip}{$\displaystyle
\left(  M^{n}\left(  \xi\ast\psi\right)  \right)  \sphat\left(  t\right)
=\xi\left(  e^{-\frac{it}{2^{n\mathstrut}}}\right)  \left(  M^{n}\psi
\right)  \sphat\left(  t\right)  .$}
\settowidth{\qedadjust}{$\displaystyle
\left(  M^{n}\left(  \xi\ast\psi\right)  \right)  \sphat\left(  t\right)
=\xi\left(  e^{-\frac{it}{2^{n\mathstrut}}}\right)  \prod_{k=1}^{n}\frac
{1}{\sqrt{2}}m_{0}\left(  \frac{t}{2^{k\mathstrut}}\right)  \hat{\psi}\left(
\frac{t}{2^{n\mathstrut}}\right)  $}
\addtolength{\qedadjust}{\textwidth}
\addtolength{\qedskip}{-0.5\qedadjust}
\rlap{\hbox to-\qedskip{\hfil$\qedsymbol$}\hss}%
\nonumber
\end{align}%
\renewcommand{\qed}{}%
\end{proof}

As a corollary, we note that

\begin{corollary}
\label{CorpLxLC}%
\begin{equation}
p\left(  L^{2}\left(  \mathbb{R}\right)  _{c}\times L^{2}\left(
\mathbb{R}\right)  _{c}\right)  =\mathbb{C}\left[  z,z^{-1}\right]  .
\label{eqOla2.70}%
\end{equation}
\end{corollary}

\begin{proof}
We have already commented on one inclusion in (\ref{eqOla2.56}), and the
second follows from%
\[
\xi=\xi\cdot p\left(  \psi^{\left(  0\right)  },\psi^{\left(  0\right)
}\right)  =p\left(  \psi^{\left(  0\right)  },\xi\ast\psi^{\left(  0\right)
}\right)  ,
\]
where $\psi^{\left(  0\right)  }\in L^{2}\left(  \mathbb{R}\right)  _{c}$ is
chosen such that $p\left(  \psi^{\left(  0\right)  },\psi^{\left(  0\right)
}\right)  =\openone$.

(Equivalently, $\operatorname*{PER}\left(  \overline{\widehat{\psi^{\left(
0\right)  }}}\widehat{\psi^{\left(  0\right)  }}\right)  =\operatorname*{PER}%
\left(  \left|  \widehat{\psi^{\left(  0\right)  }}\right|  ^{2}\right)
\equiv1$.)
\end{proof}

\begin{remark}
\label{RemPosEl}The positive elements in $\mathbb{C}\left[  z,z^{-1}\right]  $
\textup{(}i.e., pointwise nonnegative as functions on $\mathbb{T}$\textup{)}
are of the form $\left|  \xi\right|  ^{2}=\bar{\xi}\xi$ by a theorem of Riesz
\cite[Chapter 6]{Dau92}, \cite[p.\ 181]{Akh65}, and we get $\left|
\xi\right|  ^{2}=p\left(  \xi\ast\psi^{\left(  0\right)  },\xi\ast
\psi^{\left(  0\right)  }\right)  $.
\end{remark}

We are now ready to prove the main theorem of convergence of the cascade
algorithm, which is a version of Theorem 4 in \cite{Str96}.

\begin{theorem}
\label{ThmGeneral.2}Let $a_{0},a_{1},\dots,a_{N}$ be complex numbers
satisfying \textup{(\ref{eq2})--(\ref{eq3})} and let $\varphi$ be the
associated scaling function defined by \textup{(\ref{eq8}).} Identify the
Ruelle operator $R$ given in \textup{(\ref{eq25})} with its restriction to
$P\left[  -N,N\right]  $ \textup{(}or to any $P\left[  n,m\right]  $ with
$n\leq-N$, $m\geq N\,$\textup{).} The following conditions are
equivalent.\renewcommand
{\theenumi}{\roman{enumi}}

\begin{enumerate}
\item \label{ThmGeneral.2(1)}$R$ has $1$ as a simple eigenvalue and $\left|
\lambda\right|  <1$ for all other eigenvalues $\lambda$ of $R$.

\item \label{ThmGeneral.2(2)}If $\psi^{\left(  0\right)  }\in L^{2}\left(
\mathbb{R}\right)  $ is a function with compact support such that $\left\{
\psi^{\left(  0\right)  }\left(  \,\cdot\,-k\right)  \right\}  _{k\in
\mathbb{Z}}$ is an orthonormal set and $\widehat{\psi^{\left(  0\right)  }%
}\left(  0\right)  =1$ then%
\begin{equation}
\lim_{n\rightarrow\infty}\left\|  \varphi-M^{n}\psi^{\left(  0\right)
}\right\|  _{2}=0. \label{eq54}%
\end{equation}
\end{enumerate}
\end{theorem}

\begin{proof}
(\ref{ThmGeneral.2(1)}) $\Rightarrow$ (\ref{ThmGeneral.2(2)}). As an aside,
remark that $R$ having $1$ as a simple eigenvalue means here that the
corresponding eigenspace is one-dimensional. But since $\left\|  R\right\|
=1$, and hence $n\mapsto\left\|  R^{n}\right\|  ^{\prime}$ is bounded in any
equivalent norm $\left\|  \,\cdot\,\right\|  ^{\prime}$ on the linear
operators on $P\left[  -N,N\right]  \cong\mathbb{C}^{2N+1}$, it follows from
Jordan's theorem that the multiplicity of $1$ in the characteristic polynomial
is $1$ too.

Let us view $P\left[  -N,N\right]  =\left\{  \sum_{k=-N}^{N}x_{k}%
z^{k}\right\}  $ as the space of sequences $x=\left\langle x_{k}\right\rangle
_{k=-N}^{N}$. Since%
\begin{align}
\left(  R\xi\right)  \left(  1\right)   &  =\frac{1}{2}\left(  \left|
m_{0}\left(  1\right)  \right|  ^{2}\xi\left(  1\right)  +\left|  m_{0}\left(
-1\right)  \right|  ^{2}\xi\left(  -1\right)  \right) \label{eqOla2.73}\\
&  =\xi\left(  1\right) \nonumber
\end{align}
for all $\xi\in\mathbb{C}\left[  z,z^{-1}\right]  $, we have
\begin{equation}
\sum_{k}\left(  \hat{R}x\right)  \left(  k\right)  =\sum_{k}x\left(  k\right)
. \label{eqOla2.74}%
\end{equation}
But as $\hat{R}\left(  \delta_{0}\right)  =\delta_{0}$, $\delta_{0}$ is the
unique eigenvector of $\hat{R}$ corresponding to eigenvalue $1$, and as the
functional $x\mapsto\sum_{k}x\left(  k\right)  $ is preserved by $\hat{R}$, it
follows from (\ref{ThmGeneral.2(1)}) that
\begin{equation}
\lim_{n\rightarrow\infty}\hat{R}^{n}x=\left(  \sum_{k}x\left(  k\right)
\right)  \delta_{0} \label{eq55}%
\end{equation}
for all finite sequences $x$. Thus, by Lemma \ref{LemGeneral.1},%
\begin{align}
p\left(  \varphi,M^{n}\psi^{\left(  0\right)  }\right)   &  =p\left(
M^{n}\varphi,M^{n}\psi^{\left(  0\right)  }\right) \label{eqOla2.76}\\
&  =R^{n}\left(  p\left(  \varphi,\psi^{\left(  0\right)  }\right)  \right)
\nonumber\\
&  \underset{n\rightarrow\infty}{\longrightarrow}\left(  \sum_{k}%
\int_{\mathbb{R}}\overline{\varphi\left(  x-k\right)  }\psi^{\left(  0\right)
}\left(  x\right)  \,dx\right)  \cdot\openone.\nonumber
\end{align}
But the two assumptions on $\psi^{\left(  0\right)  }$ imply the so-called
Strang--Fix condition \cite{Vil94,Str96}%
\begin{equation}
\sum_{k}\psi^{\left(  0\right)  }\left(  x+k\right)  =1 \label{eqOla2.77}%
\end{equation}
for all $x$, and hence by the above,%
\begin{equation}
p\left(  \varphi,M^{n}\psi^{\left(  0\right)  }\right)  \longrightarrow
\int_{\mathbb{R}}\overline{\varphi\left(  x\right)  }\,dx=\overline
{\hat{\varphi}\left(  0\right)  }=1. \label{eqOla2.78}%
\end{equation}
In particular this means that%
\begin{equation}
\lim_{n\rightarrow\infty}\ip{\varphi}{M^{n}\psi^{\left( 0\right) }}=1.
\label{eqOla2.79}%
\end{equation}
But by (\ref{eq45}),%
\begin{equation}
p\left(  M^{n}\psi^{\left(  0\right)  },M^{n}\psi^{\left(  0\right)  }\right)
=R^{n}\left(  p\left(  \psi^{\left(  0\right)  },\psi^{\left(  0\right)
}\right)  \right)  . \label{eqOla2.80}%
\end{equation}
Since $p\left(  \psi^{\left(  0\right)  },\psi^{\left(  0\right)  }\right)
=\openone$ by orthonormality, we conclude that%
\begin{equation}
p\left(  M^{n}\psi^{\left(  0\right)  },M^{n}\psi^{\left(  0\right)  }\right)
=R^{n}\openone=\openone, \label{eqOla2.81}%
\end{equation}
and in particular,%
\begin{equation}
\left\|  M^{n}\psi^{\left(  0\right)  }\right\|  _{2}^{2}=1 \label{eqOla2.82}%
\end{equation}
for all $n$. Also%
\begin{equation}
\left\|  M^{n}\varphi\right\|  _{2}^{2}=\left\|  \varphi\right\|  _{2}^{2}=1,
\label{eqOla2.83}%
\end{equation}
so, finally,
\begin{align}
\left\|  \varphi-M^{n}\psi^{\left(  0\right)  }\right\|  _{2}^{2}  &
=\left\|  \varphi\right\|  _{2}^{2}-2\operatorname{Re}\ip{\varphi}{M^{n}%
\psi^{\left( 0\right) }}+\left\|  \psi^{\left(  0\right)  }\right\|  _{2}%
^{2}\label{eqOla2.84}\\
&  \underset{n\rightarrow\infty}{\longrightarrow}1-2+1=0.\nonumber
\end{align}

(\ref{ThmGeneral.2(2)}) $\Rightarrow$ (\ref{ThmGeneral.2(1)}). We now assume
cascade convergence in the sense (\ref{ThmGeneral.2(2)}), i.e.,
\begin{equation}
\left\|  \varphi-M^{n}\psi^{\left(  0\right)  }\right\|  _{2}\longrightarrow0
\label{eqOla2.85}%
\end{equation}
for the initial vectors $\psi^{\left(  0\right)  }$ which are specified in
(\ref{ThmGeneral.2(2)}). The object is to derive from this the spectral
picture for $R$ as specified in (\ref{ThmGeneral.2(1)}), and $R$ will be
identified with its restriction to $P\left[  -N,N\right]  $ as mentioned. Of
course $P\left[  -N,N\right]  \subset C\left(  \mathbb{T}\right)  $, and $R$
is also, by (\ref{eq28}), an operator mapping $C\left(  \mathbb{T}\right)  $
into itself. Its adjoint on the dual space of measures $M\left(
\mathbb{T}\right)  $ is given by $\left(  R^{\ast}\mu\right)  \left(
\xi\right)  =\mu\left(  R\xi\right)  =\int_{\mathbb{T}}\left(  R\xi\right)
\left(  z\right)  \,d\mu\left(  z\right)  $, $\xi\in C\left(  \mathbb{T}%
\right)  $. The Dirac point-measure $\delta_{1}\in M\left(  \mathbb{T}\right)
$, given by $\delta_{1}\left(  \xi\right)  =\xi\left(  1\right)  $, is
invariant by (\ref{eq5}), i.e., $R^{\ast}\left(  \delta_{1}\right)
=\delta_{1}$.

Consider the eigenvalue problem:%
\begin{equation}
R\xi_{0}=\lambda\xi_{0},\qquad\lambda\in\mathbb{C},\;\xi_{0}\in P\left[
-N,N\right]  \setminus\left\{  0\right\}  . \label{eqGeneral.REigen}%
\end{equation}
Then $\delta_{1}\left(  \xi_{0}\right)  =\delta_{1}\left(  R\xi_{0}\right)
=\lambda\delta_{1}\left(  \xi_{0}\right)  $, so $\delta_{1}\left(  \xi
_{0}\right)  =\xi_{0}\left(  1\right)  =0$ if $\lambda\neq1$. We assume this,
and since $\left\|  R\right\|  _{\infty\rightarrow\infty}=1$, the discussion
may be restricted to $\left|  \lambda\right|  =1$. We claim that, if
$\lambda\neq1$, $\lambda\in\mathbb{T}$, then $\xi_{0}=0$, so we cannot have
nontrivial peripheral spectrum.

By (\ref{ThmGeneral.2(2)}), $\left\|  \varphi-M^{n}\psi^{\left(  0\right)
}\right\|  _{2}\rightarrow0$, where $\psi^{\left(  0\right)  }$ is any initial
vector with the stated conditions, e.g., $\psi^{\left(  0\right)  }%
=\chi_{\left[  0,1\right]  }$. Using Lemma \ref{LemPropQuad}%
(\ref{LemPropQuad(4)}), we also get%
\begin{equation}
\left\|  \varphi-M^{n}\left(  \xi\ast\psi^{\left(  0\right)  }\right)
\right\|  _{L^{2}\left(  \mathbb{R}\right)  }\longrightarrow0,
\label{eqOla2.87}%
\end{equation}
whenever $\xi\left(  1\right)  =1$. For example, take $\xi=\openone+c\xi_{0}$
to have this satisfied. Then \linebreak $R^{n}\xi=\openone+c\lambda^{n}\xi
_{0}$, i.e., a divergent sequence if $\lambda\neq1$ and $c\neq0$, supposing
$\xi_{0}\neq0$.

Using Lemma \ref{LemPropQuad}, we will show that%
\begin{align}
R^{n}\xi &  =p\left(  M^{n}\psi^{\left(  0\right)  },M^{n}\left(  \xi\ast
\psi^{\left(  0\right)  }\right)  \right) \label{eqOla2.88}\\
&  \underset{n\rightarrow\infty}{\longrightarrow}p\left(  \varphi
,\varphi\right)  ,\nonumber
\end{align}
where the last convergence is in the finite-dimensional subspace of
$\mathbb{C}\left[  z,z^{-1}\right]  $ and thus in any norm. This will
contradict the divergence of $R^{n}\xi$. The formula (\ref{eqOla2.88}) can be
verified in two ways: since $\varphi$, $M^{n}\left(  \psi^{\left(  0\right)
}\right)  $ and $M^{n}\left(  \xi\ast\psi^{\left(  0\right)  }\right)  $ all
have support inside a common compact set, the convergence is immediate from
the finite sum (\ref{eq44}), (\ref{eqOla2.85}), and (\ref{eqOla2.87}).
Alternatively one can use Lemma \ref{LemPropQuad}(\ref{LemPropQuad(1)}), and
$\operatorname*{PER}\left|  \widehat{\psi^{\left(  0\right)  }}\right|
^{2}=\openone$. In checking the conditions in Lemma \ref{LemPropQuad}%
(\ref{LemPropQuad(1)}), we note that $M^{n}\left(  \xi\ast\psi^{\left(
0\right)  }\right)  \underset{n\rightarrow\infty}{\longrightarrow}\xi\left(
1\right)  \varphi$, in $L^{2}\left(  \mathbb{R}\right)  $, so we must verify
that%
\[
\left\|  \operatorname*{PER}\left(  \left|  \widehat{M^{n}\psi^{\left(
0\right)  }}\right|  ^{2}\right)  \right\|  _{\infty}%
\]
is bounded in $n$. But the function inside $\left\|  \,\cdot\,\right\|
_{\infty}$ equals%
\[
p\left(  M^{n}\psi^{\left(  0\right)  },M^{n}\psi^{\left(  0\right)  }\right)
=R^{n}\left(  p\left(  \psi^{\left(  0\right)  },\psi^{\left(  0\right)
}\right)  \right)  =R^{n}\openone=\openone.
\]
This contradiction completes the first part of the proof of
(\ref{ThmGeneral.2(2)}) $\Rightarrow$ (\ref{ThmGeneral.2(1)}).

It remains to show that (\ref{ThmGeneral.2(2)}) of Theorem \ref{ThmGeneral.2}
implies that $\lambda=1$ has multiplicity one in the spectrum of $R$, where
again $R$ is identified with its restriction to $P\left[  -N,N\right]  $.
Since $R\openone=\openone$, we need only exclude that the multiplicity is $2$
or more. But $\left\|  M^{n}\psi^{\left(  0\right)  }\right\|  _{2}=1$ for all
$n$ by (\ref{eq45}) and $R\openone=\openone$ and it follows from (\ref{eq54})
that $\left\|  \varphi\right\|  _{2}=1$. By (\ref{eq18}) $\Rightarrow$
(\ref{eq21}) it follows that $\lambda=1$ has multiplicity $1$.
\end{proof}

\section{\label{Some}Some examples}

Our interest in the subject of cascade approximation was ignited when using
the cascade algorithm in \cite{BEJ99} to draw the scaling function associated
with the low-pass wavelet filter%
\begin{equation}
m_{0}^{\left(  \theta\right)  }\left(  z\right)  =\sum_{k=0}^{3}a_{k}^{\left(
\theta\right)  }z^{k} \label{eqSome.1}%
\end{equation}
where%
\begin{equation}%
\begin{aligned} a_{0}^{\left
( \theta\right) } & =\frac{1}{2\sqrt{2}}\left( 1-\cos\theta+\sin\theta
\right) ,\\ a_{1}^{\left( \theta\right) }& =\frac{1}{2\sqrt{2}}\left
( 1-\cos\theta-\sin\theta\right) ,\\ a_{2}^{\left( \theta\right) } & =\frac
{1}{2\sqrt{2}}\left( 1+\cos\theta-\sin\theta\right) ,\\ a_{3}^{\left
( \theta\right) } & =\frac{1}{2\sqrt{2}}\left( 1+\cos\theta+\sin\theta
\right) ; \end{aligned}%
\label{eqSome.2}%
\end{equation}
and $\theta$ varies over the circle; see \cite{Pol89,Pol90}. This family is
discussed in detail in \cite{BEJ99}. The symmetry $\varphi^{\left(  \pi
-\theta\right)  }\left(  x\right)  =\varphi^{\left(  \theta\right)  }\left(
3-x\right)  $ (see \cite[Proposition 4.1]{BEJ99}) implies that we may limit
the discussion to the interval $-\frac{\pi}{2}\leq\theta\leq\frac{\pi}{2}$. In
this interval the conditions (\ref{eq17})--(\ref{eq22}) are fulfilled with the
one exception of $\theta=\frac{\pi}{2}$, where%
\[
m_{0}^{\left(  \frac{\pi}{2}\right)  }\left(  z\right)  =\frac{1}{\sqrt{2}%
}\left(  1+z^{3}\right)
\]
and the scaling function $\varphi$ is given by%
\[
\varphi^{\left(  \frac{\pi}{2}\right)  }\left(  x\right)  =%
\begin{cases}
\frac{1}{3} & \text{for }0\leq x\leq3,  \\
0 & \text{otherwise.}
\end{cases}%
\]

It has been observed by several authors that at this point the cascade
approximants converge merely weakly; see \cite[Figure 6]{BEJ99}, \cite[Figure
3.3]{CoRy95}, \cite[Note 4]{Str96}, \cite[Figure 4]{Coh92}. When using the
cascade algorithm to depict $\varphi^{\left(  \theta\right)  }$ for $\theta$
near $\frac{\pi}{2}$, this weak convergence seems to persist; see Figures
\ref{Res045_00}--\ref{Res045_stages}. So let us compute the Ruelle operator in
(\ref{eq25}) $R$ as a $7\times7$ matrix on $P\left[  -3,3\right]  $. The
result is the following slant-Toeplitz matrix:%
\[
R=%
\begin{bmatrix}
b & 0 & 0 & 0 & 0 & 0 & 0\\
c & 0 & b & 0 & 0 & 0 & 0\\
c & 1 & c & 0 & b & 0 & 0\\
b & 0 & c & 1 & c & 0 & b\\
0 & 0 & b & 0 & c & 1 & c\\
0 & 0 & 0 & 0 & b & 0 & c\\
0 & 0 & 0 & 0 & 0 & 0 & b
\end{bmatrix}
\]
where $b=a_{3}a_{0}$ and $c=a_{1}a_{0}+a_{2}a_{1}+a_{3}a_{2}$, and we have
used (\ref{eq2}). Viewing $a=\left(  a_{0},a_{1},a_{2},a_{3}\right)  $ as a
function on $\mathbb{Z}_{4}$, and letting $T$ denote cyclic translation on
$\mathbb{Z}_{4}$,
\[
T\left(  a_{0},a_{1},a_{2},a_{3}\right)  =\left(  a_{1},a_{2},a_{3}%
,a_{0}\right)  ,
\]
the relations (\ref{eq2}) take the form%
\[
\left\langle a,a\right\rangle  =1,\qquad
\left\langle a,T^{2}a\right\rangle  =0.
\]
Furthermore%
\[
\left\langle a,Ta\right\rangle =\left\langle a,T^{3}a\right\rangle =b+c,
\]
and by (\ref{eq3}),%
\[
\sum_{k=0}^{3}T^{k}a=\left(  \sum a_{i},\sum a_{i},\sum a_{i},\sum
a_{i}\right)  =\sqrt{2}\left(  1,1,1,1\right)  .
\]
Thus%
\[
\left\langle a,\sum_{k=0}^{3}T^{k}a\right\rangle =\left\langle \left(
a_{0},a_{1},a_{2},a_{3}\right)  ,\sqrt{2}\left(  1,1,1,1\right)  \right\rangle
=\sqrt{2}\sum_{i}a_{i}=2.
\]
On the other side,%
\[
\left\langle a,\sum_{k=0}^{3}T^{k}a\right\rangle =\sum_{k=0}^{3}\left\langle
a,T^{k}a\right\rangle =1+\left(  b+c\right)  +0+\left(  b+c\right)
=1+2\left(  b+c\right)  ,
\]
so%
\[
b+c={1}/{2}.
\]

\begin{figure}[ptb]
\includegraphics
[bb=0 55 288 243,width=360pt]{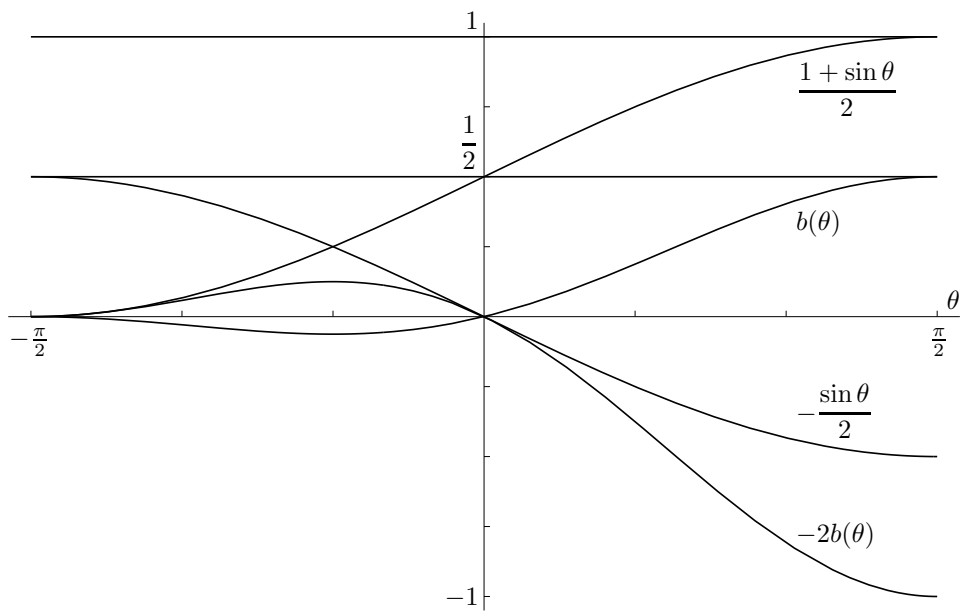} \raisebox{8pt}[0pt][0pt]{\makebox
[360pt]{\hskip8pt\raisebox{104pt}[0pt][0pt]{\makebox[0pt]{\hss$-\frac{\pi}%
{2}$\hss}}\hskip170pt\raisebox
{6pt}[0pt][0pt]{\llap{$-1$}}\raisebox{178pt}[0pt][0pt]{\llap{$\displaystyle
\frac{1}{2}$}}\raisebox{224pt}[0pt][0pt]{\llap{$1$}}\hskip120pt\raisebox
{30pt}[0pt][0pt]{\rlap{$-2b(\theta
)$}}\raisebox{76pt}[0pt][0pt]{\rlap{$\displaystyle-\frac{\sin\theta
}{2}$}}\raisebox{148pt}[0pt][0pt]{\rlap{$b(\theta
)$}}\raisebox{198pt}[0pt][0pt]{\rlap{$\displaystyle\frac{1+\sin\theta
}{2}$}}\hskip54pt\raisebox{104pt}[0pt][0pt]{\makebox[0pt]{\hss$\frac{\pi
}{2}$\hss}}\hskip8pt\raisebox{118pt}[0pt][0pt]{\llap{$\theta
$}}}}\caption{Eigenvalues $\lambda\left(  \theta\right)  $ of $R$}%
\label{boftheta}%
\end{figure}

\noindent
Since $b\left(  \theta\right)  =\frac{1}{8}\left[  1+2\sin\theta-\cos\left(
2\theta\right)  \right]  $, $b\left(  \theta\right)  $ ranges over $\left[
-\frac{1}{16},\frac{1}{2}\right]  $ as $\theta$ ranges over $\left[
-\frac{\pi}{2},\frac{\pi}{2}\right]  $. The characteristic equation of $R$ has
the obvious roots $1$ and $b$, where $b$ has multiplicity $2$. Using the flip
symmetry of the matrix under the joint flip around the fourth column and
fourth row, the remaining fourth-order polynomial factors into two
second-order polynomials, and thus one deduces that the eigenvalues of $R$ are%
\newcommand
{\eigenstrut}
{\vphantom{$\displaystyle\frac{1-\sqrt{1+16b}}{4}=\left\{
\begin{array}
[c]{ll}\frac{1+\sin\theta}{2_{\mathstrut}} & \left(  -\frac{\pi}{2}\leq\theta
\leq-\frac{\pi}{6}\right)  \\
-\frac{\sin\theta}{2} & \left(  -\frac{\pi}{6}\leq\theta\leq\frac{\pi}%
{2}\right)
\end{array}
\right.  $}}%
\[%
\begin{tabular}
[c]{c|c}%
Eigenvalue $\lambda\left(  \theta\right)  $ & Range of $\lambda\left(
\theta\right)  $\\\hline
\eigenstrut$1$ & $\left\{  1\right\}  $\\
\eigenstrut$b$ (multiplicity two) & $\left[  -\frac{1}{16},\frac{1}{2}\right]
$\\
\eigenstrut$\displaystyle\frac{1}{2}$ & $\left\{  \frac{1}{2}\right\}  $\\
\eigenstrut$-2b$ & $\left[  -1,\frac{1}{8}\right]  $\\
\eigenstrut$\displaystyle\frac{1+\sqrt{1+16b}}{4}=\left\{
\begin{array}
[c]{ll}%
-\frac{\sin\theta}{2_{\mathstrut}} & \left(  -\frac{\pi}{2}\leq\theta
\leq-\frac{\pi}{6}\right) \\
\frac{1+\sin\theta}{2_{\mathstrut}} & \left(  -\frac{\pi}{6}\leq\theta
\leq\frac{\pi}{2}\right)
\end{array}
\right.  $ & $\left[  \frac{1}{4},1\right]  $\\
\eigenstrut$\displaystyle\frac{1-\sqrt{1+16b}}{4}=\left\{
\begin{array}
[c]{ll}%
\frac{1+\sin\theta}{2_{\mathstrut}} & \left(  -\frac{\pi}{2}\leq\theta
\leq-\frac{\pi}{6}\right) \\
-\frac{\sin\theta}{2} & \left(  -\frac{\pi}{6}\leq\theta\leq\frac{\pi}%
{2}\right)
\end{array}
\right.  $ & $\left[  -\frac{1}{2},\frac{1}{4}\right]  $%
\end{tabular}
\]
We see that unless $b=\frac{1}{2}$, i.e., $\theta=\frac{\pi}{2}$, $1$ is the
unique peripheral eigenvalue, and it has multiplicity $1$. By Theorem
\ref{ThmGeneral.2}, the cascade approximants then do indeed converge to the
scaling function. If $\theta=\frac{\pi}{2}$, then $1$ is an eigenvalue of
multiplicity $2$ and $-1$ is an eigenvalue of multiplicity $1$, and by Theorem
\ref{ThmGeneral.2} the cascades do not converge in $L^{2}$-norm. The slow
convergence near $\theta=\frac{\pi}{2}$ can be explained by the fact that $R$
then has eigenvalues differing from $\pm1$ by
$\smash[b]{\mathrm{O}
\left(  \left(  \theta -\frac{\pi}{2}\right)  ^{2}\right)  }$,
so the approximation after $n$ steps is
like
$\left(  
1-\text{const.}\left(  \theta-\frac{\pi}{2}\right)  ^{2}\right)
^{n}$. This slow rate of convergence near $\frac{\pi}{2}$ is also clear from
Figure \ref{Res045_stages}(a)--(l). We should emphasize that the semiregular
layers displayed by these pictures dissolve more and more when doing further
iterations and a plot after $1000$ iterations with a resolution of $2^{-10}$
shows virtually no discernible small-scale regularity. This is shown in
Figures \ref{MPsiPlus}--\ref{MPsiDiff}, and let us explain how these were
produced. First note that by (\ref{eq9}), the value of $\left(  M\psi\right)
\left(  x\right)  $ at a point 
$x\in2^{-N}\mathbb{Z}$ 
only depends on the values of
$\psi\left(y\right)$ at points $y\in 2^{-(N-1)}\mathbb{Z}$,
and hence the exact values $\psi^{(n)}$ for
$x\in2^{-N}\mathbb{Z}$ 
can be determined
exactly by an iterated matrix scheme. Also the jumps of $\psi^{\left(
n\right)  }$ at the points in $2^{-N}\mathbb{Z}$ can be determined exactly by
the following cascade scheme with fixed $N$ (e.g., $N=10$): We define%
\[
x_{n}=n\cdot2^{-N},\qquad n=0,1,\dots,3\cdot2^{N},
\]
and define
\[
\psi_{+}^{\left(  m\right)  }\left(  n\right)   =\lim_{x\rightarrow
n\cdot2^{-N}+}\psi^{\left(  m\right)  }\left(  x\right)  ,\qquad
\psi_{-}^{\left(  m\right)  }\left(  n\right)   =\lim_{x\rightarrow
n\cdot2^{-N}-}\psi^{\left(  m\right)  }\left(  x\right)  ,
\]
where%
\[
\psi^{\left(  0\right)  }\left(  x\right)  =%
\begin{cases}
1 & 0\leq x<1,  \\
0 & 1\leq x\leq3,
\end{cases}
\]
and%
\[
\psi^{\left(  m\right)  }\left(  x\right)  =\left(  M^{m}\psi^{\left(
0\right)  }\right)  \left(  x\right)  =M\psi^{\left(  m-1\right)  }\left(
x\right)  ,
\]
and%
\begin{equation}
M\psi\left(  x\right)  =\sqrt{2}\sum_{k=0}^{3}a_{k}\psi\left(  2x-k\right)  .
\label{eqScheme.star}%
\end{equation}
We start the recursion with%
\[
\psi_{+}^{\left(  0\right)  }\left(  n\right)  =%
\begin{cases}
1 & \text{for }n=0,1,\dots,2^{N}-1,  \\
0 & \text{for }n=2^{N},2^{N}+1,\dots,3\cdot2^{N},
\end{cases}
\]
and%
\[
\psi_{-}^{\left(  0\right)  }\left(  n\right)  =%
\begin{cases}
0 & \text{for }n=0,  \\
1 & \text{for }n=1,2,\dots,2^{N},  \\
0 & \text{for }n=2^{N}+1,\dots,3\cdot2^{N},
\end{cases}
\]
and from (\ref{eqScheme.star}),%
\begin{equation}%
\begin{cases}
& \psi_{+}^{\left( m\right) }\left( n\right) =\sqrt{2}\sum_{k=0}^{3}%
a_{k}\psi_{+}^{\left( m-1\right) }\left( 2n-k\cdot2^{N}\right) , \\
& \psi_{-}^{\left( m\right) }\left( n\right) =\sqrt{2}\sum_{k=0}^{3}%
a_{k}\psi_{-}^{\left( m-1\right) }\left( 2n-k\cdot2^{N}\right) ,
\end{cases}
\label{eqScheme.starstar}%
\end{equation}
where we use the convention that
\[
\psi_{\pm}^{\left(  m-1\right)  }\left(  l\right)  =0
\]
if $l<0$ or $l>3\cdot2^{N}$.

In Figures \ref{MPsiPlus}--\ref{MPsiDiff}, we have used this algorithm to plot
$\psi_{+}^{\left(  1000\right)  }\left(  n\right)  $, $\psi
_{-}^{\left(  1000\right)  }\left(  n\right)  $, $\left(  \psi
_{+}^{\left(  1000\right)  }-\psi_{+}^{\left(  1000\right)  }\right)
\left(  n\right)  $ for $N=10$, $\theta=9\pi/20$. If $1000$ is replaced by
$100$ the large-scale plots of $\psi_{-}^{\left(  100\right)  }$ and
$\psi_{+}^{\left(  100\right)  }$ look similar, but
$\smash{\left\|  \psi
_{+}^{\left(  100\right)  }-\psi_{+}^{\left(  100\right)  }\right\|
_{\infty}}$ is much larger than $\left\|  \psi_{+}^{\left(  1000\right)
}-\psi_{+}^{\left(  1000\right)  }\right\|  _{\infty}<5\cdot10^{-6}$. Thus
the plot shows an amazing amount of local continuity even though the
larger-scale behaviour is quite irregular. However, the scaling function at
$\theta=9\pi/20$ is indeed discontinuous by the discussion below.
This apparent non-compatibility is explained in detail
in the Appendix.

\begin{figure}[ptb]
\begin{picture}(358,267)(0,-24)
\put(0,0){\includegraphics
[bb=0 55 288 233,width=354bp]{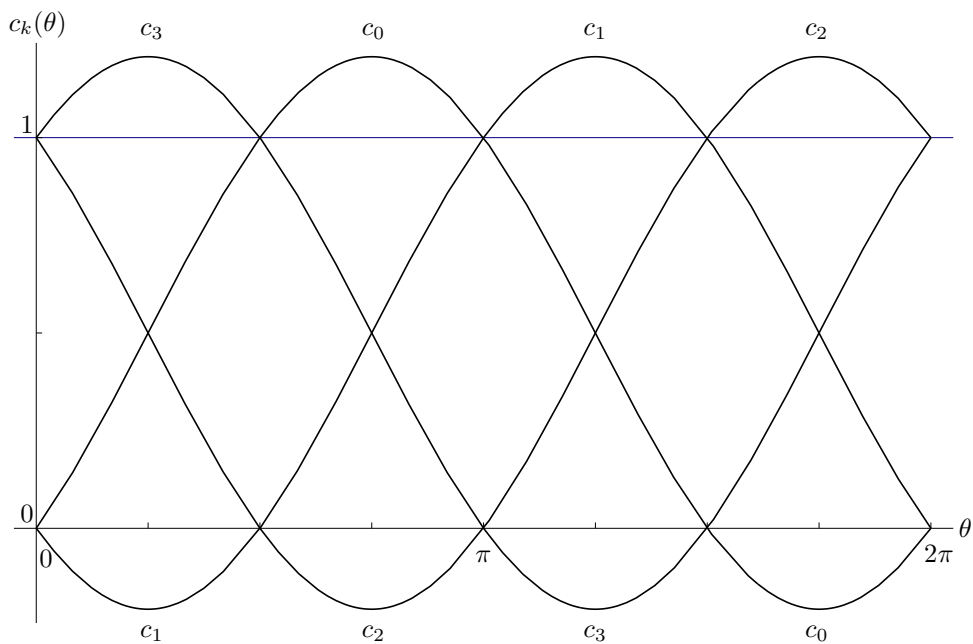}}
\put(9,226){\makebox(0,0){$c_{k}(\theta)$}}
\put(52,223){\makebox(0,0){$c_{3}$}}
\put(135.5,223){\makebox(0,0){$c_{0}$}}
\put(219,223){\makebox(0,0){$c_{1}$}}
\put(302.5,223){\makebox(0,0){$c_{2}$}}
\put(52,-4){\makebox(0,0){$c_{1}$}}
\put(135.5,-4){\makebox(0,0){$c_{2}$}}
\put(219,-4){\makebox(0,0){$c_{3}$}}
\put(302.5,-4){\makebox(0,0){$c_{0}$}}
\put(7.5,185){\makebox(0,0)[br]{$1$}}
\put(7.5,38){\makebox(0,0)[br]{$0$}}
\put(9.5,21){\makebox(0,0)[bl]{$0$}}
\put(177,23){\makebox(0,0)[b]{$\pi$}}
\put(343,23){\makebox(0,0)[bl]{$2\pi$}}
\put(356,36){\makebox(0,0)[l]{$\theta$}}
\end{picture}
\caption{$c_{k}\left( \theta \right) =\sqrt{2}a_{k}^{(\theta)}$,
$k=0,1,2,3$. See (\ref{eqSomeNew.7}).}%
\label{ckplot}%
\end{figure}

\begin{remark}
\label{RemNeighborhood}In \cite{BEJ99}, we noted that there are two
neighborhoods on the circle, one near $\theta_{1}=\frac{7\pi}{6}$, and a
symmetric one near $\theta_{2}=-\frac{\pi}{6}$, such that each of the scaling
functions $\varphi^{\left(  \theta\right)  }$ has $x\mapsto\varphi^{\left(
\theta\right)  }\left(  x\right)  $ continuous when $\theta$ is in the union
of the two neighborhoods. On the other hand if $0<\theta<\frac{\pi}{2}$, then
$\sqrt{2}a_{3}^{\left(  \theta\right)  }>1$ so at the right-hand endpoint $x$
of each of the dyadic partitions, we will have $\left(  M^{n}\psi^{\left(
0\right)  }\right)  \left(  x\right)  \underset{n\rightarrow\infty
}{\longrightarrow}$ $\infty$. It then follows from \cite[Proposition 6.5.2 and
footnote 9]{Dau92} or \cite{DaLa92} that $x\mapsto\varphi^{\left(
\theta\right)  }\left(  x\right)  $ cannot be continuous for $\theta$ in the
first quarter circle. By the symmetry $\varphi^{\left(  \pi-\theta\right)
}\left(  x\right)  =\varphi^{\left(  \theta\right)  }\left(  3-x\right)  $, it
is then also not continuous in the second quarter of the $\theta$-circle.

The question of continuity of $\varphi^{\left(  \theta\right)  }$ has been
considered in even more detail for our examples in the papers \cite{CoHe92},
\cite[Section 4.8.1]{CoHe94}, \cite{Wan95}, \cite{Wan96}. They use the real
coefficients%
\begin{equation}
c_{n}=\sqrt{2}a_{n}, \label{eqSomeNew.3}%
\end{equation}
which then satisfy%
\begin{equation}%
\begin{aligned}
\sum_{k}c_{k}c_{k+2l} &=2\delta_{l},  \\
\sum_{k}c_{k} &=2,
\end{aligned}
\label{eqSomeNew.4}%
\end{equation}
or, equivalently,%
\begin{equation}
c_{0}+c_{2}=1=c_{1}+c_{3}, \label{eqSomeNew.5}%
\end{equation}
so the scaling operator $M$ in \textup{(\ref{eq9})} becomes%
\begin{equation}
\left(  M\psi\right)  \left(  x\right)  =\sum_{k=0}^{3}c_{k}\psi\left(
2x-k\right)  , \label{eqSomeNew.6}%
\end{equation}
where%
\begin{equation}%
\begin{aligned}
c_{0} & =\frac{1}{2}\left( 1-\cos\theta+\sin\theta\right) ,\\
c_{1} & =\frac{1}{2}\left( 1-\cos\theta-\sin\theta\right) ,\\
c_{2} & =\frac{1}{2}\left( 1+\cos\theta-\sin\theta\right) ,\\
c_{3} & =\frac{1}{2}\left( 1+\cos\theta+\sin\theta\right) ;
\end{aligned}%
\label{eqSomeNew.7}%
\end{equation}
see an illustration in Figure \textup{\ref{ckplot}.} In \cite[Proposition
2.1]{Wan96}, it is stated that if there is a continuous scaling function,
then $\left|  c_{0}\right|  <1$ and $\left|  c_{3}\right|  <1$, that is,
$\pi<\theta<2\pi$. Since our movie reel shows that the scaling function
clearly is discontinuous for $0\leq\theta\leq\frac{\pi}{2}$ and then by
symmetry for $0\leq\theta\leq\pi$, this is consistent with the movie reel.
(It is still an open question, for $0<\theta <\pi $,
as to ``how discontinuous''
$x\mapsto \varphi^{\left( \theta \right) }\left( x\right) $
then is. Based on
graphics, and
analogies (see Section \ref{Con})
to iterated function
systems, it is likely
that the discontinuous
cases have interesting
fractal structure,
but that will be
postponed to a later
paper.) On
the other hand, the comment to Figure \textup{4.3} in \cite[page 193]{CoHe94}
indicates that the condition $\theta\in\left\langle \pi,3\pi/2\right\rangle
\cup\left\langle 3\pi/2,2\pi\right\rangle $ is necessary and sufficient for a
continuous scaling function.
This is indeed consistent with our movie reel, but the reel shows
extremely singular behaviour of the scaling function at
some dyadic rationals even in the domain of continuity.
The effect is most pronounced for Figures
\ref{Movie}(b) and \ref{Movie}(j).
We refer to the papers above and \cite{DaLa92} for the actual methods
used to establish continuity.
\end{remark}

The flip symmetry of the example above holds more generally whenever the
Ruelle matrix is defined from a low-pass filter function $m_{0}\left(
z\right)  =a_{0}+a_{1}z+\dots+a_{N}z^{N}$ with real coefficients $a_{i}$.
Specifically let $R$ denote the associated Ruelle operator. We then have

\begin{lemma}
\label{LemSome.1}%
\[
\left(  R\xi\right)  \left(  z^{-1}\right)  =\left(  R\xi\spcheck\right)
\left(  z\right)  ,
\]
where $\xi\spcheck\left(  z\right)  :=\xi\left(  z^{-1}\right)  $, $\xi
\in\mathbb{C}\left[  z,z^{-1}\right]  $, $z\in\mathbb{T}$.
\end{lemma}

\begin{proof}
{}From the definition of $R$ we have%
\begin{align*}
\left(  R\xi\right)  \left(  z^{-1}\right)   &  =\frac{1}{2}\sum_{w^{2}%
=z^{-1}}\left|  m_{0}\left(  w\right)  \right|  ^{2}\xi\left(  w\right) \\
&  =\frac{1}{2}\sum_{w^{2}=z}\left|  m_{0}\left(  w^{-1}\right)  \right|
^{2}\xi\left(  w^{-1}\right) \\
&  =\frac{1}{2}\sum_{w^{2}=z}\left|  \overline{m_{0}\left(  w\right)
}\right|  ^{2}\xi\spcheck\left(  w\right) \\
&  =R\left(  \xi\spcheck\right)  \left(  z\right)  ,
\end{align*}
where we used the reality assumption in the form $m_{0}\left(  w^{-1}\right)
=\overline{m_{0}\left(  w\right)  }$.
\end{proof}

We now turn to the graphics which illustrate the cascade approximation
(\ref{eq9})--(\ref{eq10}). It follows from (\ref{eq1}) and \cite[pp.~205--6]%
{Dau92} that this can be based on the matrix $\setlength{\arraycolsep}%
{0pt}\left( \begin{array}{cccccc}
a_0 \;& &\; 0 \;&\; 0 \;& &\; 0  \\ \cline{2-5}
a_2 \;& \vline&\; a_1 \;&\; a_0 \;& \vline&\; 0  \\
0 \;& \vline&\; a_3 \;&\; a_2 \;& \vline&\; a_1  \\ \cline{2-5}
0 \;& &\; 0 \;&\; 0 \;& &\; a_3
\end{array}\right) $
, or alternatively on just the framed $2\times2$ matrix $
\begin{pmatrix}
a_{1} & a_{0}\\
a_{3} & a_{2}%
\end{pmatrix}
$. The so-called cascade algorithm \cite[p.~205]{Dau92} starts with
initial points on an integral grid, and the $n$'th step fills in points on
$2^{-n}\mathbb{Z}$ places according to matrix multiplication and use of the
weights from the respective matrix entries. We include more technical points
in the captions of the pictures, and in the Appendix. 

\begin{remark}
\label{RemOrth}
While the results above concern primarily orthogonality and $L^{2}$-cascade
approximation, there is a direct connection between the $L^{2}$-theory and
pointwise features of the approximation, as is pointed out in \cite[p.~204]%
{Dau92}. Proposition 6.5.2 in \cite{Dau92} makes that explicit when the
\emph{a~priori} assumption is made that $\varphi^{\left(  \theta\right)
}\left(  x\right)  $ is continuous in $x$.
As noted in Remark \ref{RemNeighborhood} and Section \ref{Con}, and
in \cite{BEJ99},
such continuity is only known when $\theta$ is restricted to
certain subintervals of $\left(  -\pi,\pi\right]  $, and the pictures serve to
illustrate the features when $\theta$ is in the complement of the ``good'' regions.

We stress that the fast algorithm used for
some of
the
graphics (see the Appendix) does in fact
depend on the orthogonality of
the family
$\left\{ \varphi^{\left( \theta \right) }\left( x-n\right) 
\mid n\in \mathbb{Z}\right\} $ in
$L^{2}\left( \mathbb{R}\right) $,
and this orthogonality we verified
in \cite{BEJ99} to be satisfied
for all values of $\theta$
except for $\theta =
\frac{\pi}{2}$.
The significance of this orthogonality
is also directly related to the
assumption made in Theorem \ref{ThmGeneral.2}(\ref{ThmGeneral.2(2)}),
i.e., (\ref{eq24}), on the starting function
$\psi^{\left( 0\right) }$ for
the cascade approximation (\ref{eq12}), but
the fast algorithm of the Appendix is
different from (\ref{eq9})--(\ref{eq12}).
Details on the comparison of the two
are to be found in the Appendix.
\end{remark}

\section*{List of Figures}

\begin{figure}[ptb]
\begin{picture}(351,480)
\put(0,264){\includegraphics
[bb=99 0 333 432,height=216bp,width=117bp]{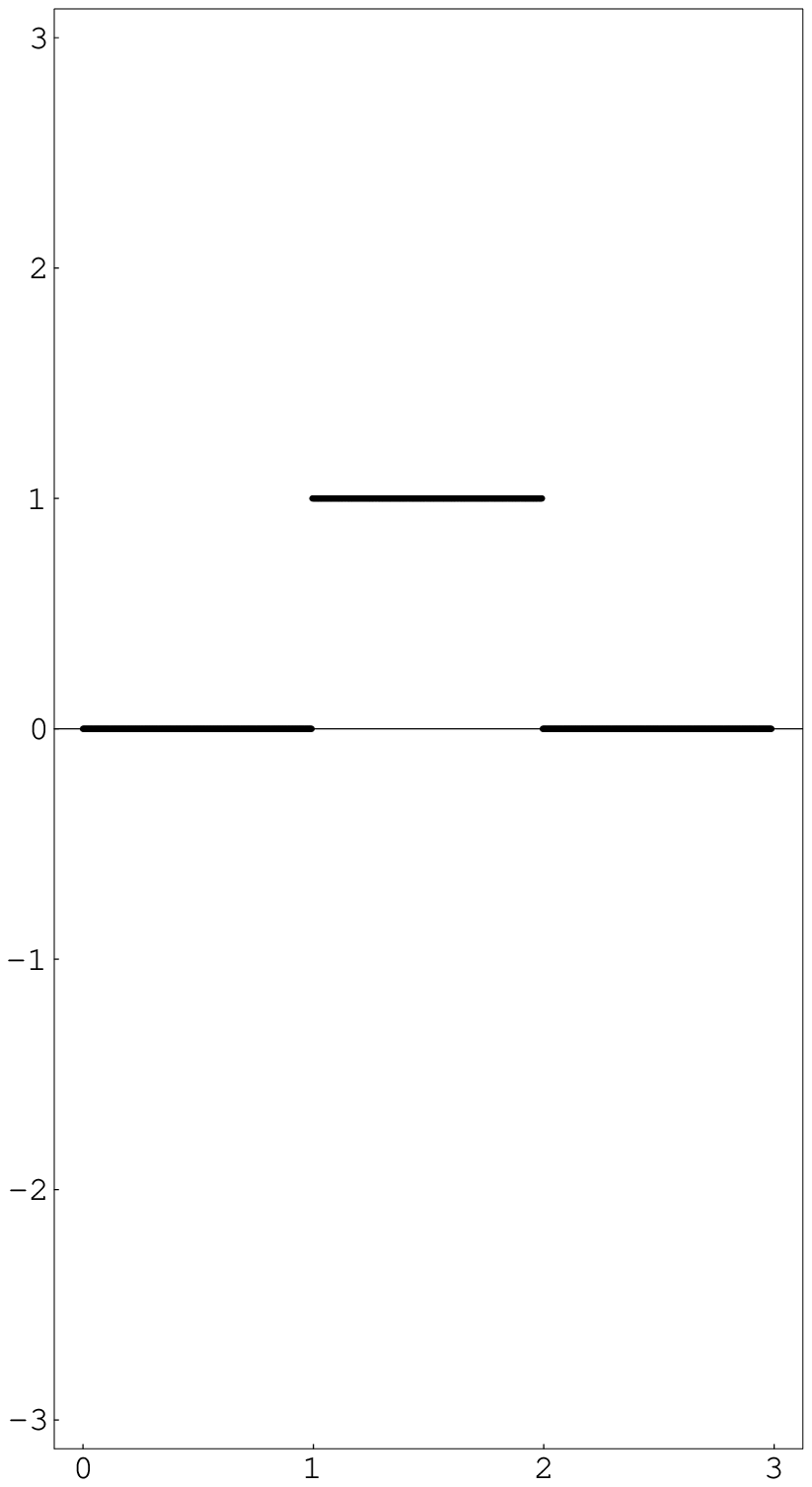}}
\put(117,264){\includegraphics
[bb=189 166 422 596,height=216bp,width=117bp]{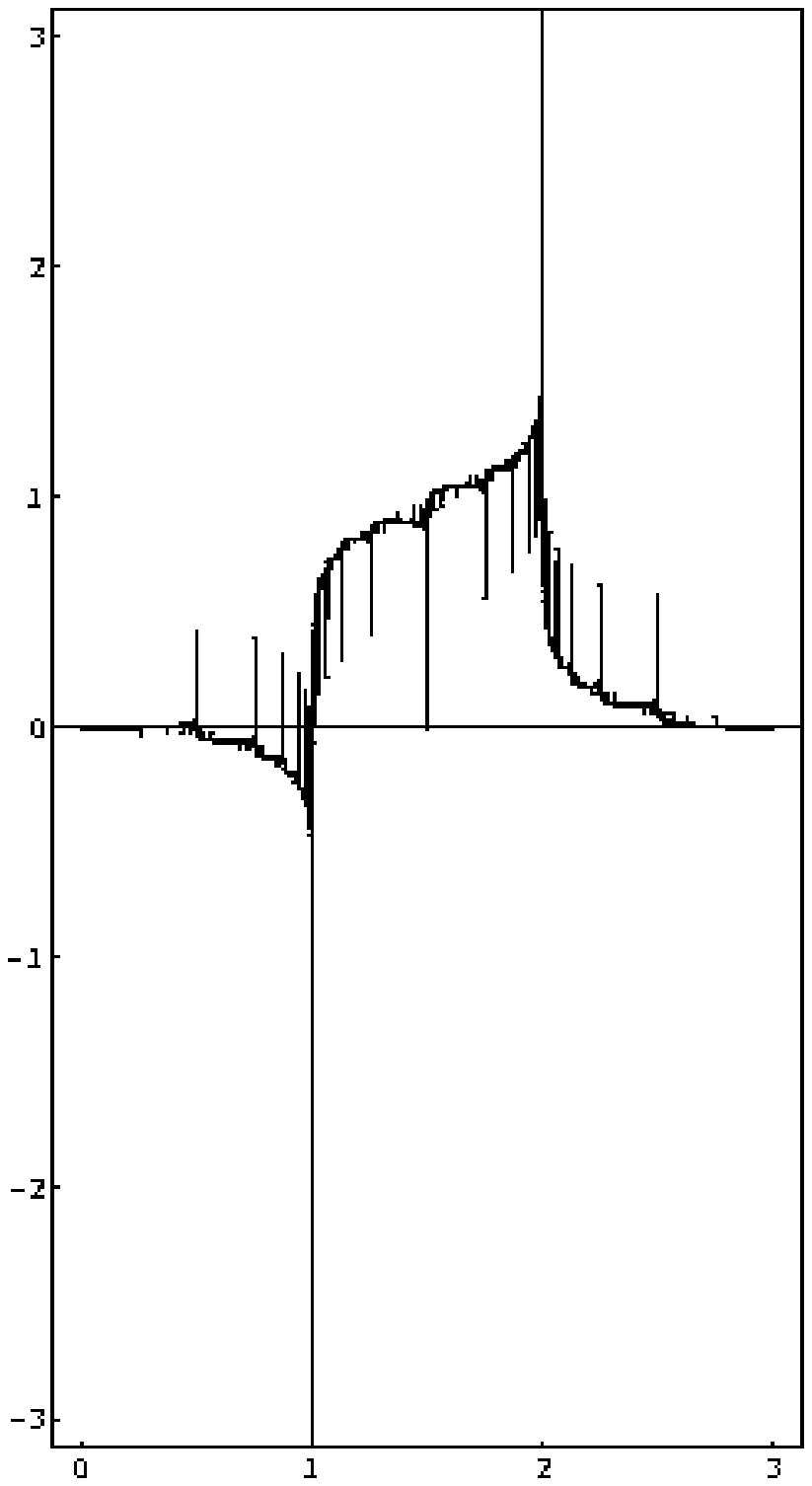}}
\put(234,264){\includegraphics
[bb=189 166 422 596,height=216bp,width=117bp]{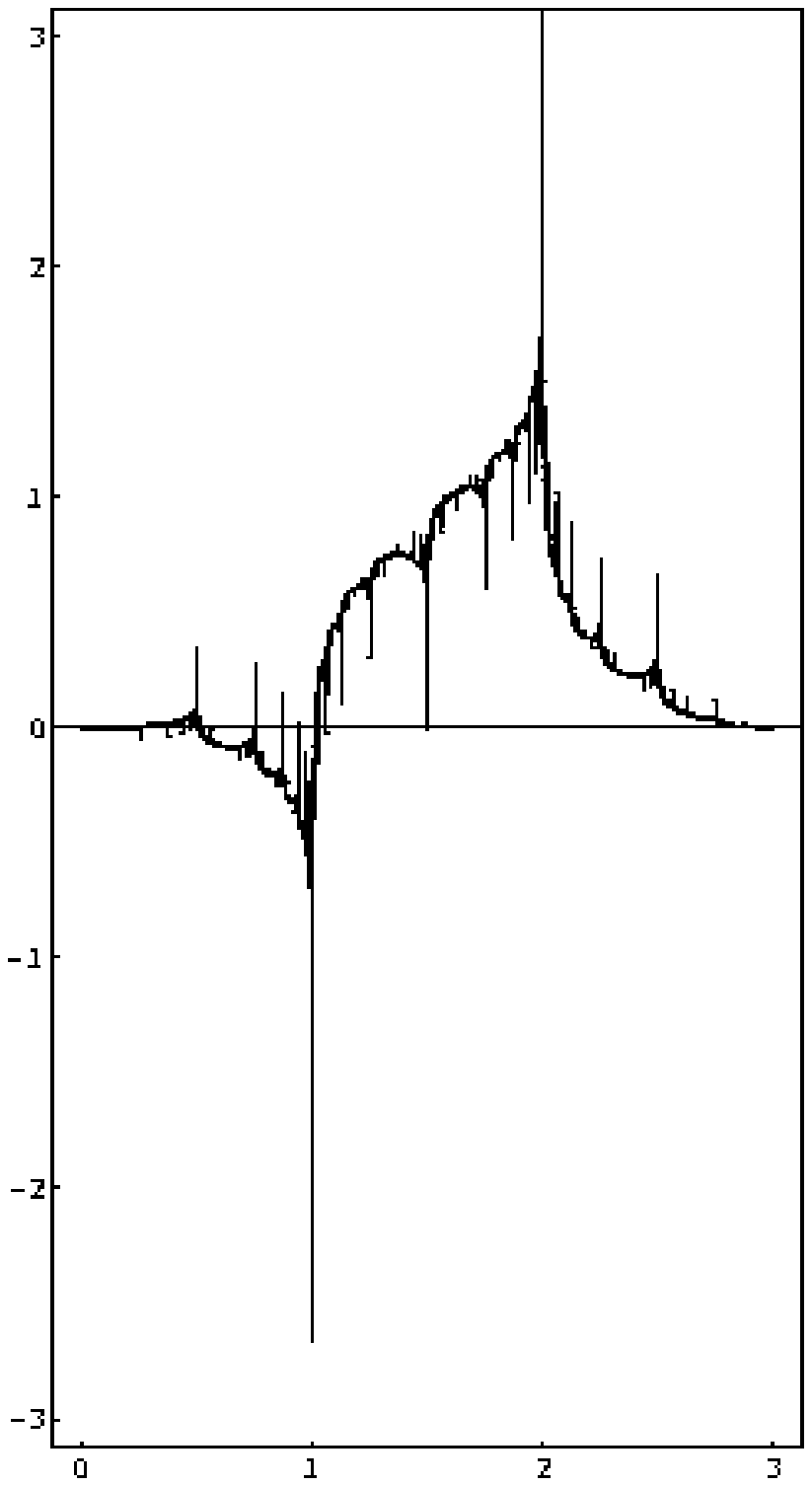}}
\put(0,249){\makebox(117,12){a: $\theta=-\frac{\pi}{2}$}}
\put(117,249){\makebox(117,12){b: $\theta=-\frac{9\pi}{20}$}}
\put(234,249){\makebox(117,12){c: $\theta=-\frac{2\pi}{5}$}}
\put(0,18){\includegraphics
[bb=189 166 422 596,height=216bp,width=117bp]{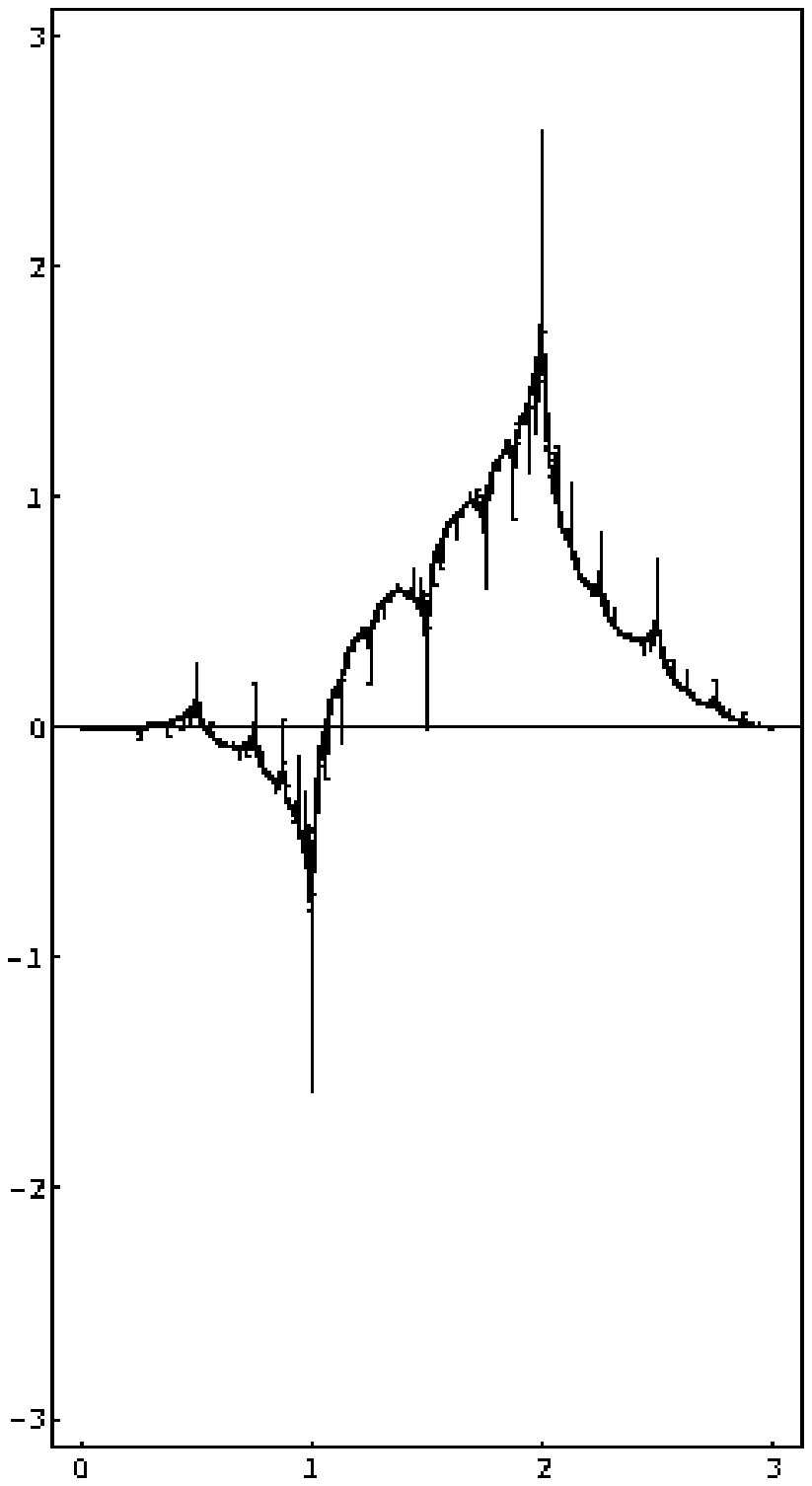}}
\put(117,18){\includegraphics
[bb=189 166 422 596,height=216bp,width=117bp]{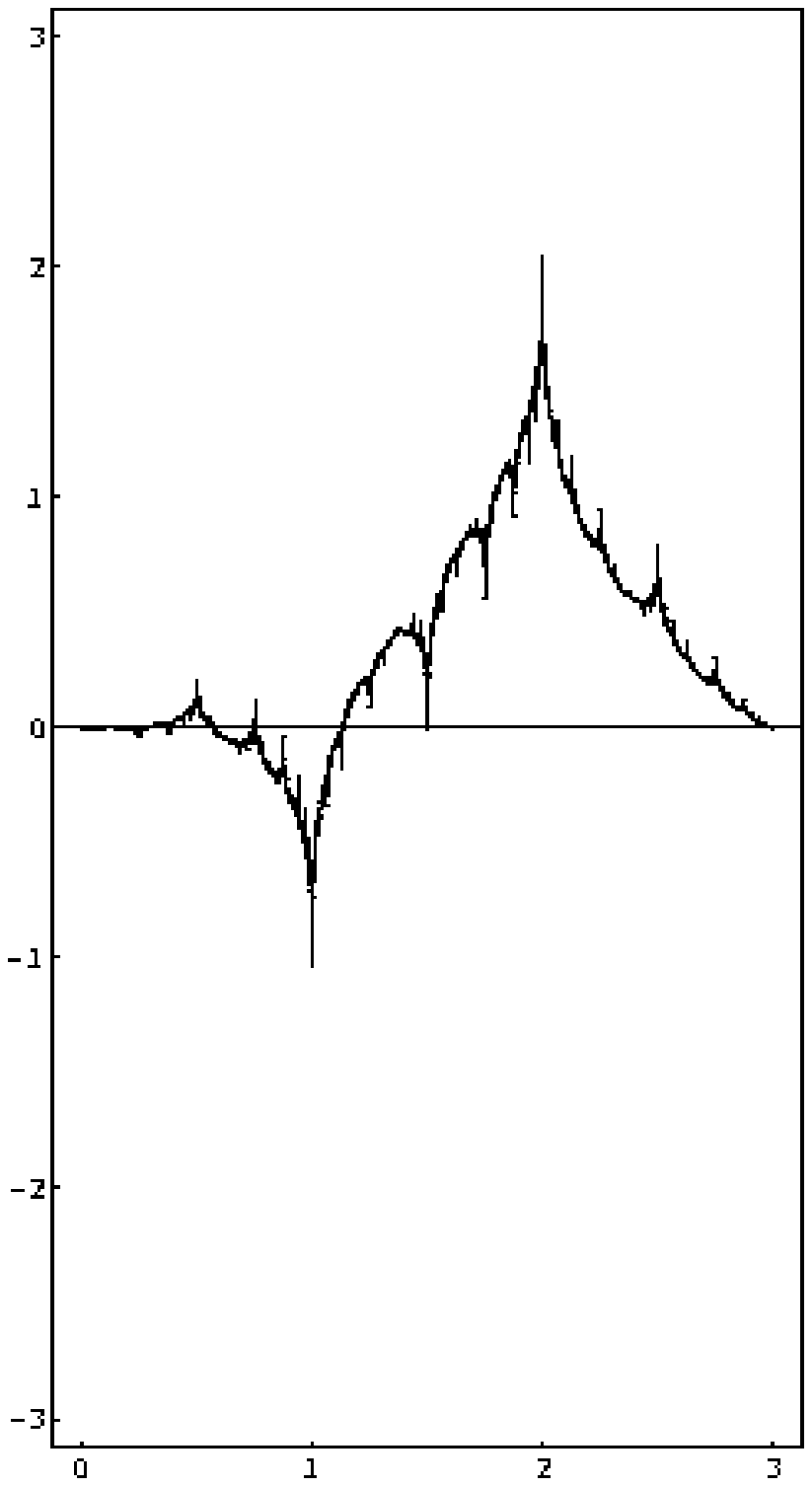}}
\put(234,18){\includegraphics
[bb=189 166 422 596,height=216bp,width=117bp]{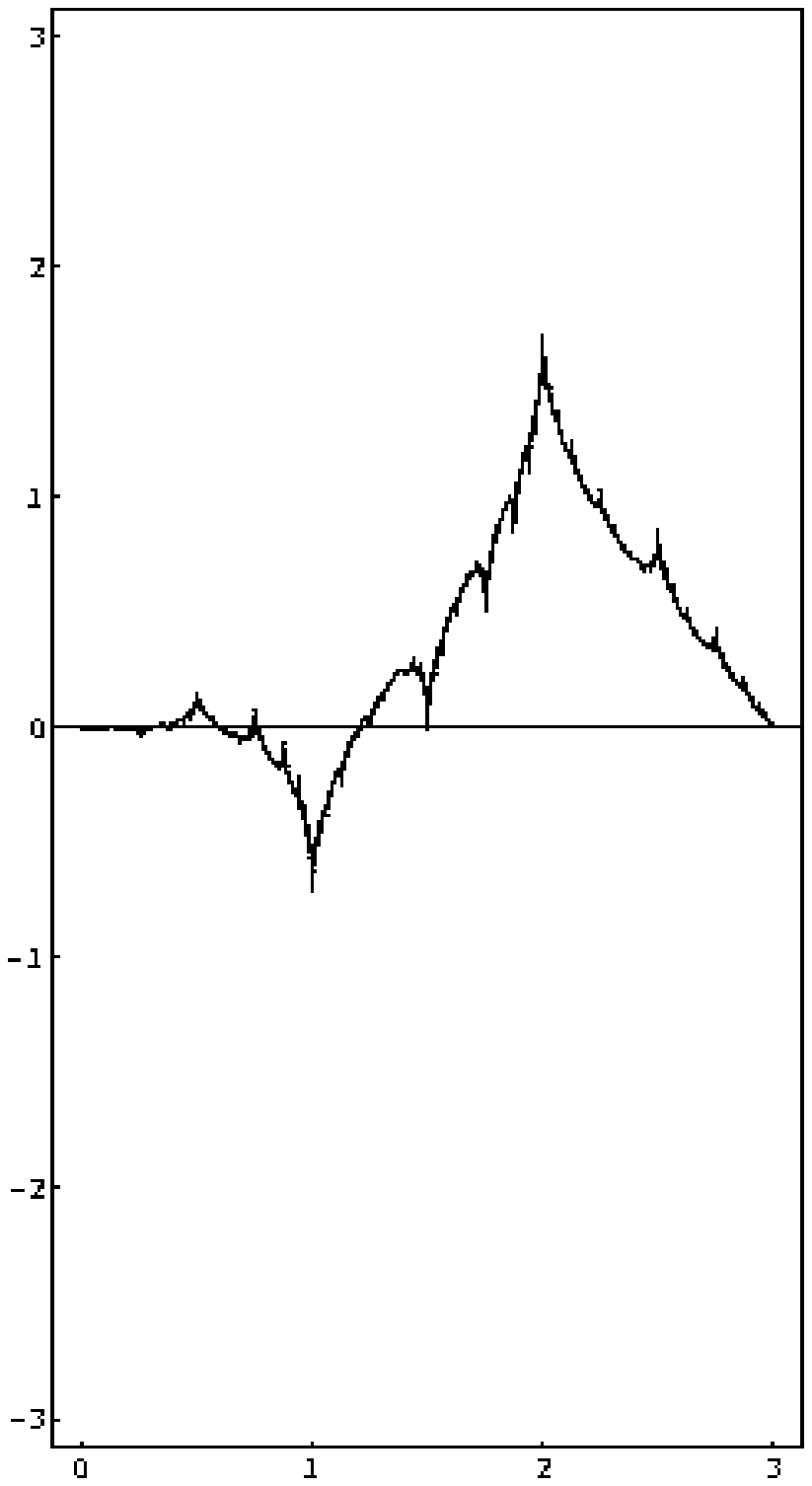}}
\put(0,3){\makebox(117,12){d: $\theta=-\frac{7\pi}{20}$}}
\put(117,3){\makebox(117,12){e: $\theta=-\frac{3\pi}{10}$}}
\put(234,3){\makebox(117,12){f: $\theta=-\frac{\pi}{4}$}}
\end{picture}
\caption{Scaling function ``movie reel'', $\theta$ from $\frac{-\pi}{2}$ to
$\frac{\pi}{2}$ in twenty-one frames: Frames a--f}%
\label{Movie}%
\end{figure}

\contentsline{figure}{\numberline{\protect\ref{boftheta}}{\ignorespaces
Eigenvalues $\lambda\left(\theta\right)$ of $R$}}{\protect\pageref{boftheta}}
\contentsline{figure}{\numberline{\protect\ref{ckplot}}{\ignorespaces
$c_{k}\left( \theta \right) =\sqrt{2}a_{k}^{(\theta)}$,
$k=0,1,2,3$
(see comments below)}}{\protect\pageref{ckplot}}
\contentsline{figure}{\numberline{\protect\ref{Movie}}{\ignorespaces
Scaling function ``movie reel'', $\theta$ from $\frac{-\pi}{2}$ to $\frac{\pi
}{2}$ in twenty-one frames
(see comments below)}}{\protect\pageref{Movie}--\protect\pageref
{Movie_end}}
\contentsline{figure}{\numberline{\protect\ref{Giraffe}%
}{\ignorespaces The $m\rightarrow \infty $ limit
of $\psi^{\left( m\right) }\left( x\right) $
for $\theta = -\frac{9\pi}{20}$}}{\protect\pageref
{Giraffe}}
\contentsline{figure}{\numberline{\protect\ref{Res045_00}%
}{\ignorespaces Initial function $\psi^{(0)}(x)=\chi_{[0,1]}%
(x)$ for the cascade series (Haar scaling function)}}{\protect\pageref
{Res045_00}} \contentsline{figure}{\numberline{\protect\ref{Res045_stages}%
}{\ignorespaces Cascade stages of scaling function, $\theta=\frac{9\pi}{20}$
(see comments below)}%
}{\protect\pageref{Res045_stages}--\protect\pageref{Res045_stages_end}}
\contentsline{figure}{\numberline{\protect\ref{MPsiPlus}}{\ignorespaces
$\psi_{+}^{(1000)}\left( n\right) $ at 
$n=0,1,\dots ,3\cdot 2^{10}$,
$x_{n}=n\cdot 2^{-10}$}}{\protect\pageref{MPsiPlus}}
\contentsline{figure}{\numberline{\protect\ref{MPsiMinus}}{\ignorespaces
$\psi_{-}^{(1000)}\left( n\right) $ at 
$n=0,1,\dots ,3\cdot 2^{10}$,
$x_{n}=n\cdot 2^{-10}$}}{\protect\pageref{MPsiMinus}}
\contentsline{figure}{\numberline{\protect\ref{MPsiDiff}}{\ignorespaces
$\left( 
\psi_{+}^{(1000)}\left( n\right) 
-\psi_{-}^{(1000)}\left( n\right) 
\right) $ at 
$n=0,1,\dots ,3\cdot 2^{10}$,
$x_{n}=n\cdot 2^{-10}$}}{\protect\pageref{MPsiDiff}}
\contentsline{figure}{\numberline{\protect\ref{InfinityPsi}}{\ignorespaces
$\psi_{\pm}^{(\infty)}\left( n\right) $ at 
$n=0,1,\dots ,3\cdot 2^{10}$,
$x_{n}=n\cdot 2^{-10}$}}{\protect\pageref{InfinityPsi}}
\bigskip

In Figure \ref{Res045_stages}(a)--(l), we study the cascade approximation at
$\theta=\frac{9\pi}{20}$ with the Haar function as starting point. Special
attention will be given to the asymptotic properties of the jumps at certain
dyadic rationals.
The
rightmost term in the sum expansion for $M^{n}\psi^{(0)}(x)=\sum_{i} h_{i}
(n)\psi^{(0)}(2^{n} x-i)$ is
\[
h_{N_{n}}(n)\psi^{(0)}(2^{n} x-N_{n}),
\]
where $N_{n}=3(2H{n-1}+\dots+2+1)=3\cdot(2Hn -1)$, and where $\sqrt2 a_{3}
>1$, $h_{N_{n}}(n)=2^{\frac{n}{2}}a_{3}^{n}\underset{n\rightarrow\infty
}{\longrightarrow}\infty$, and each of the previous subpartition highpoints
$h_{N_{n}-k3}(n)$ contains a term with a factor $2^{\frac{n}{2}}a_{3}%
^{n-k}\underset{n\rightarrow\infty}{\longrightarrow}\infty$, $k=0,1,2,\dots$.
For example, note that $\psi^{(n)}=M^{n}\psi^{(0)}$ is supported in $\left[
0,x_{n}\right)  $, where $x_{n}=3-2^{1-n}$, and
\[
\lim_{x\rightarrow x_{n}-}\psi^{(n)}\left(  x\right)  =\left(  \sqrt{2}%
a_{3}\right)  ^{n} \underset{n\rightarrow\infty}{\longrightarrow}\infty.
\]

\begin{description}
\raggedright

\item[Figure \ref{Res045_stages}(a)] $\psi^{(1)}(x)= M\psi^{(0)}(x)=\sqrt
2\sum_{i} a_{i} \psi^{(0)}(2x-i)$; cf.\ (\ref{eq9}). The $a_{i}$'s are given
by (\ref{eqSome.2}).

\item[Figure \ref{Res045_stages}(b)] $\psi^{(2)}(x)= M^{2}\psi^{(0)}%
(x)=2\sum_{k} \sum_{i} a_{i} a_{k-2i}\psi^{(0)}(4x-k)=$ $\sum_{k} h_{k} (2)
\psi^{(0)} (4x-k)$.

\item[Figure \ref{Res045_stages}(c)] $\psi^{(3)}(x)= M^{3}\psi^{(0)}%
(x)=2^{\frac32}\sum_{k} \sum_{ij} a_{i} a_{j} a_{k-4i-2j}\psi^{(0)}(8x-k)
=\sum_{k} h_{k} (3) \psi^{(0)} (8x-k)$.

\item[Figure \ref{Res045_stages}(d)] $\psi^{(4)}(x)= M^{4}\psi^{(0)}%
(x)=\sum_{k}h_{k}(4)\psi^{(0)}(2^{4}x-k)$, where the height of the rightmost
column is $h_{45}(4)=(\sqrt{2}a_{3})^{4}$, and further high points at
$h_{42}(4)$, $h_{39}(4)$, $h_{36}(4)$, $\dots$, corresponding to rightmost
``bumps'' in the subpartition intervals.

\noindent$\vdots$

\item[Figure \ref{Res045_stages}(l)] The $n=12$ case of $\psi^{(n)}(x)=
M^{n}\psi^{(0)}(x)=$ $\sum_{k}h_{k}(n)\psi^{(0)}(2^{n}x-k)$, where $2^{-n}%
\sum_{k}h_{k}(n)^{2}=\Vert M^{n}\psi^{(0)}\Vert_{2}^{2}\Vert\psi^{(0)}%
\Vert_{2}^{2}=1$ by (\ref{eqOla2.82}), despite the divergence at the high
rightmost ``bumps'', $n\rightarrow\infty$.
\end{description}
\bigskip

\noindent
\textit{Additional comments on the Figures:}
The numerical sizes of the $y$-coordinates
of the cusp points in Figure \ref{Movie}
are computed and displayed in the Appendix.
Note that figures \ref{Movie}(b)--\ref{Movie}(j)
represent continuous functions, while
all other figures represent discontinuous functions
\cite[p.~193]{CoHe94}. See more details of Figure
\ref{Movie}(b) in Figure~\ref{Giraffe}.

   See \texttt{http://cm.bell-labs.com/who/wim/cascade/} for more pictures, but
without emphasis of the singularities at dyadic points.

In Figure \ref{ckplot}, note that everywhere except
at the four points $\theta =k\frac{\pi}{2}$, $k=0,1,2,3$,
we have precisely three positive $a_{i}$'s and
one negative one. The significance of that is discussed
in Section \ref{Con} below.

\addtocounter{figure}{-1}

\begin{figure}[ptb]
\begin{picture}(351,480)
\put(0,264){\includegraphics
[bb=189 166 422 596,height=216bp,width=117bp]{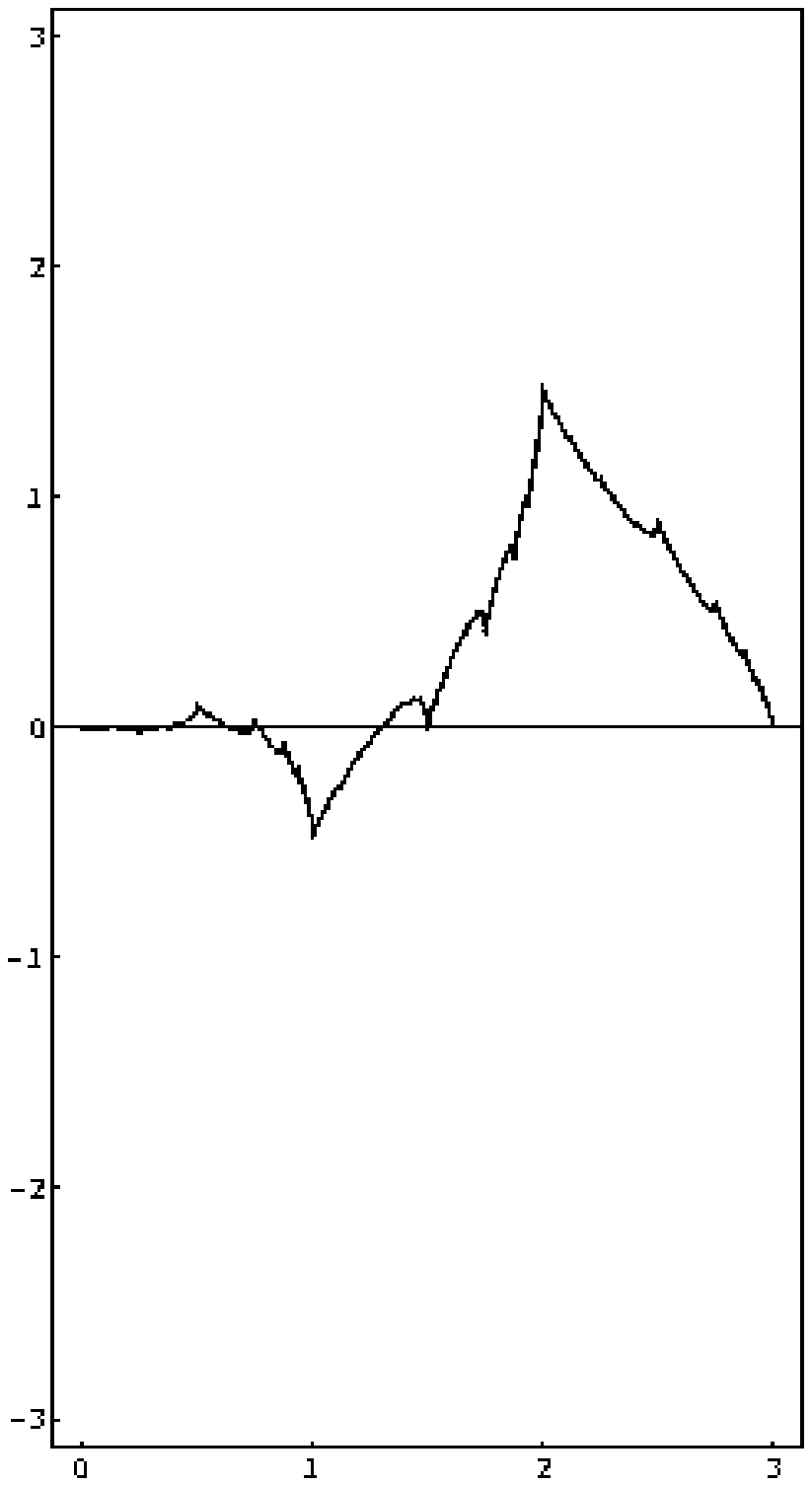}}
\put(117,264){\includegraphics
[bb=189 166 422 596,height=216bp,width=117bp]{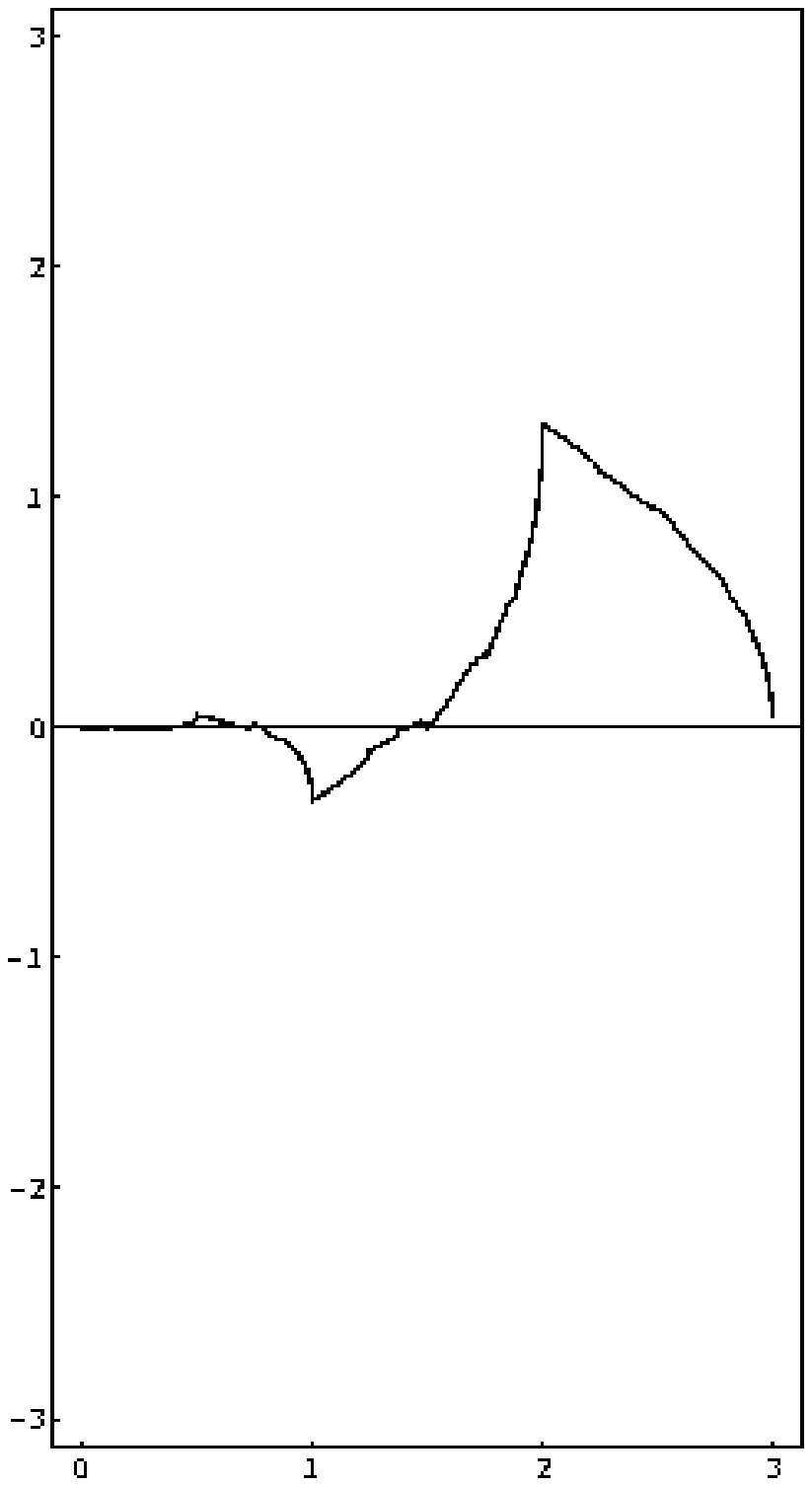}}
\put(234,264){\includegraphics
[bb=189 166 422 596,height=216bp,width=117bp]{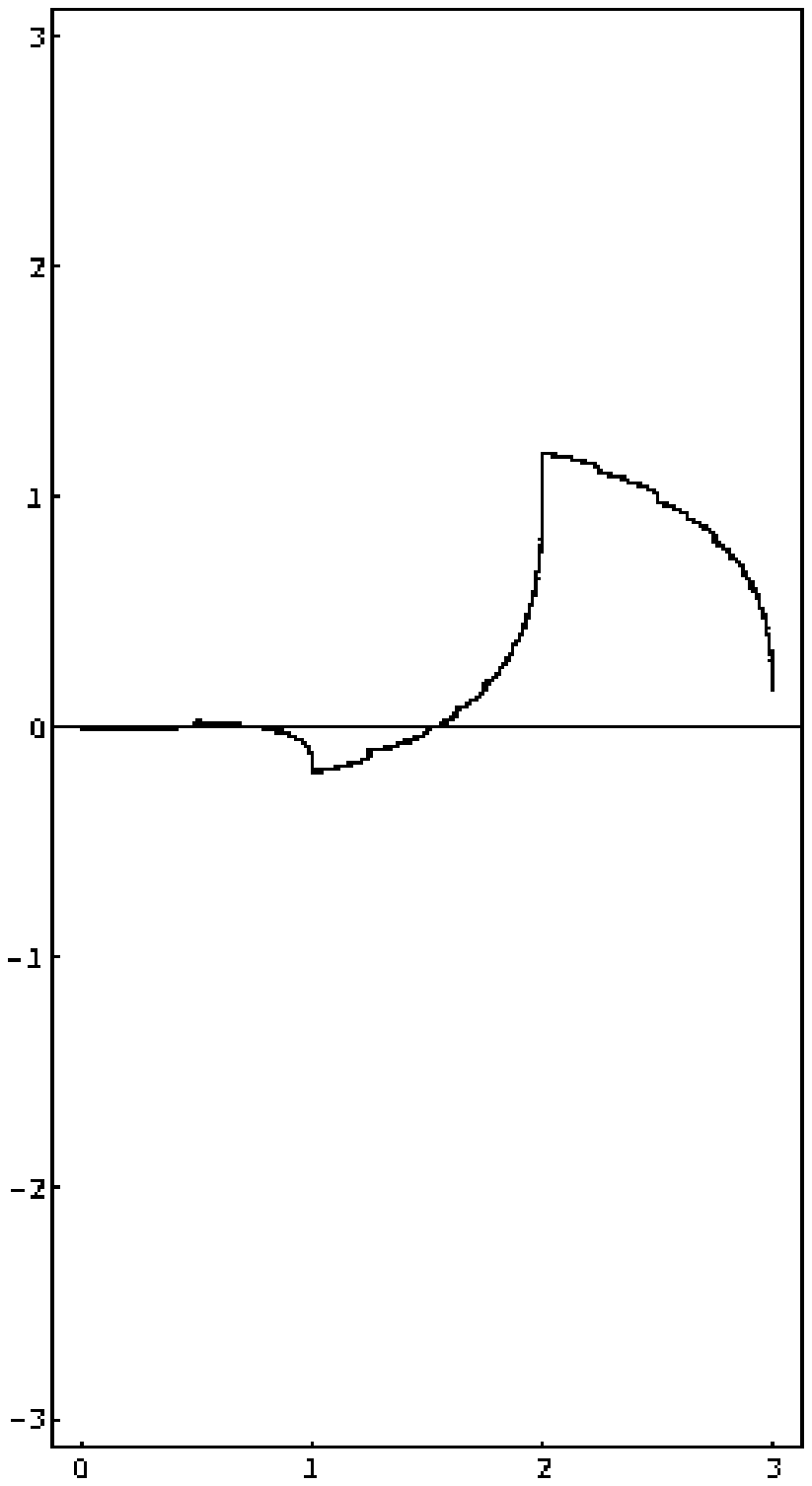}}
\put(0,249){\makebox(117,12){g: $\theta=-\frac{\pi}{5}$}}
\put(117,249){\makebox(117,12){h: $\theta=-\frac{3\pi}{20}$}}
\put(234,249){\makebox(117,12){i: $\theta=-\frac{\pi}{10}$}}
\put(0,18){\includegraphics
[bb=189 166 422 596,height=216bp,width=117bp]{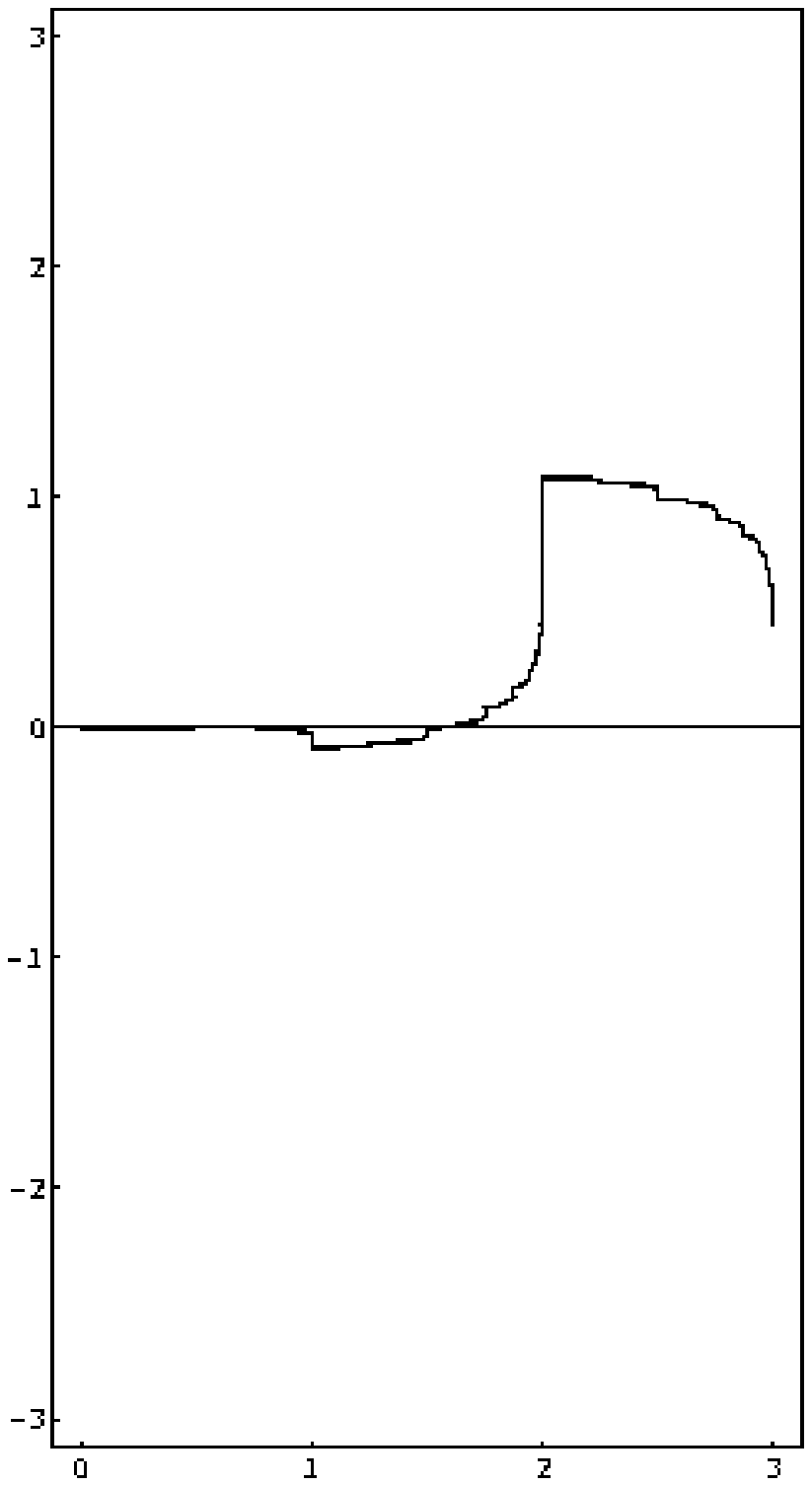}}
\put(117,18){\includegraphics
[bb=99 0 333 432,height=216bp,width=117bp]{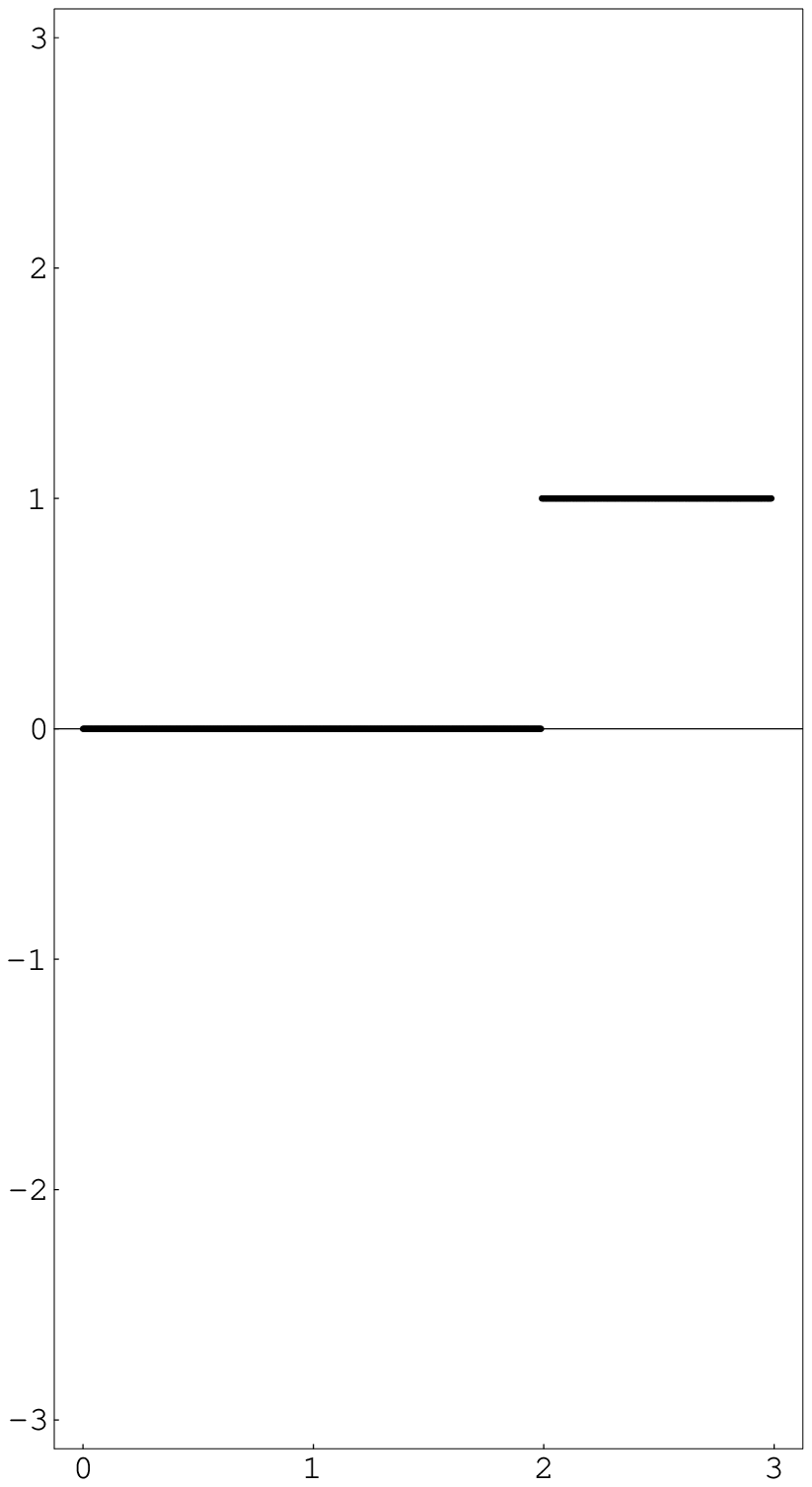}}
\put(234,18){\includegraphics
[bb=99 0 333 432,height=216bp,width=117bp]{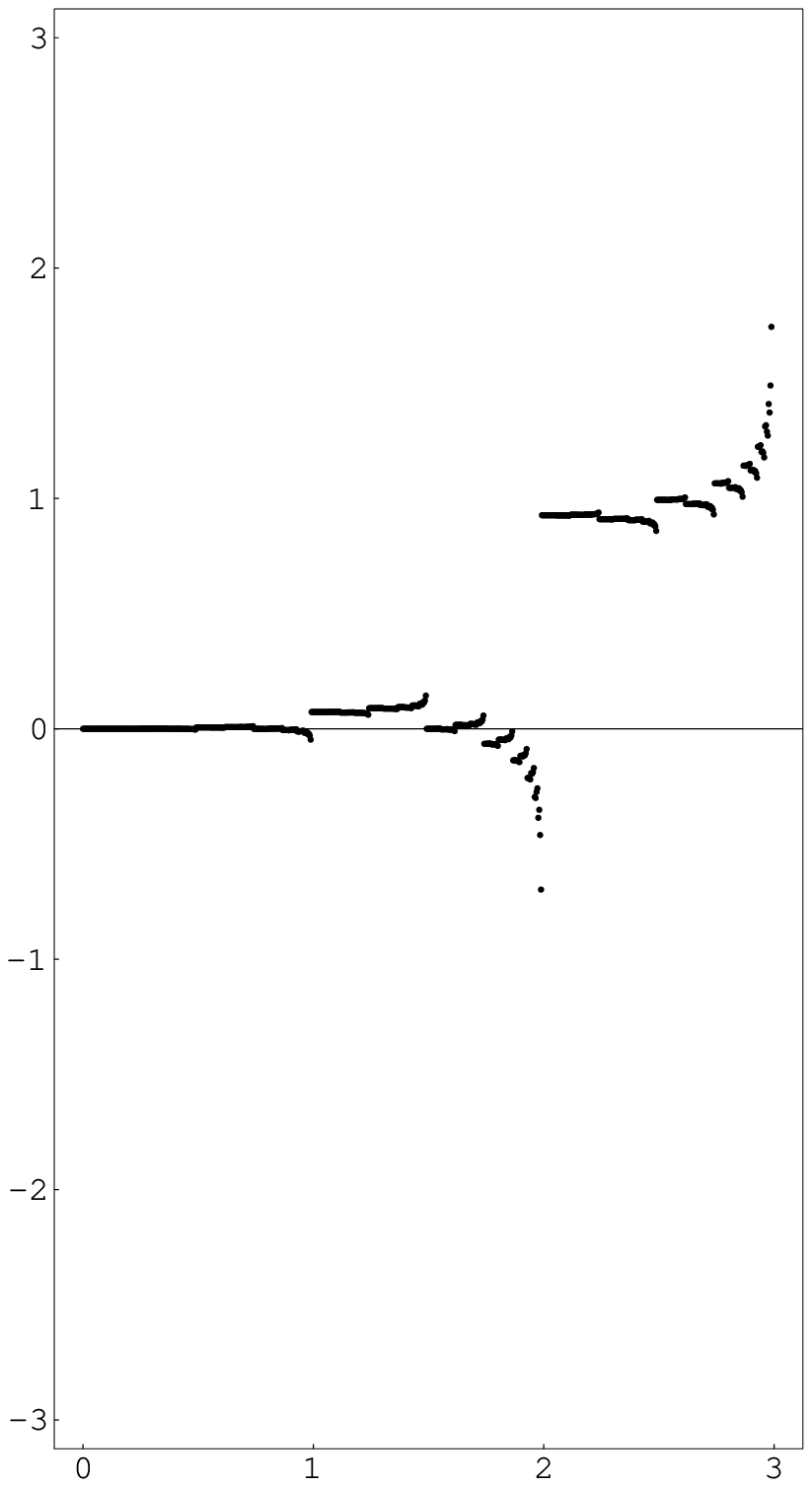}}
\put(0,3){\makebox(117,12){j: $\theta=-\frac{\pi}{20}$}}
\put(117,3){\makebox(117,12){k: $\theta=0$}}
\put(234,3){\makebox(117,12){l: $\theta=\frac{\pi}{20}$}}
\end{picture}
\caption{Scaling function ``movie reel'', $\theta$ from $\frac{-\pi}{2}$ to
$\frac{\pi}{2}$ in twenty-one frames: Frames g--l}%
\end{figure}

\addtocounter{figure}{-1}

\begin{figure}[ptb]
\begin{picture}(351,480)
\put(0,264){\includegraphics
[bb=99 0 333 432,height=216bp,width=117bp]{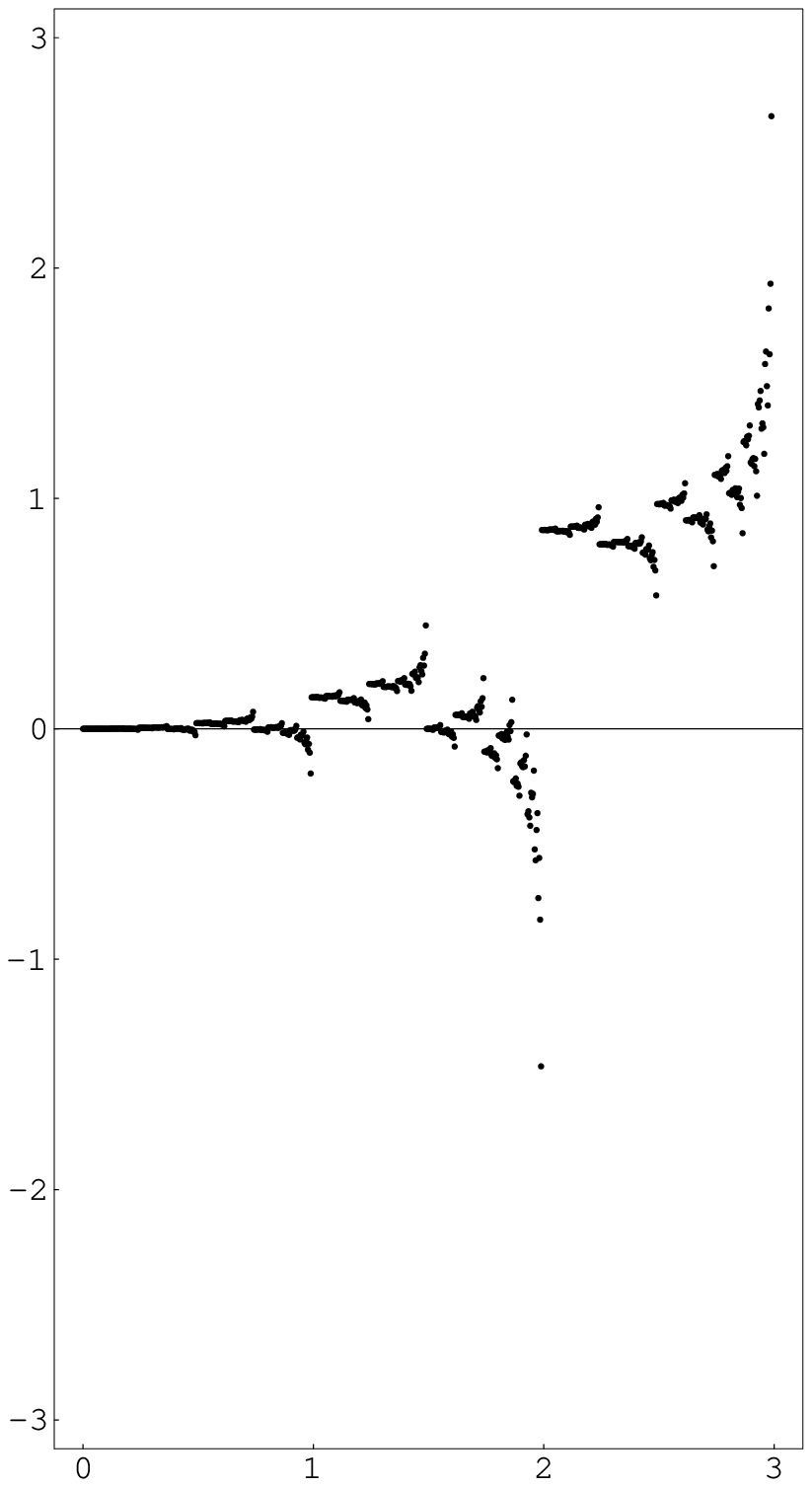}}
\put(117,264){\includegraphics
[bb=99 0 333 432,height=216bp,width=117bp]{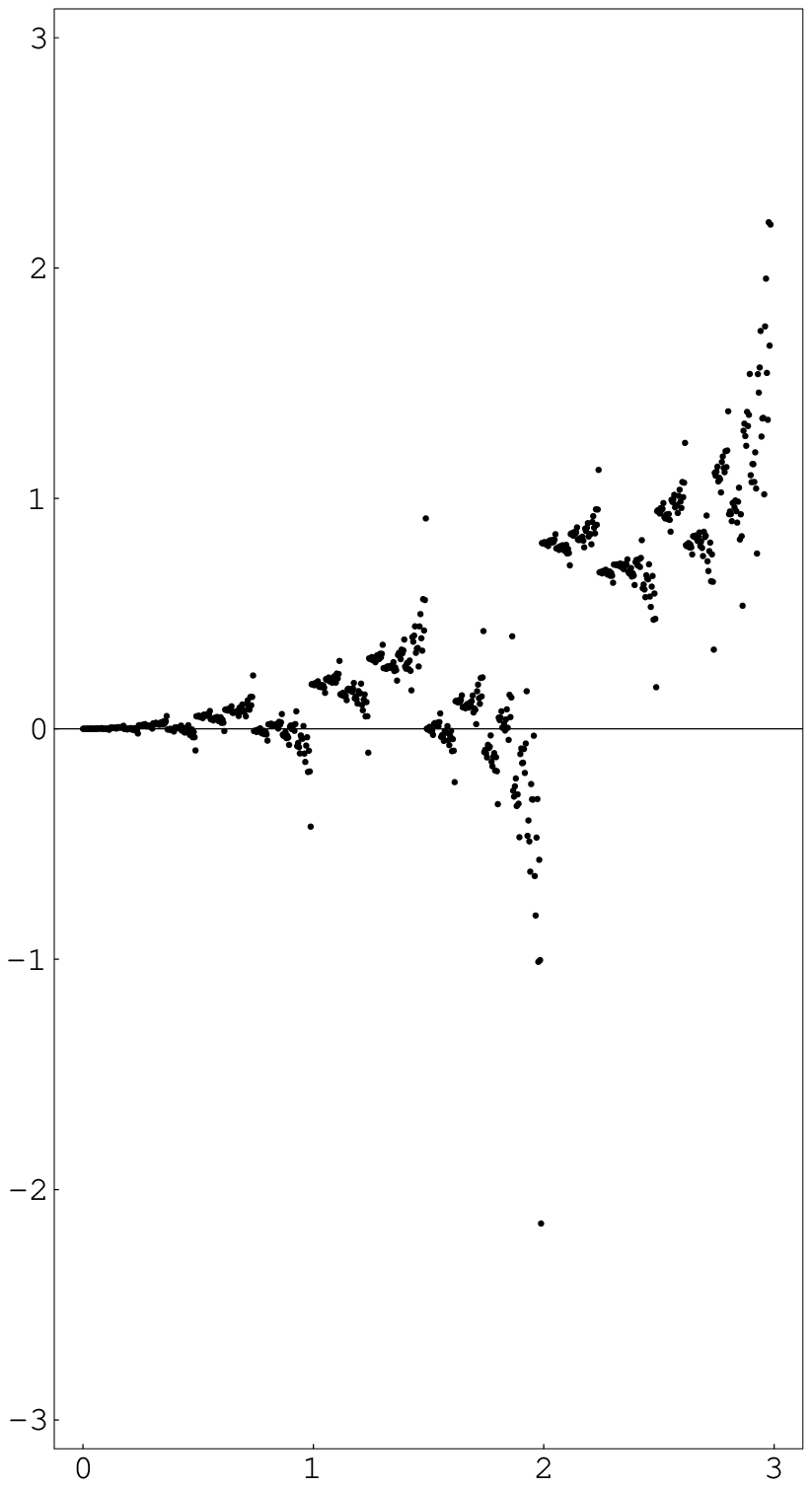}}
\put(234,264){\includegraphics
[bb=99 0 333 432,height=216bp,width=117bp]{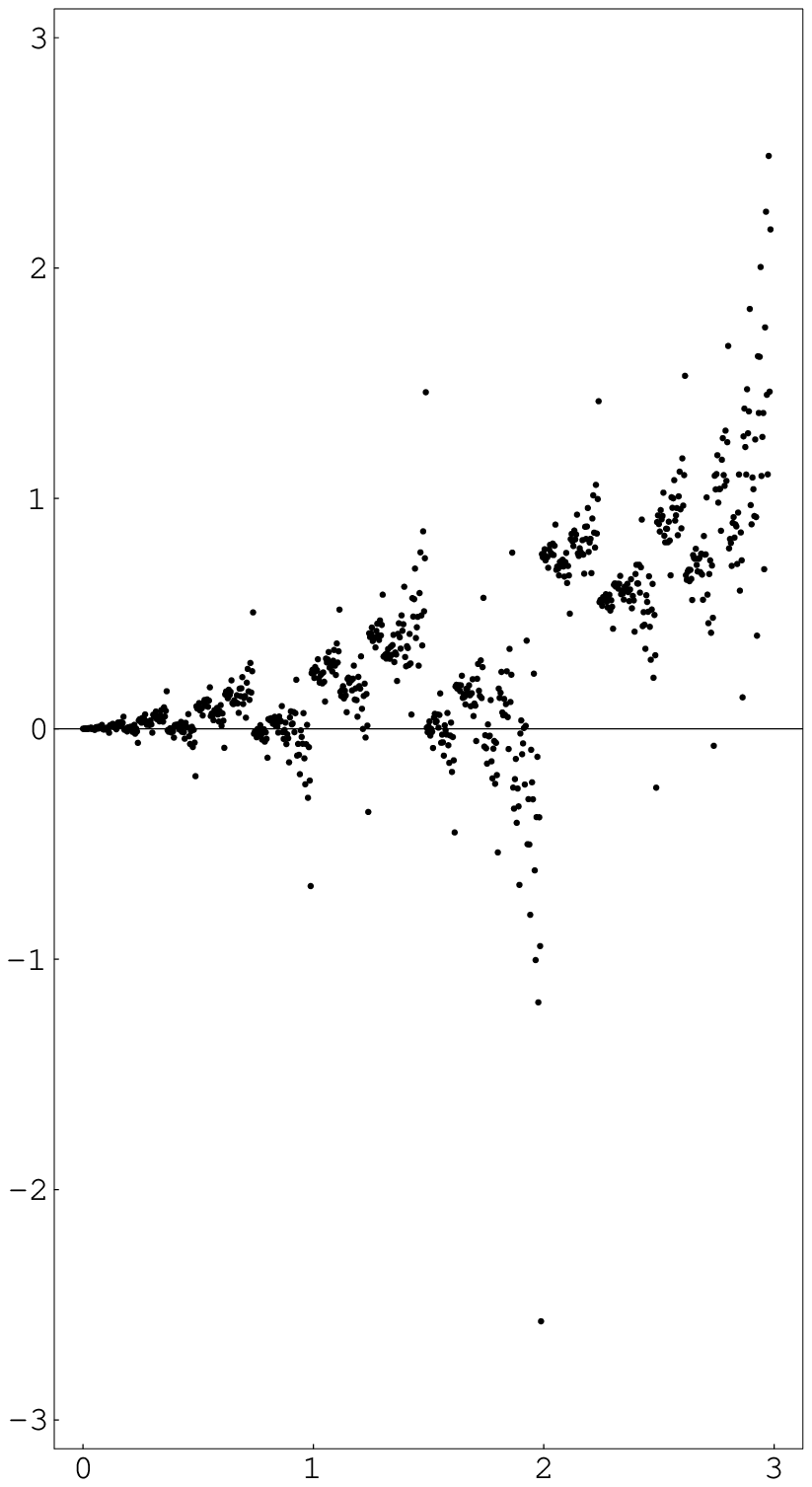}}
\put(0,249){\makebox(117,12){m: $\theta=\frac{\pi}{10}$}}
\put(117,249){\makebox(117,12){n: $\theta=\frac{3\pi}{20}$}}
\put(234,249){\makebox(117,12){o: $\theta=\frac{\pi}{5}$}}
\put(0,18){\includegraphics
[bb=99 0 333 432,height=216bp,width=117bp]{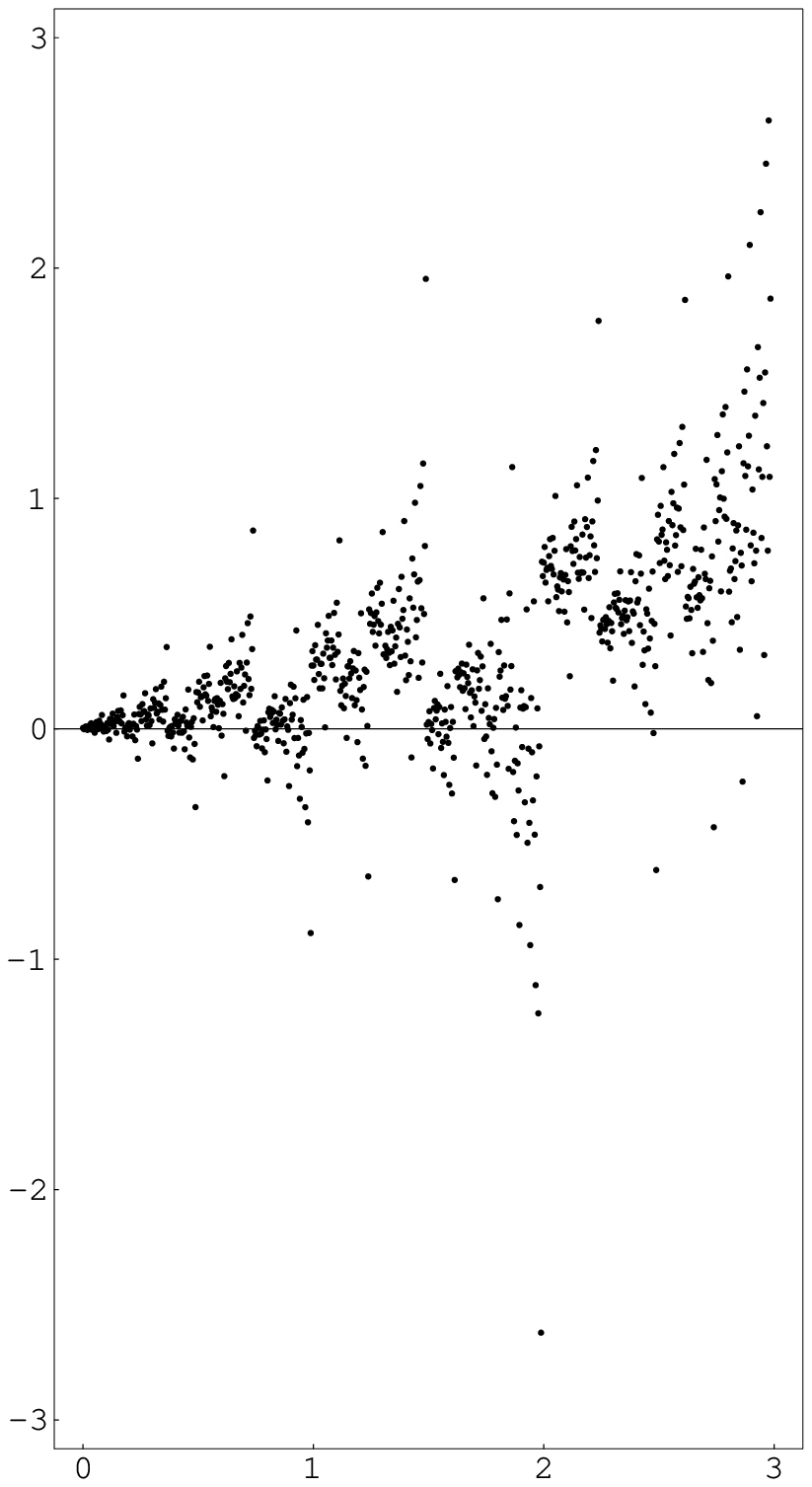}}
\put(117,18){\includegraphics
[bb=99 0 333 432,height=216bp,width=117bp]{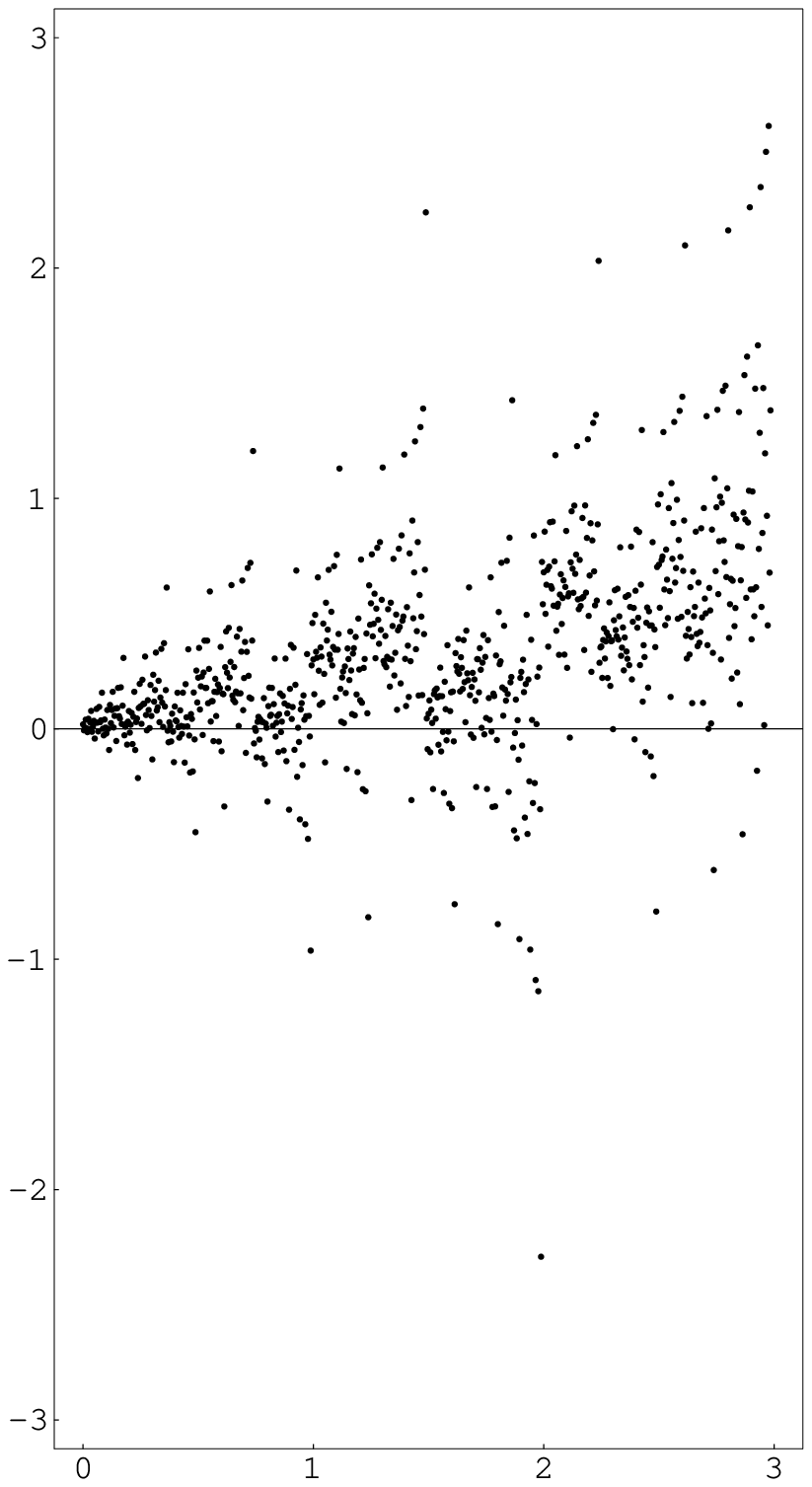}}
\put(234,18){\includegraphics
[bb=99 0 333 432,height=216bp,width=117bp]{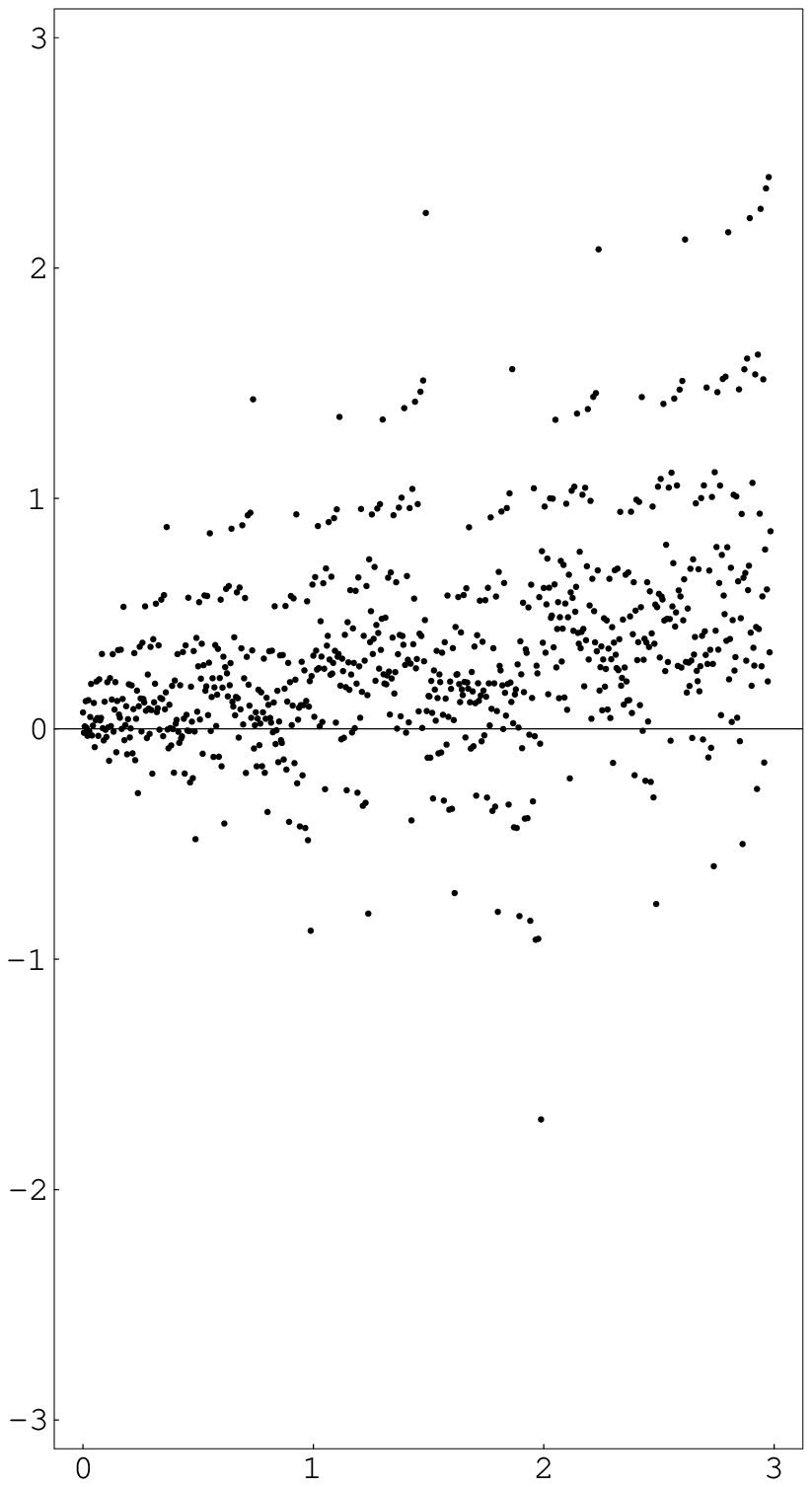}}
\put(0,3){\makebox(117,12){p: $\theta=\frac{\pi}{4}$}}
\put(117,3){\makebox(117,12){q: $\theta=\frac{3\pi}{10}$}}
\put(234,3){\makebox(117,12){r: $\theta=\frac{7\pi}{20}$}}
\end{picture}
\caption{Scaling function ``movie reel'', $\theta$ from $\frac{-\pi}{2}$ to
$\frac{\pi}{2}$ in twenty-one frames: Frames m--r}%
\end{figure}

\addtocounter{figure}{-1}

\begin{figure}[ptb]
\begin{picture}(351,243)(0,-9)
\put(0,18){\includegraphics
[bb=99 0 333 432,height=216bp,width=117bp]{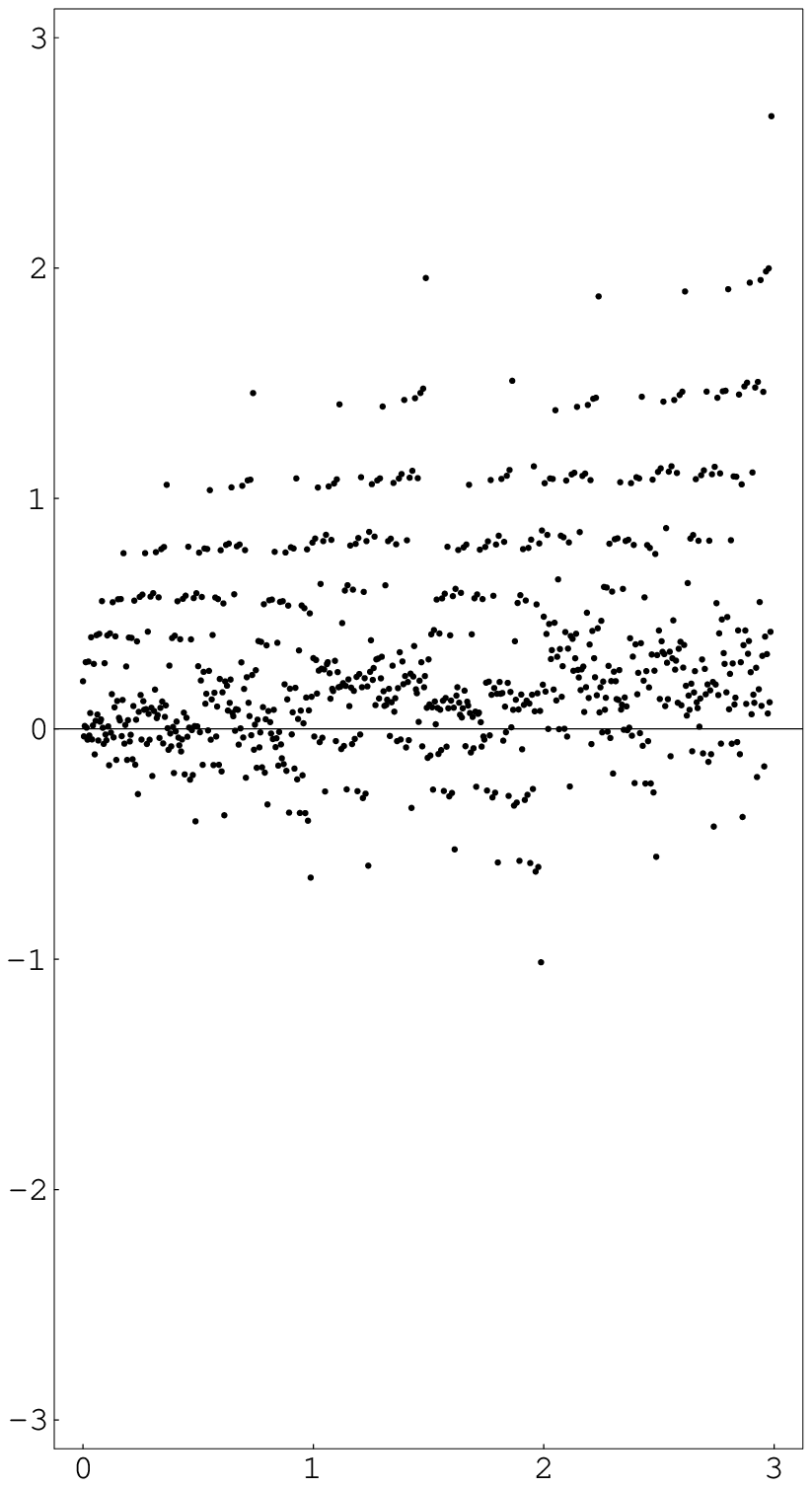}}
\put(117,18){\includegraphics
[bb=99 0 333 432,height=216bp,width=117bp]{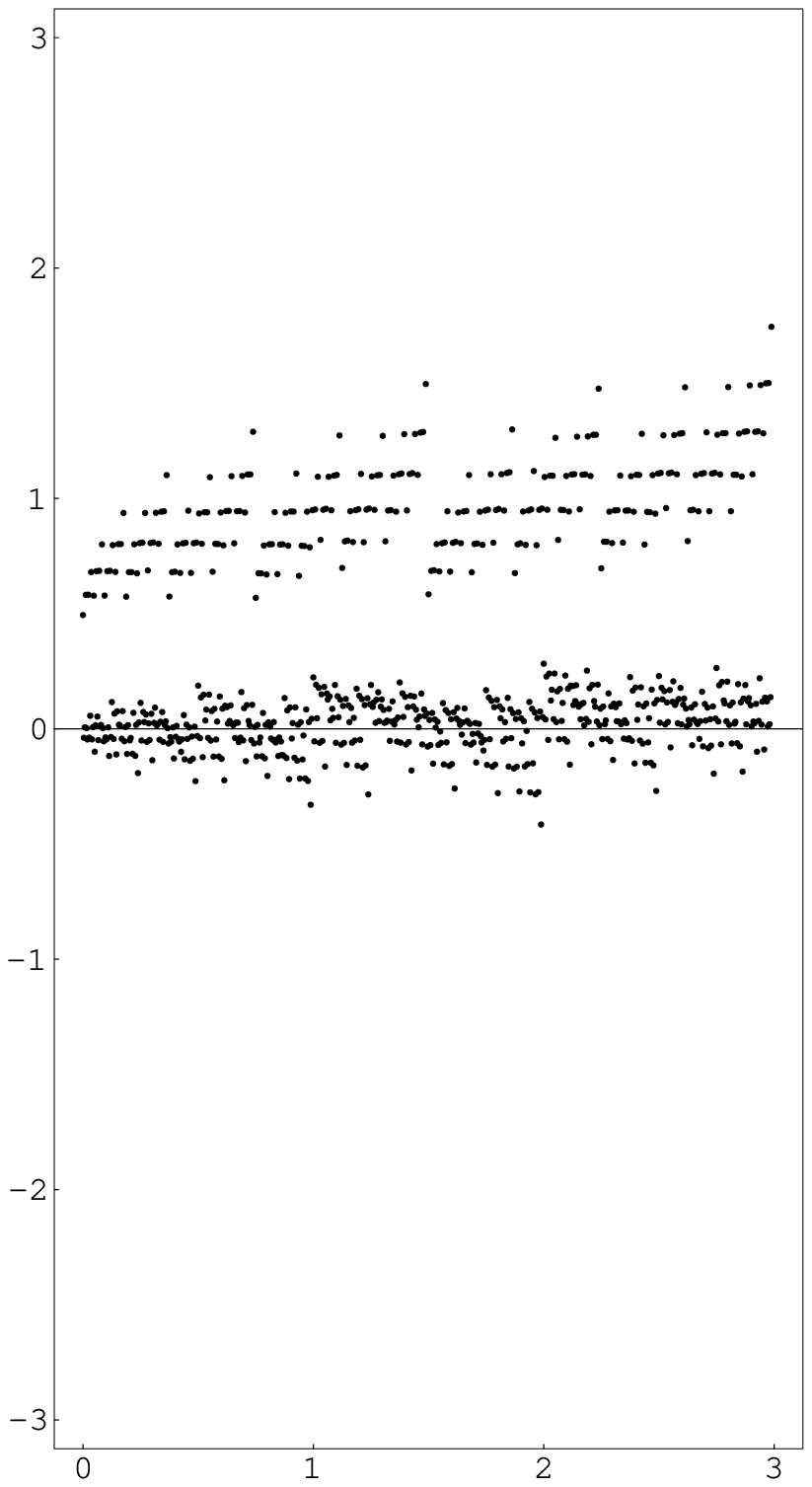}}
\put(234,18){\includegraphics
[bb=99 0 333 432,height=216bp,width=117bp]{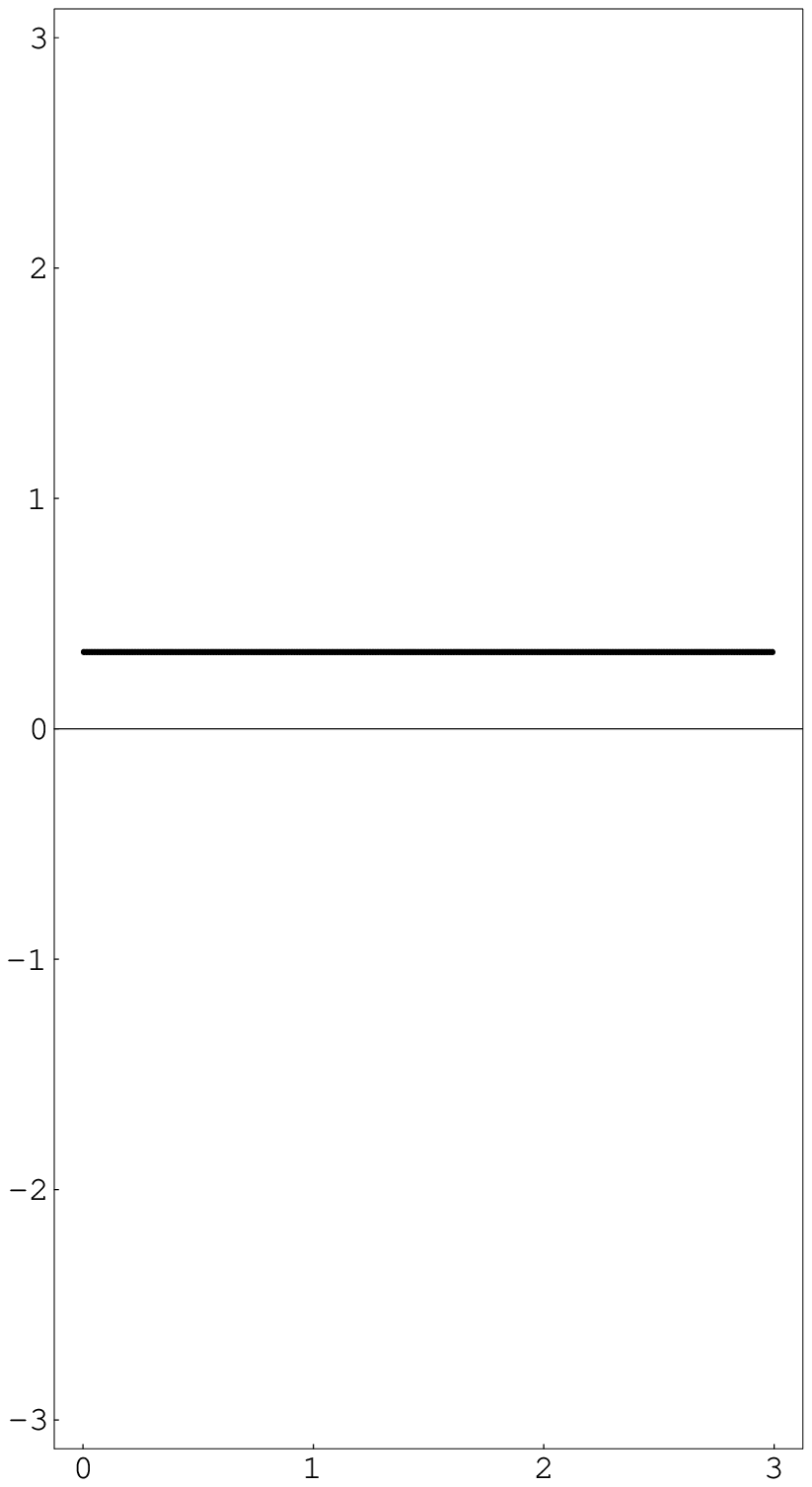}}
\put(0,3){\makebox(117,12){s: $\theta=\frac{2\pi}{5}$}}
\put(117,3){\makebox(117,12){t: $\theta=\frac{9\pi}{20}$}}
\put(117,-9){\makebox(117,12){(see also Figure \ref{Res045_stages}(h))}}
\put(234,3){\makebox(117,12){u: $\theta=\frac{\pi}{2}$}}
\put(234,-9){\makebox(117,12){(averaged)}}
\end{picture}
\caption{Scaling function ``movie reel'', $\theta$ from $\frac{-\pi}{2}$ to
$\frac{\pi}{2}$ in twenty-one frames: Frames s--u}%
\label{Movie_end}%
\end{figure}

\begin{figure}[ptb]
\begin{picture}(360,354)(0,-9)
\put(0,18){\includegraphics
[bb=40 58 571 704,height=330bp]{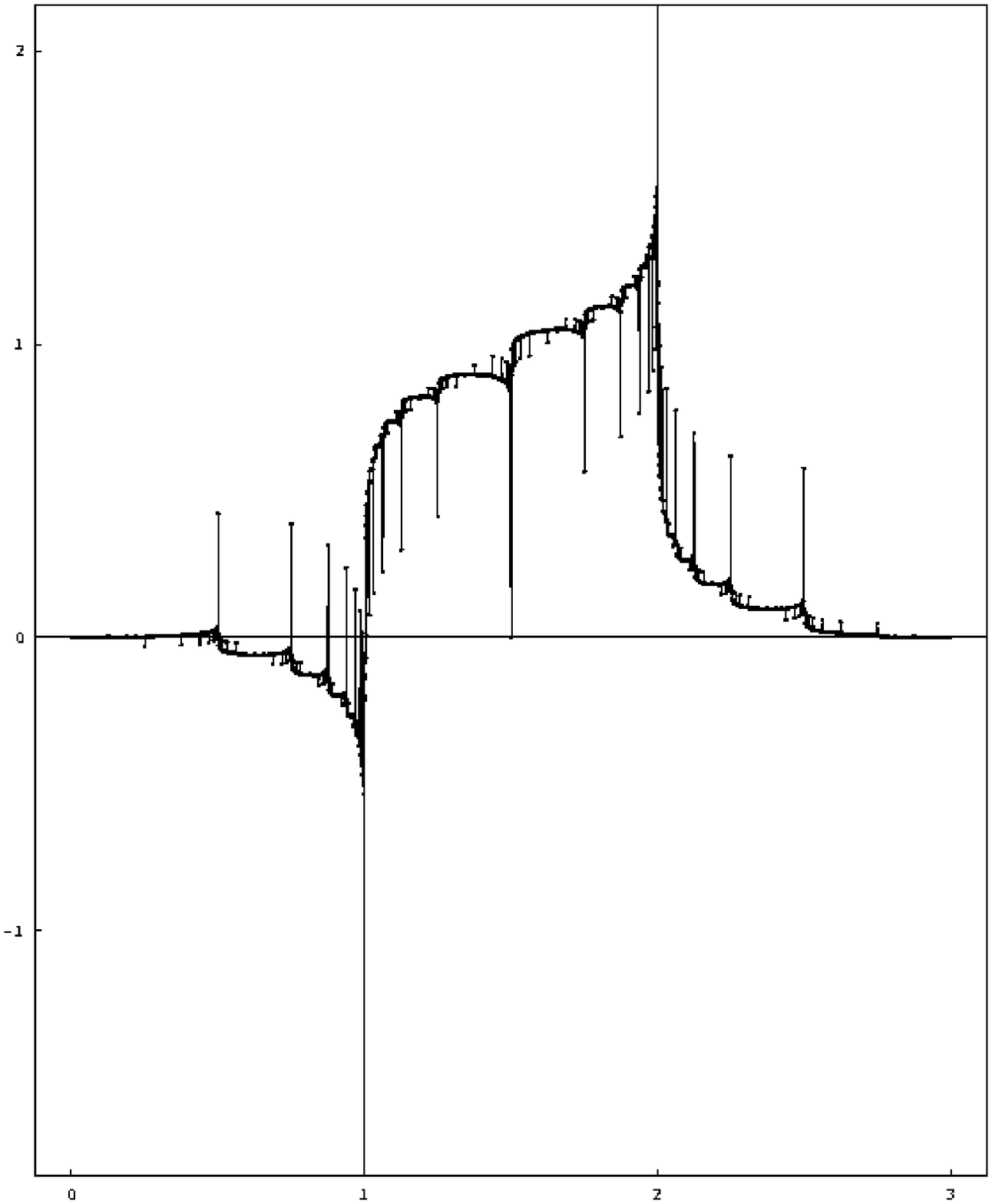}}
\put(273,18){\includegraphics
[bb=222 58 390 704,height=330bp]{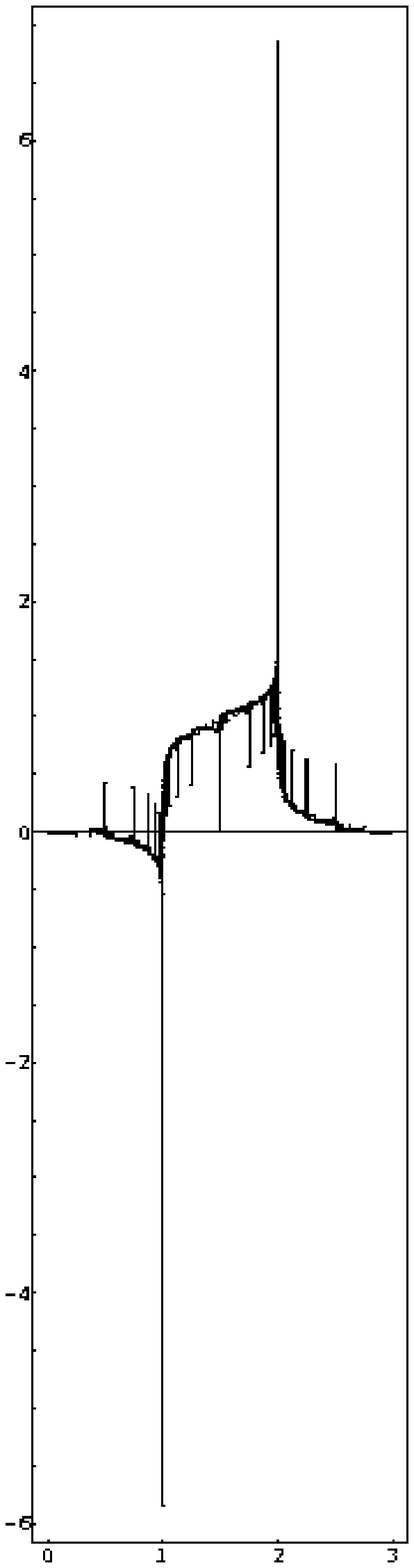}}
\put(0,3){\makebox(273,12){(detail)}}
\end{picture}
\caption{The $m\rightarrow \infty $ limit
of $\psi^{\left( m\right) }\left( x\right) $
for $\theta = -\frac{9\pi}{20}$}%
\label{Giraffe}
\end{figure}

\begin{figure}[ptb]
\includegraphics
[bb=29 0 316 360,height=216bp]{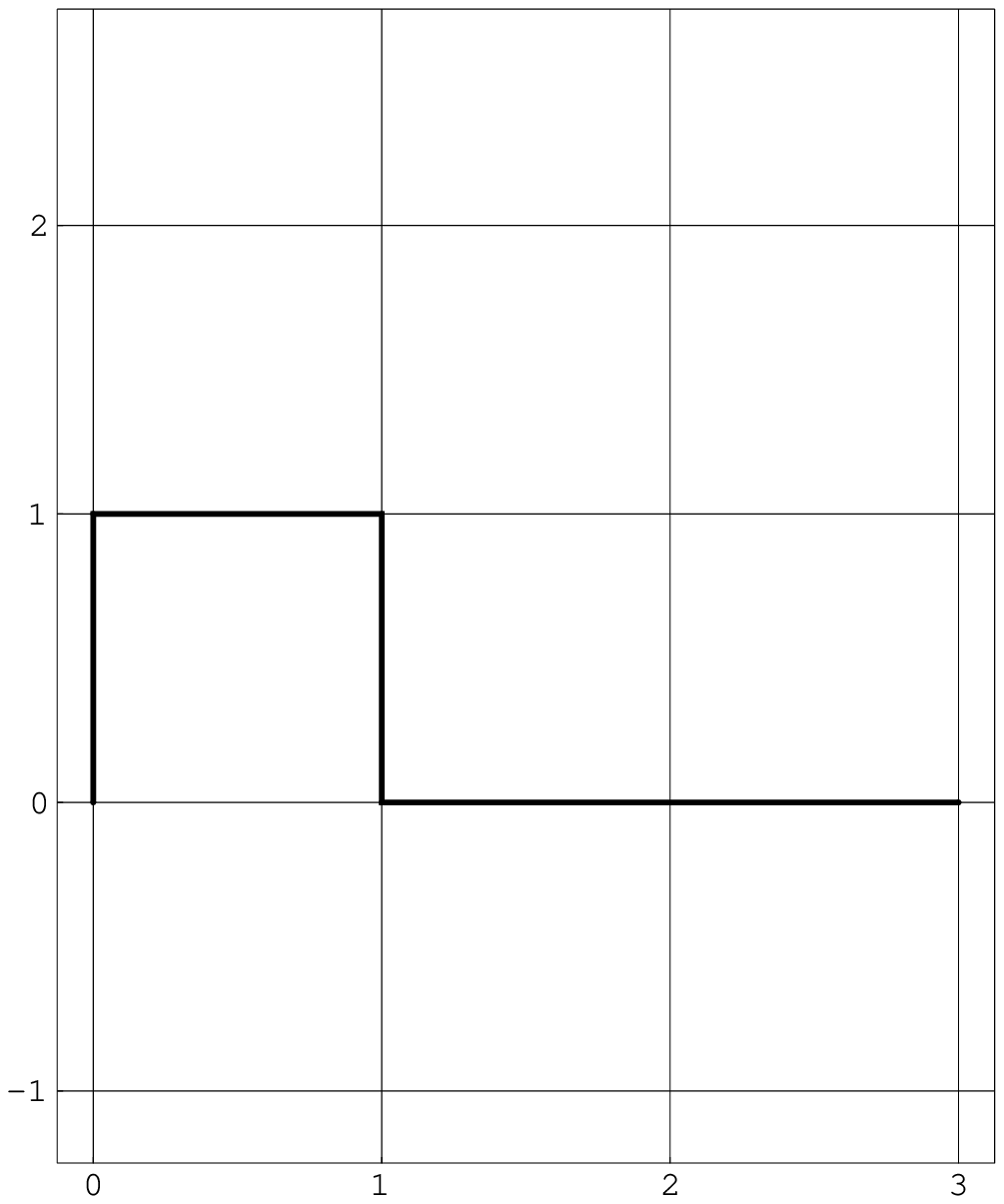}\caption{Initial function
$\psi^{(0)}(x)=\chi_{[0,1]}(x)$ for the cascade series (Haar scaling
function)}%
\label{Res045_00}%
\end{figure}

\begin{figure}[ptb]
\begin{picture}(344,462)
\put(0,255){\includegraphics
[bb=29 0 316 345,height=207bp]{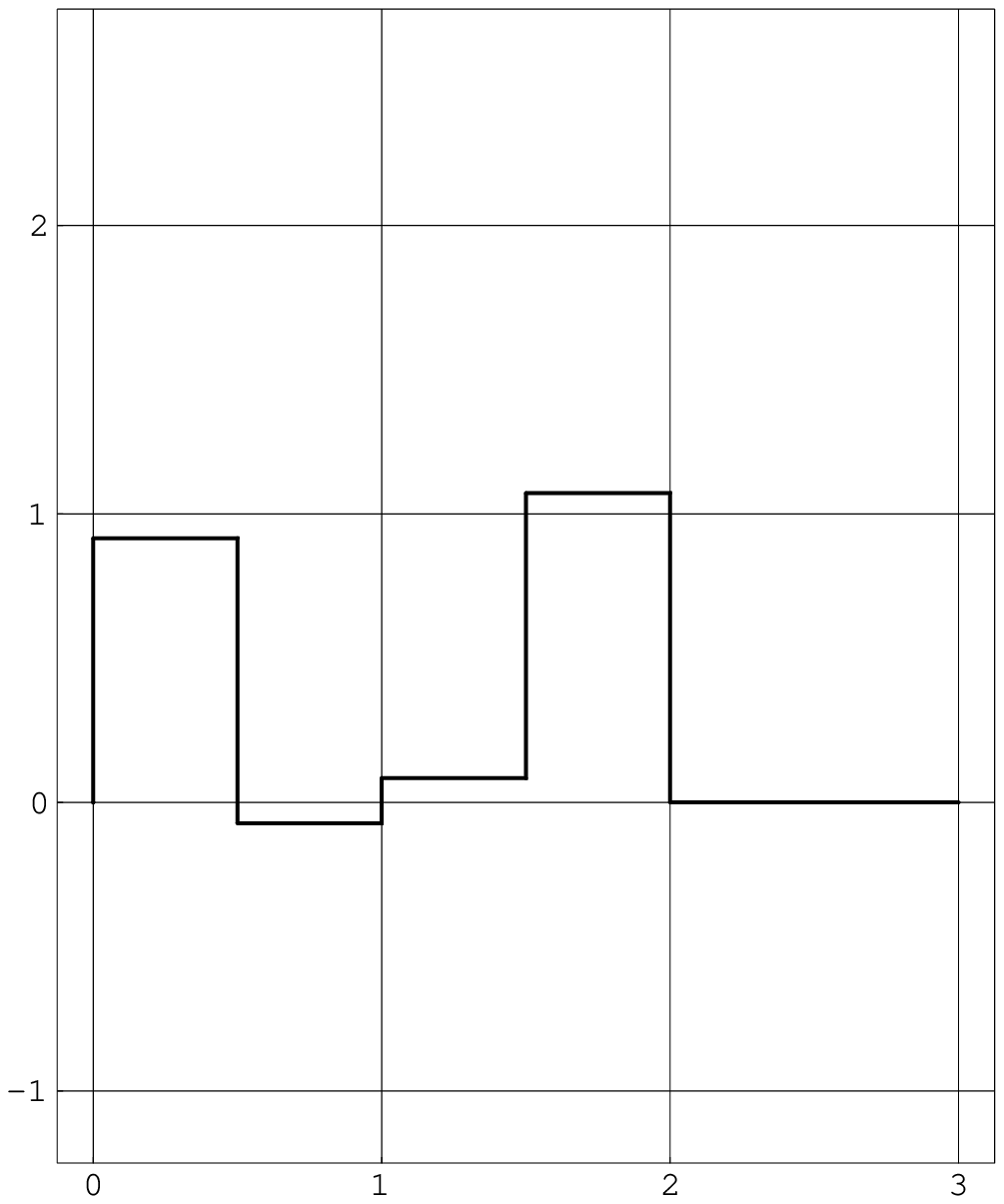}}
\put(28,364){\makebox(0,0){$\sqrt{2}a_0$}}
\put(53,314){\makebox(0,0){$\sqrt{2}a_1$}}
\put(78,336){\makebox(0,0){$\sqrt{2}a_2$}}
\put(103,386){\makebox(0,0){$\sqrt{2}a_3$}}
\put(172,255){\includegraphics
[bb=29 0 316 345,height=207bp]{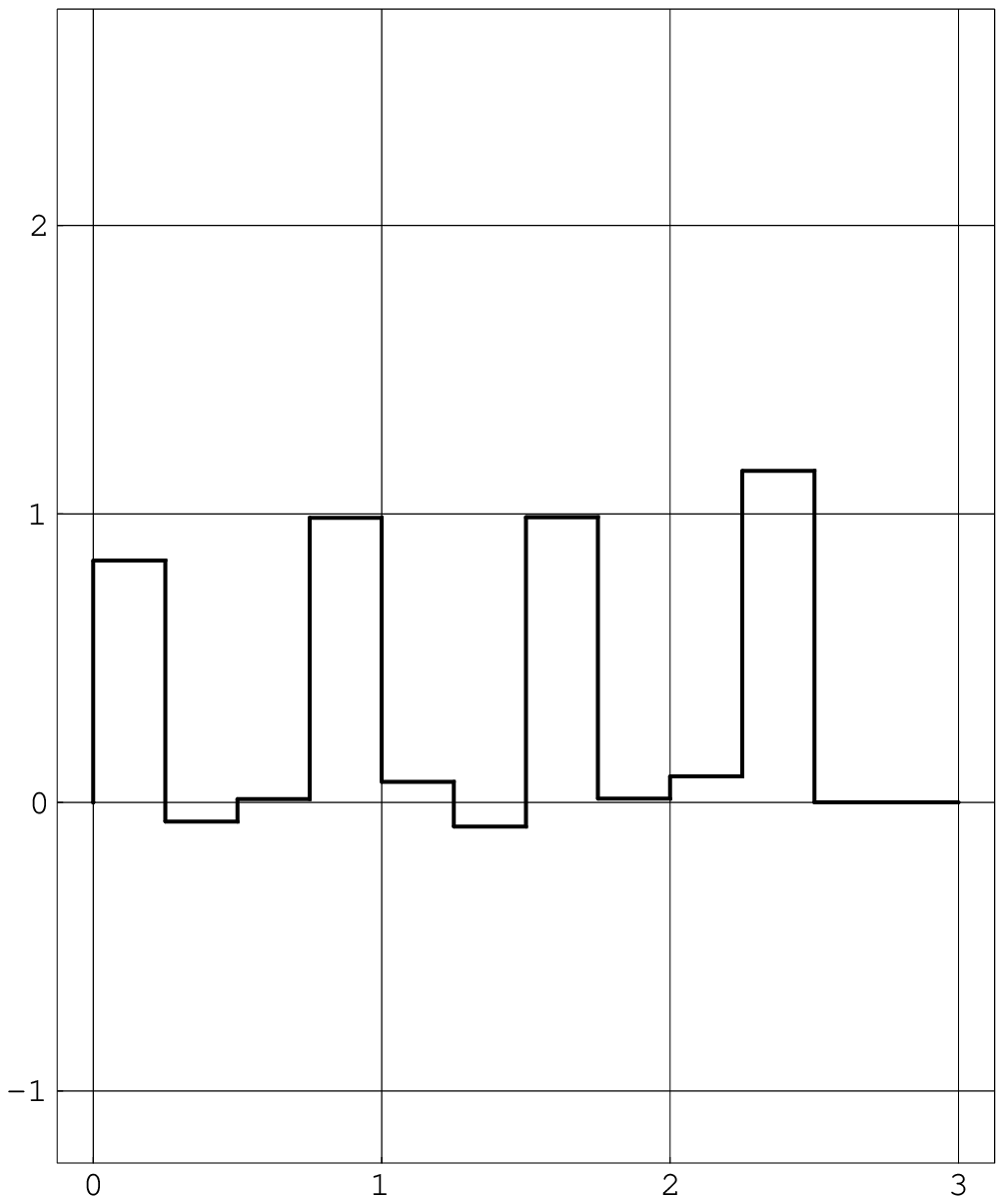}}
\put(0,240){\makebox(172,12){a: Cascade stage $1$}}
\put(172,240){\makebox(172,12){b: Cascade stage $2$}}
\put(0,18){\includegraphics
[bb=29 0 316 345,height=207bp]{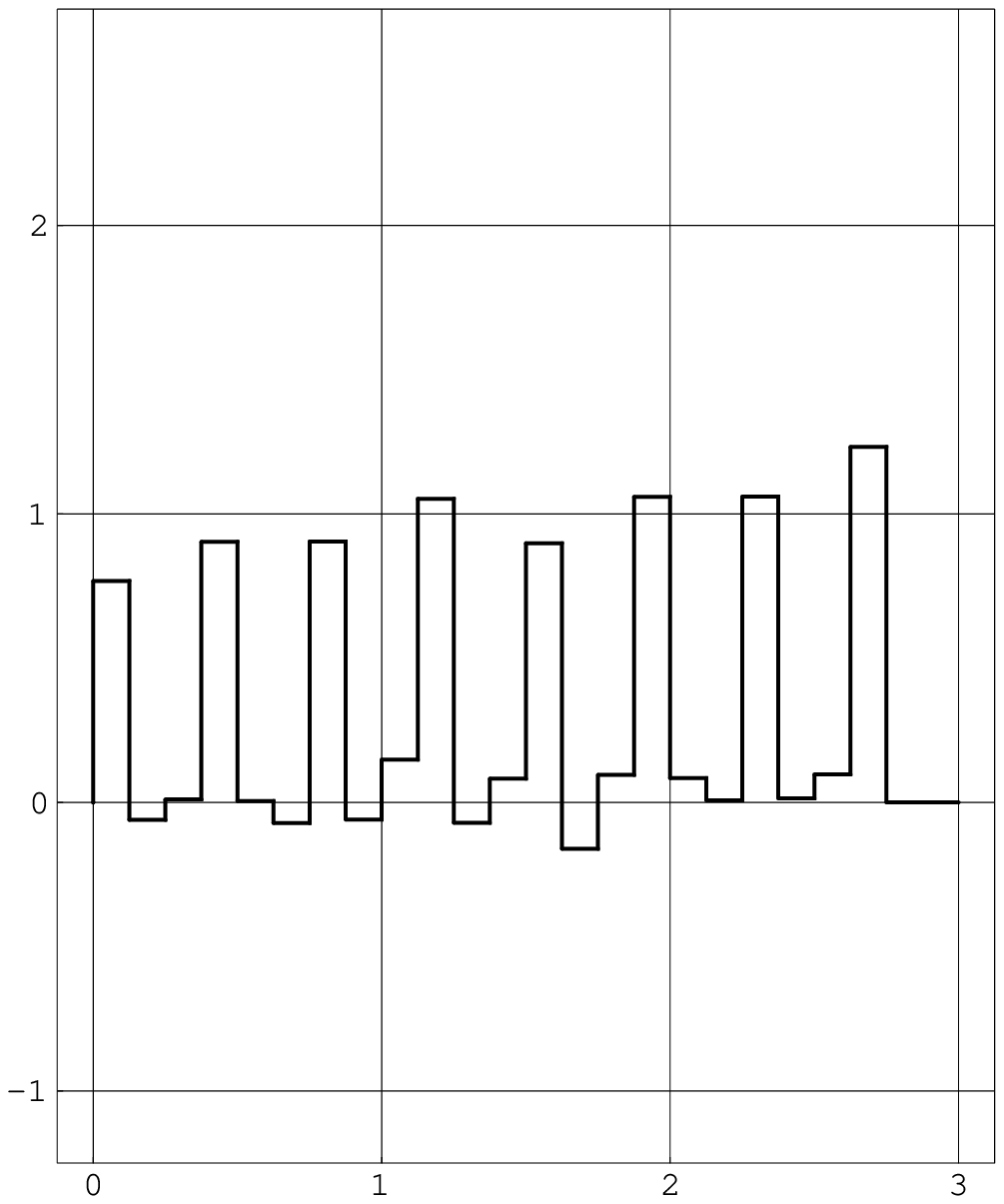}}
\put(172,18){\includegraphics
[bb=29 0 316 345,height=207bp]{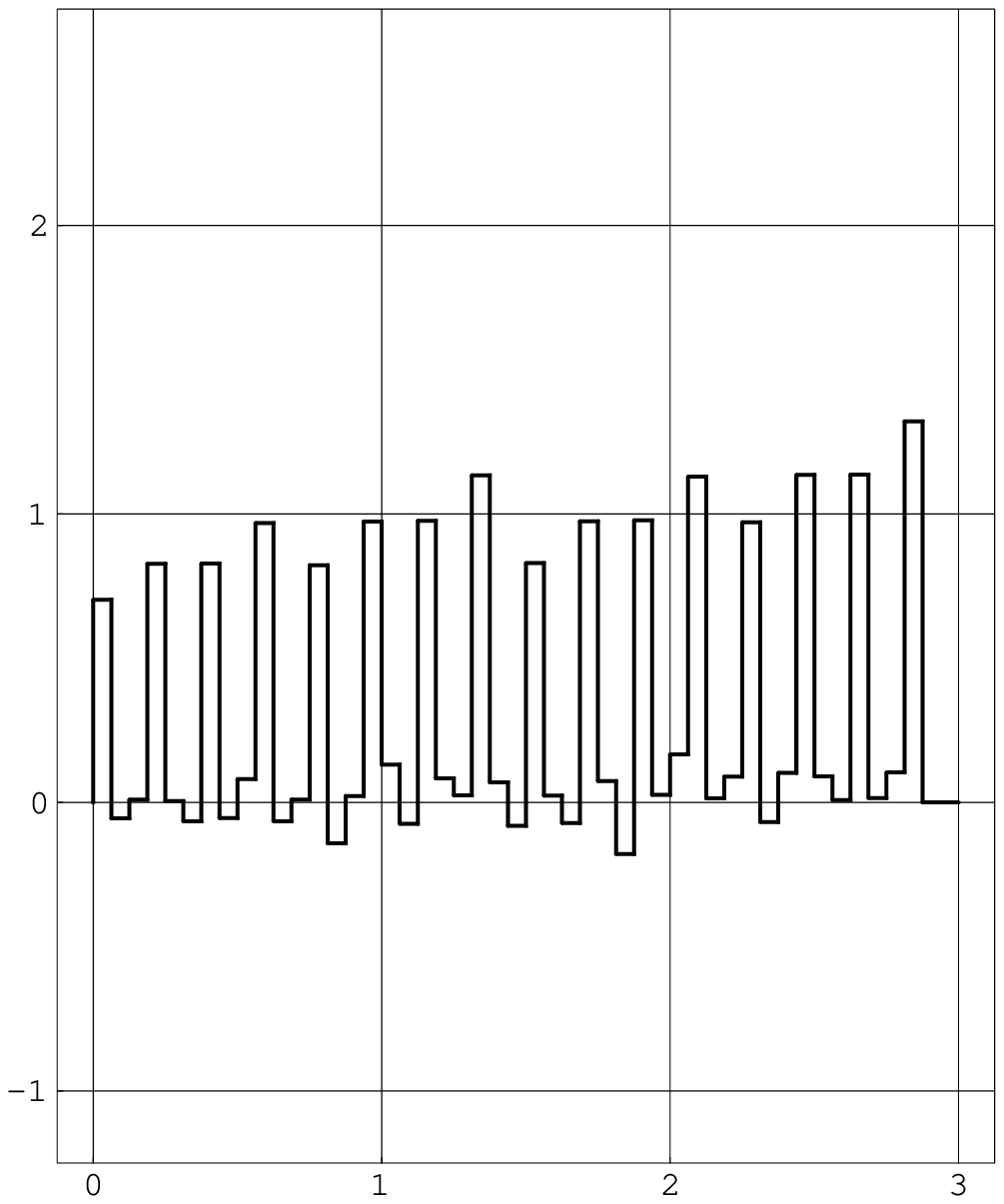}}
\put(0,3){\makebox(172,12){c: Cascade stage $3$}}
\put(172,3){\makebox(172,12){d: Cascade stage $4$}}
\end{picture}
\caption{Cascade stages of scaling function, $\theta=\frac{9\pi}{20}$: Stages
$1$--$4$}%
\label{Res045_stages}%
\end{figure}

\addtocounter{figure}{-1}

\begin{figure}[ptb]
\begin{picture}(344,471)(0,-9)
\put(0,255){\includegraphics
[bb=29 0 316 345,height=207bp]{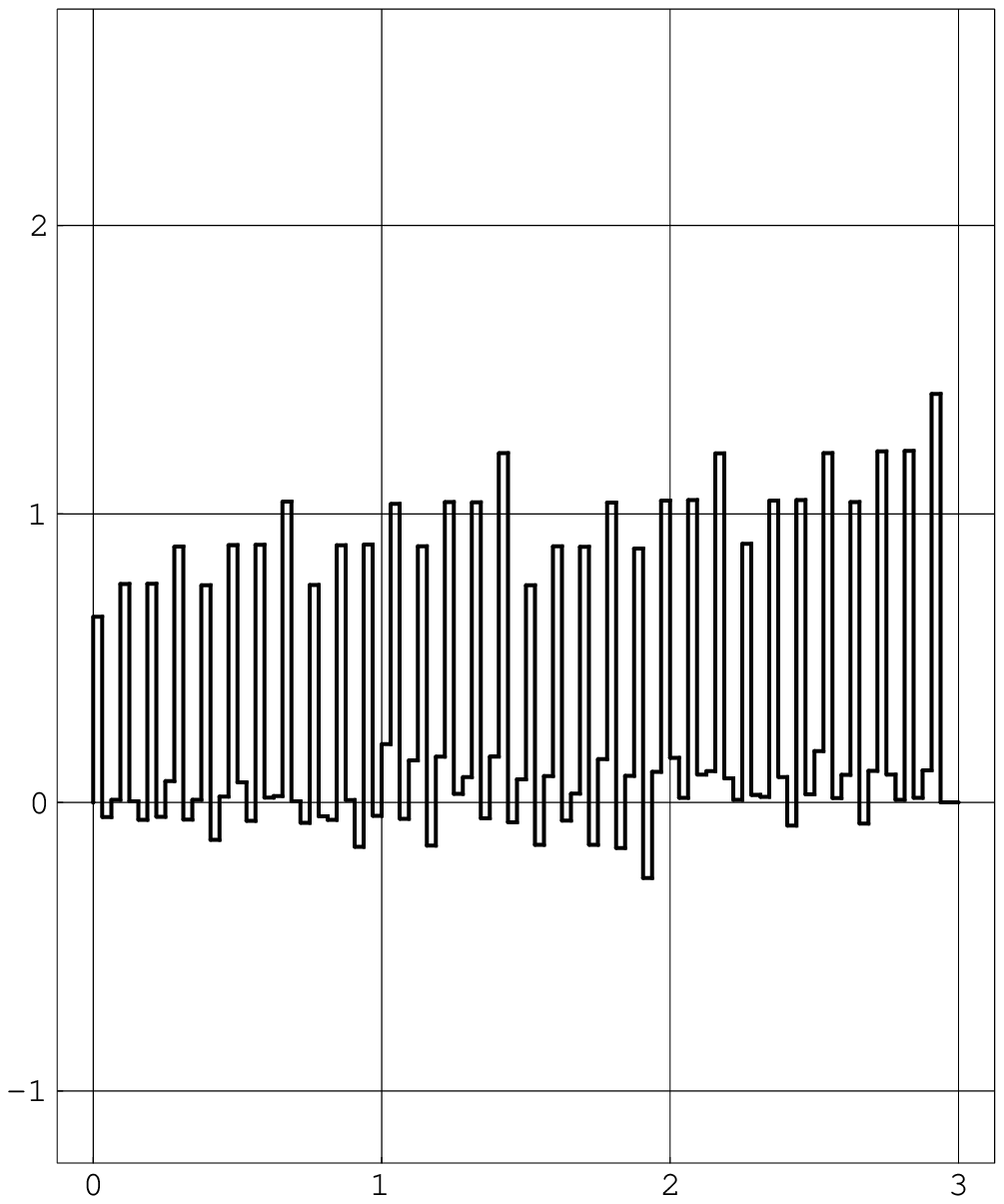}}
\put(172,255){\includegraphics
[bb=29 0 316 345,height=207bp]{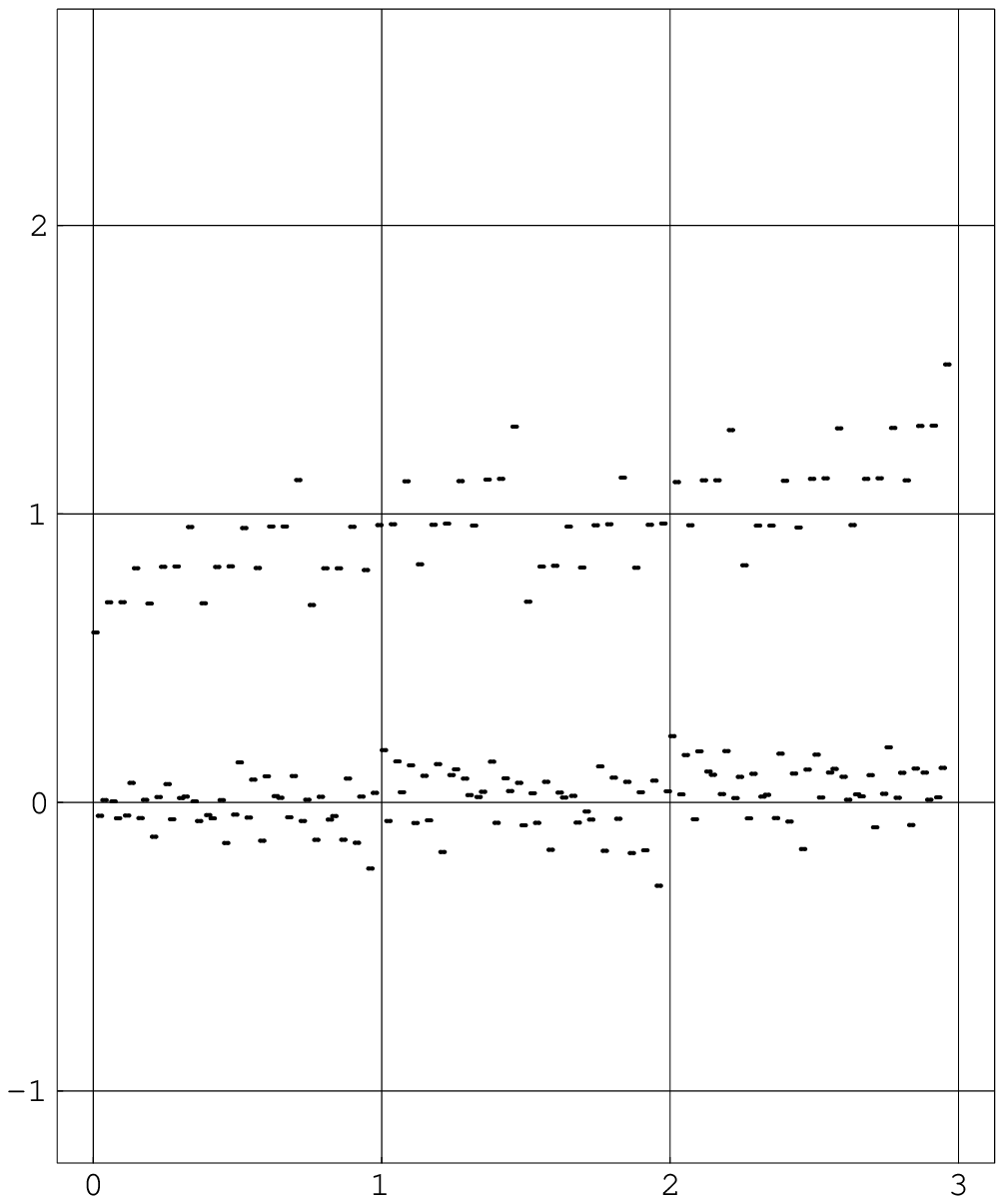}}
\put(0,240){\makebox(172,12){e: Cascade stage $5$}}
\put(172,240){\makebox(172,12){f: Cascade stage $6$}}
\put(0,18){\includegraphics
[bb=29 0 316 345,height=207bp]{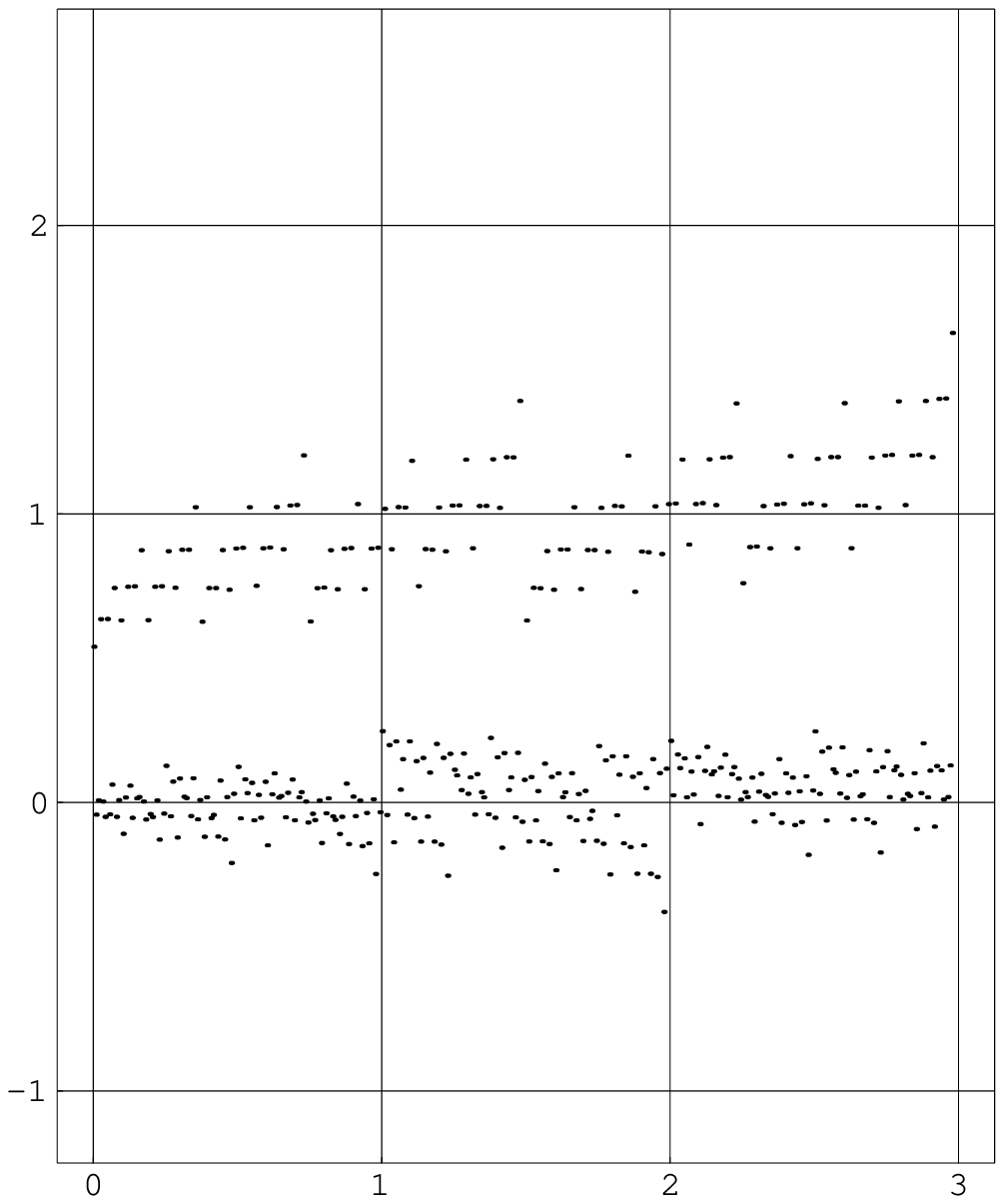}}
\put(172,18){\includegraphics
[bb=29 0 316 345,height=207bp]{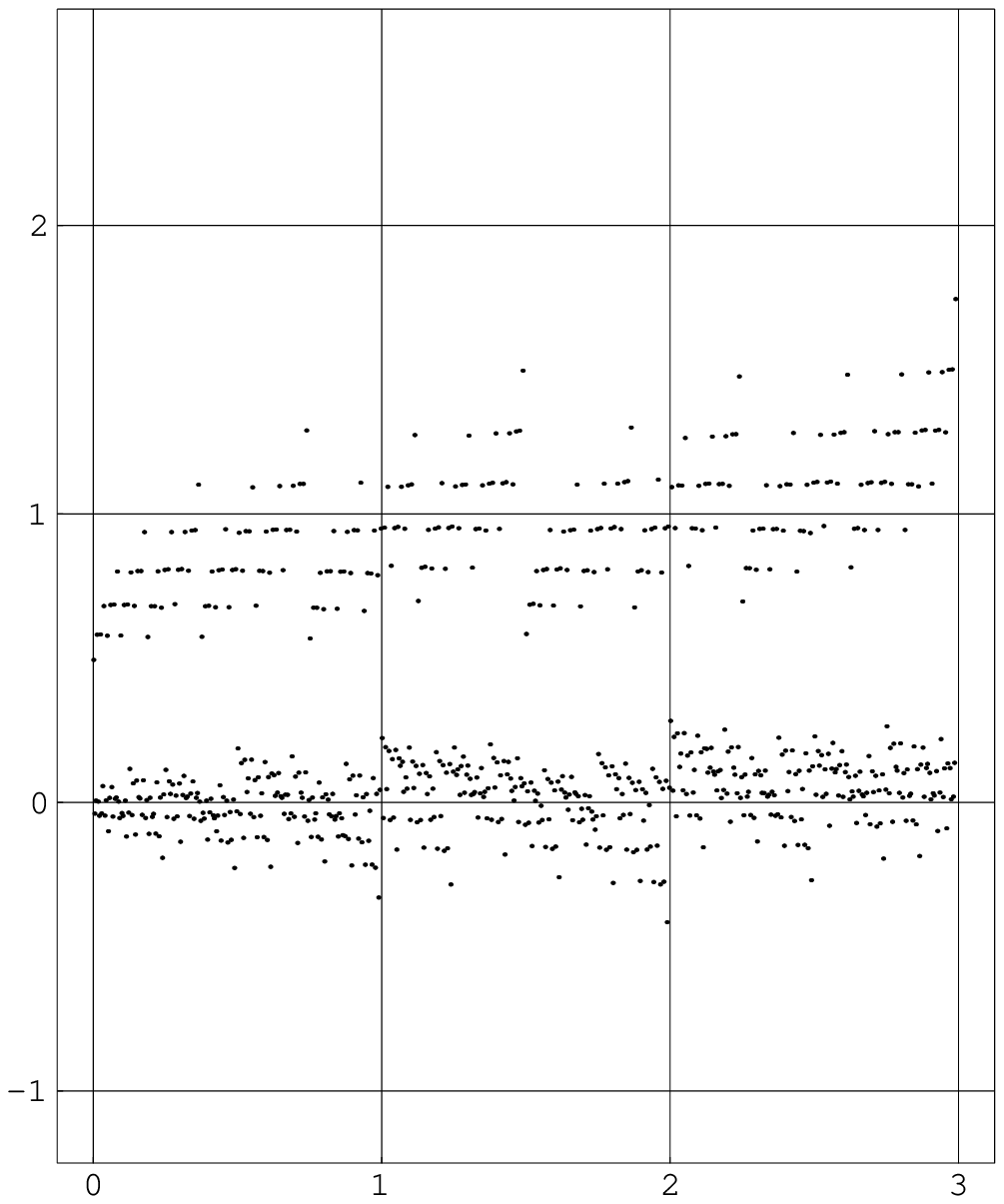}}
\put(0,3){\makebox(172,12){g: Cascade stage $7$}}
\put(172,3){\makebox(172,12){h: Cascade stage $8$}}
\put(172,-9){\makebox(172,12){(This is another version of
Figure \ref{Movie}(t).)}}
\end{picture}
\caption{Cascade stages of scaling function, $\theta=\frac{9\pi}{20}$: Stages
$5$--$8$}%
\end{figure}

\addtocounter{figure}{-1}

\begin{figure}[ptb]
\begin{picture}(344,462)
\put(0,255){\includegraphics
[bb=163 209 448 552,height=207bp]{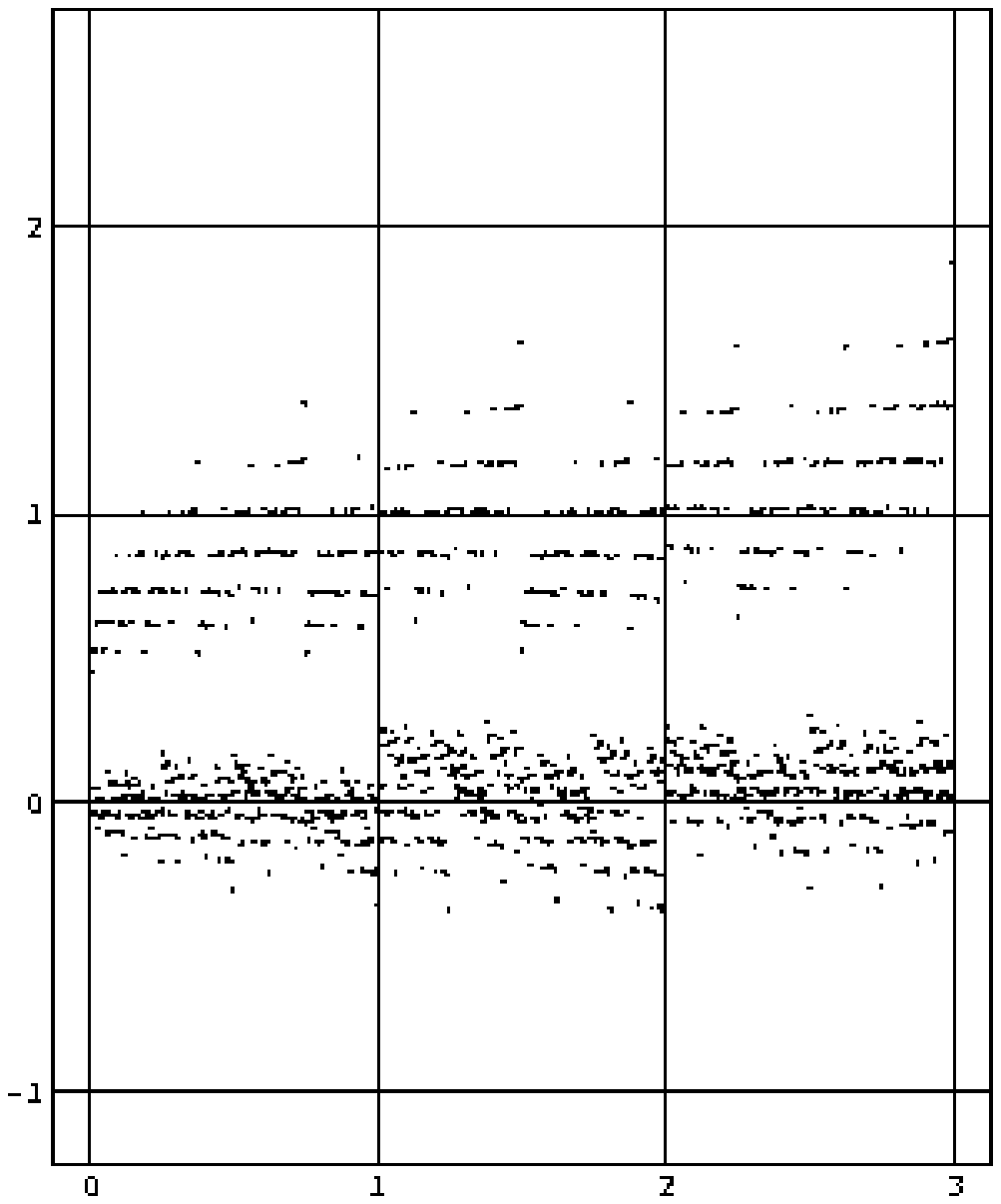}}
\put(172,255){\includegraphics
[bb=163 209 448 552,height=207bp]{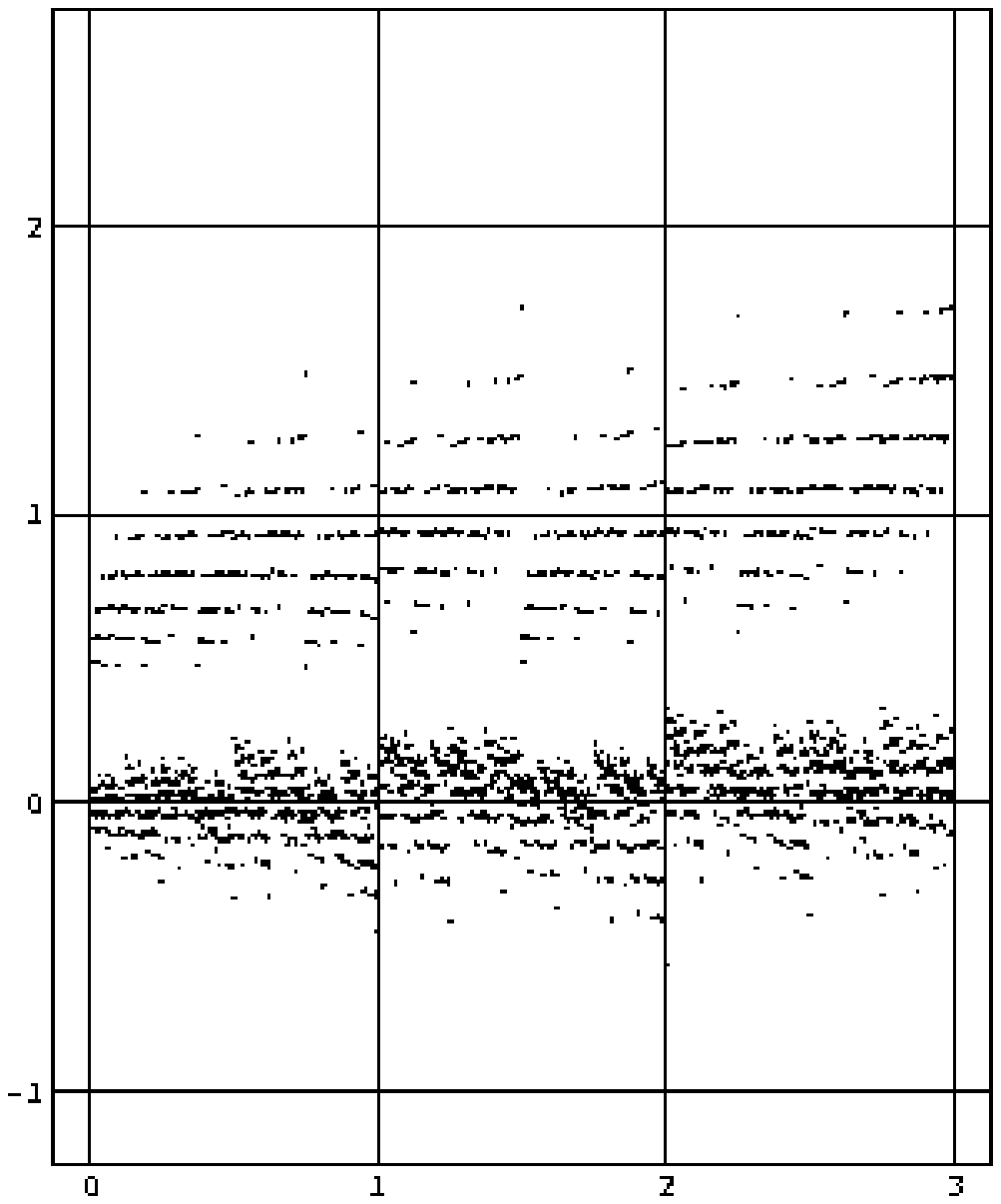}}
\put(0,240){\makebox(172,12){i: Cascade stage $9$}}
\put(172,240){\makebox(172,12){j: Cascade stage $10$}}
\put(0,18){\includegraphics
[bb=163 209 448 552,height=207bp]{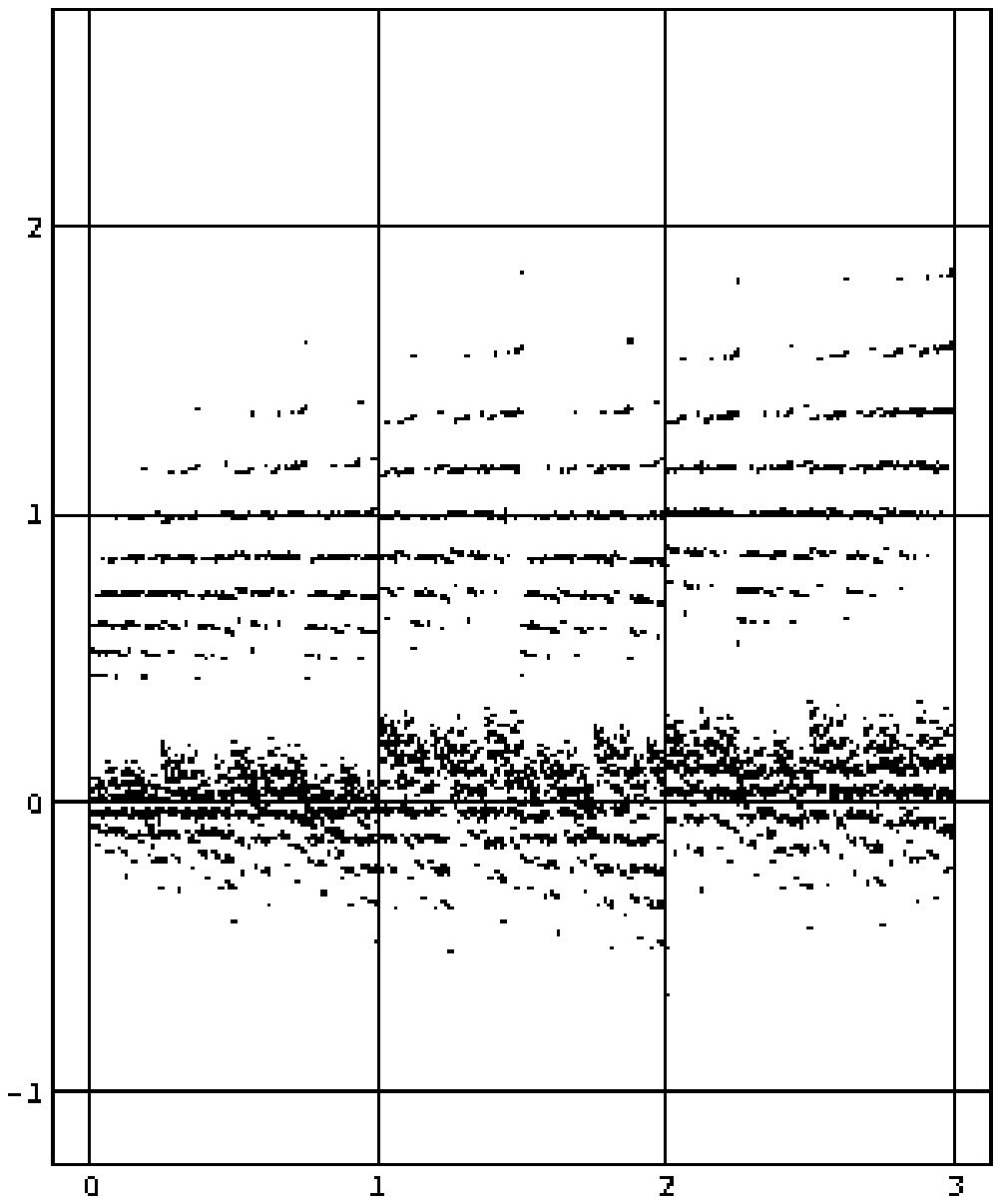}}
\put(172,18){\includegraphics
[bb=163 209 448 552,height=207bp]{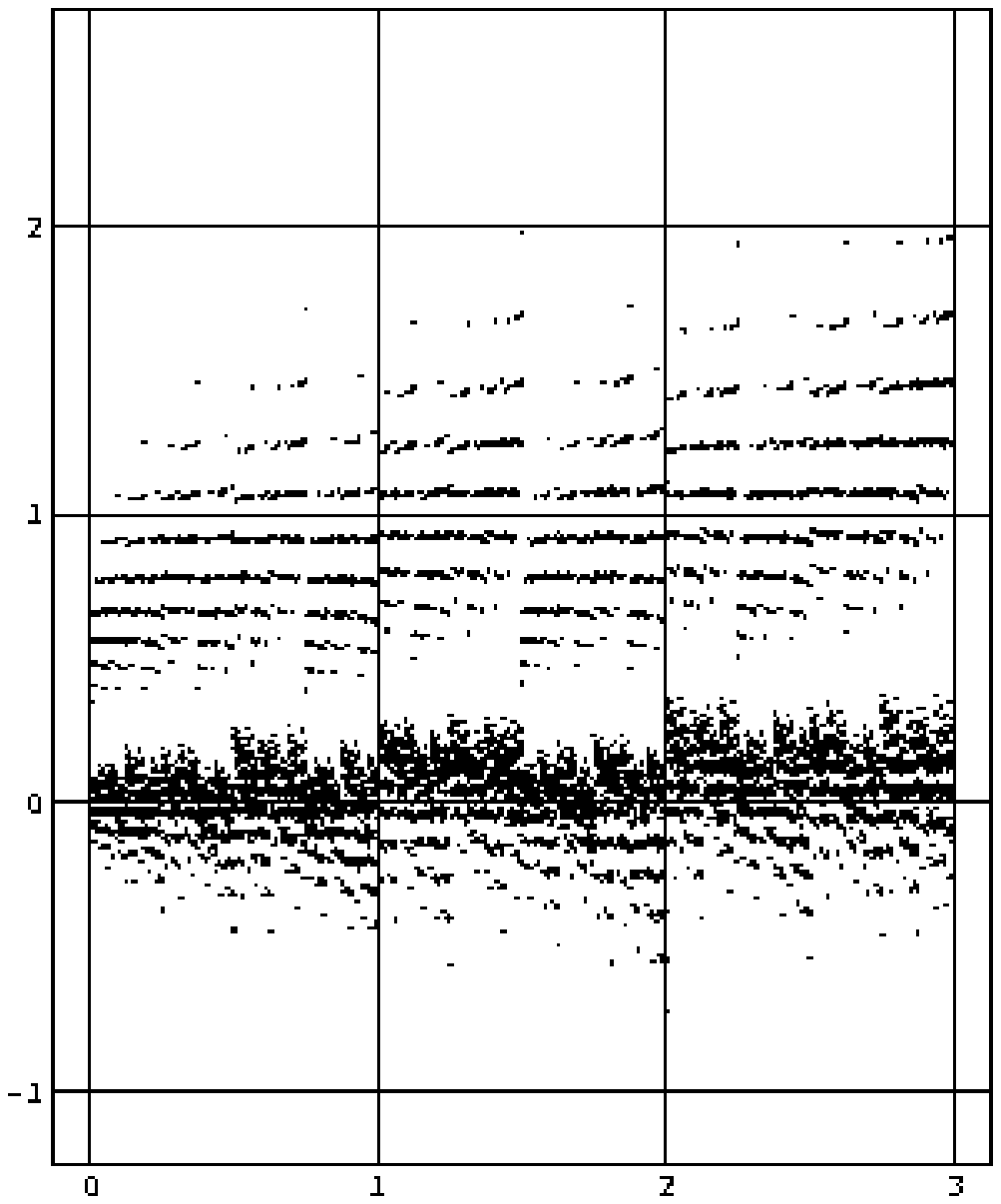}}
\put(0,3){\makebox(172,12){k: Cascade stage $11$}}
\put(172,3){\makebox(172,12){l: Cascade stage $12$}}
\end{picture}
\caption{Cascade stages of scaling function, $\theta=\frac{9\pi}{20}$: Stages
$9$--$12$}%
\label{Res045_stages_end}%
\end{figure}

\begin{figure}[ptb]
\includegraphics
[bb=90 250 521 511,width=354bp]{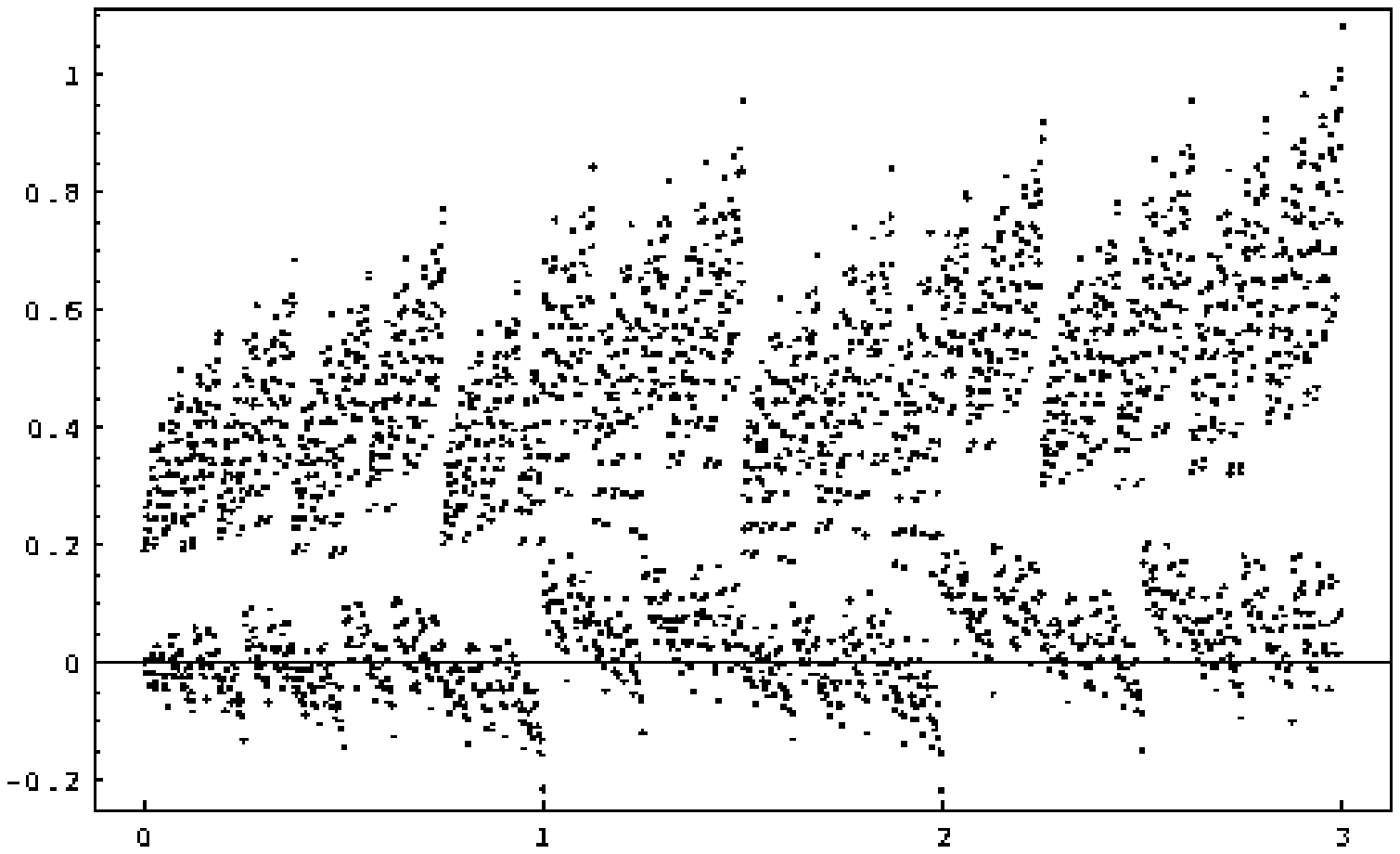}
\caption{The figure above shows
$\psi_{+}^{(1000)}\left( n\right) $ at 
$n=0,1,\dots ,3\cdot 2^{10}$,
$x_{n}=n\cdot 2^{-10}$ ($\theta =\frac{9\pi}{20}$).}%
\label{MPsiPlus}%
\end{figure}

\begin{figure}[ptb]
\includegraphics
[bb=90 250 521 511,width=354bp]{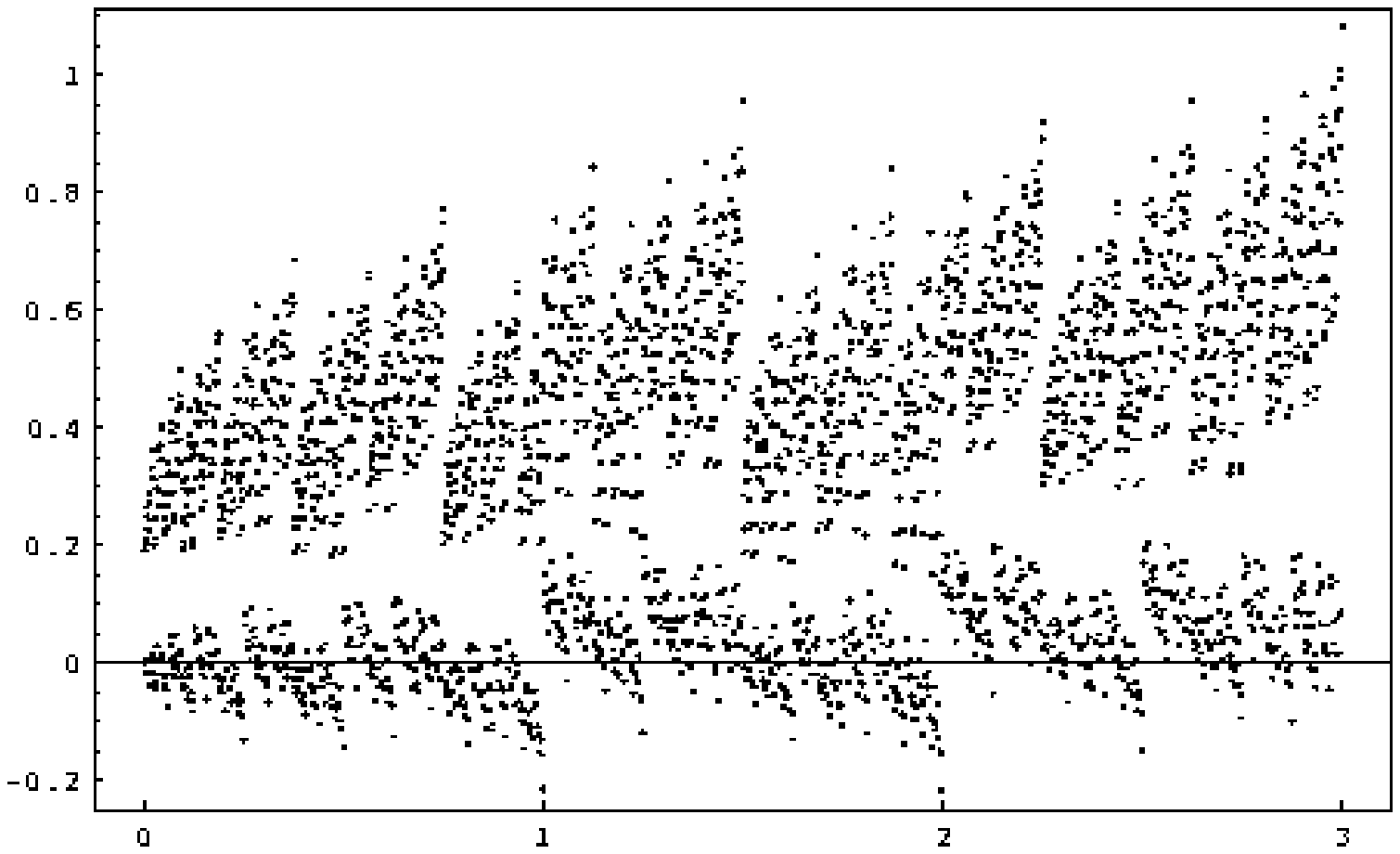}
\caption{The figure above shows
$\psi_{-}^{(1000)}\left( n\right) $ at 
$n=0,1,\dots ,3\cdot 2^{10}$,
$x_{n}=n\cdot 2^{-10}$ ($\theta =\frac{9\pi}{20}$).}%
\label{MPsiMinus}%
\end{figure}

\begin{figure}[ptb]
\begin{picture}(354,217)
\put(-24,0){\includegraphics
[bb=90 258 521 504,width=378bp]{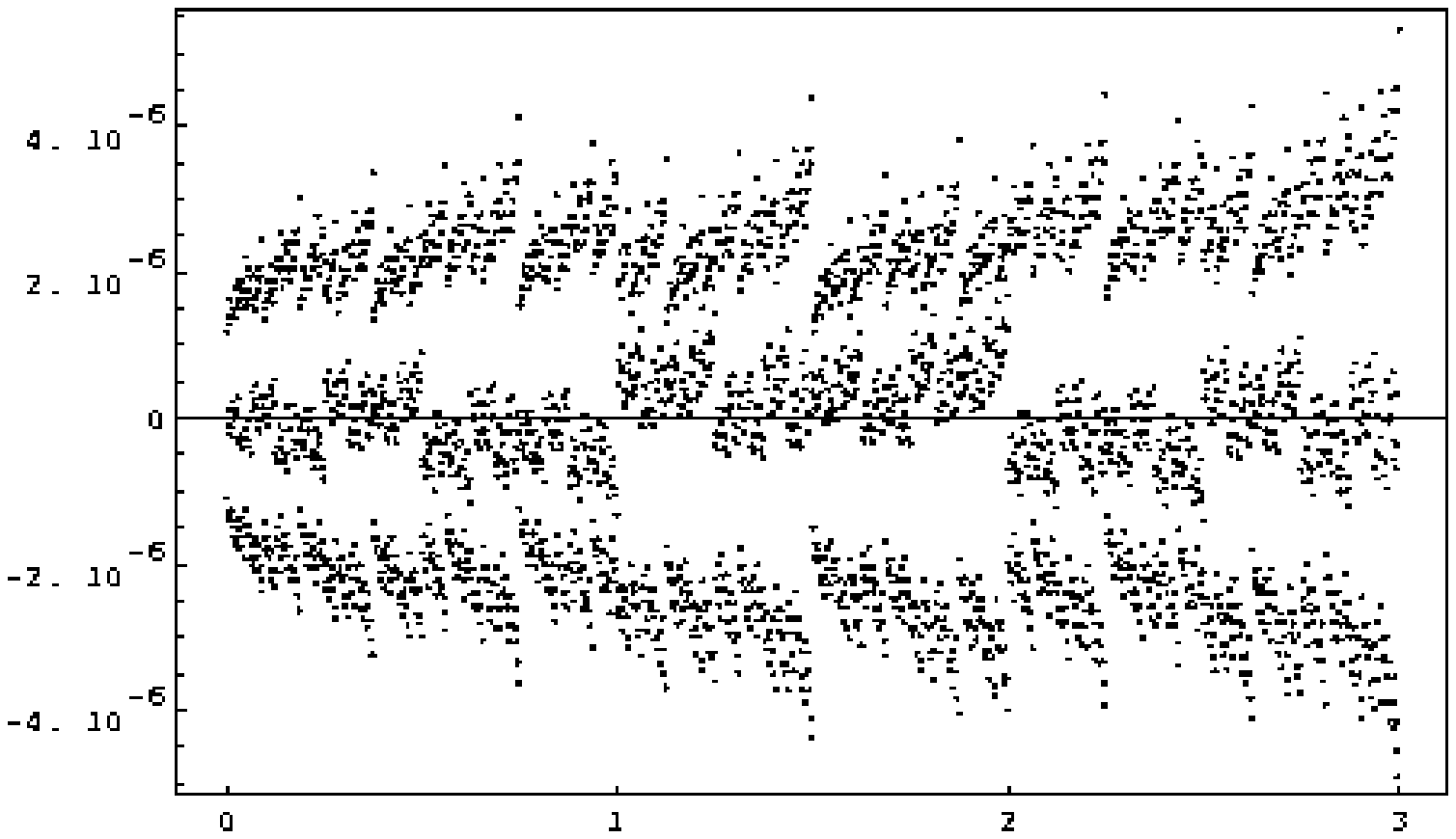}}
\end{picture}
\caption{The figure above shows
$\left( 
\psi_{+}^{(1000)}\left( n\right) 
-\psi_{-}^{(1000)}\left( n\right) 
\right) $ at 
$n=0,1,\dots ,3\cdot 2^{10}$,
$x_{n}=n\cdot 2^{-10}$.
Note the small scale at the $y$-axis, showing
$\left|
\psi_{+}^{(1000)}\left( n\right) 
-\psi_{-}^{(1000)}\left( n\right) 
\right| <\frac{5}{1000000}$. See the Appendix
for a discussion of this convergence of the differences to zero.}%
\label{MPsiDiff}%
\end{figure}

\begin{figure}[ptb]
\includegraphics
[bb=90 250 521 511,width=354bp]{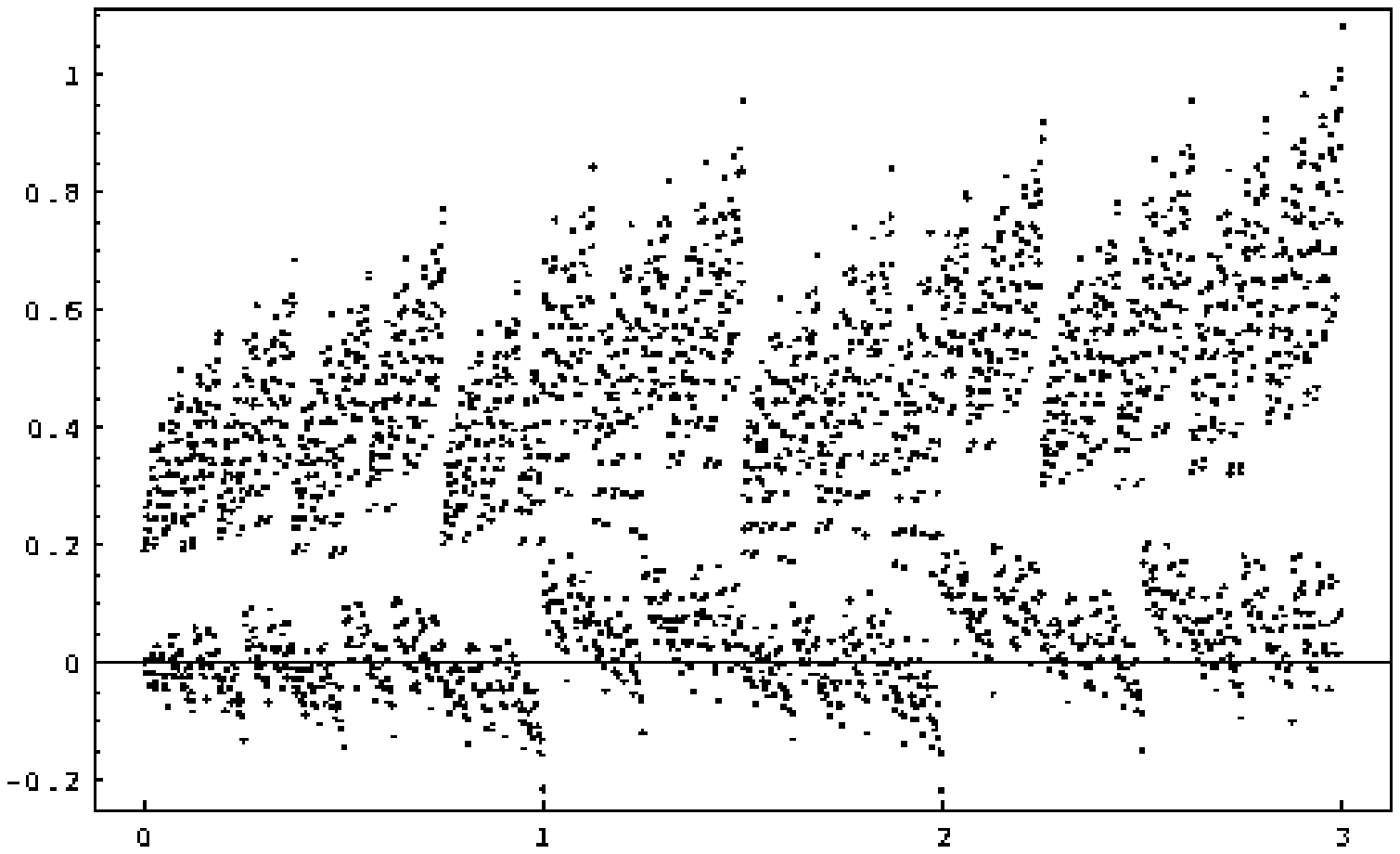}
\caption{The figure above shows
$\psi_{\pm}^{(\infty)}\left( n\right) $ at 
$n=0,1,\dots ,3\cdot 2^{10}$,
$x_{n}=n\cdot 2^{-10}$ ($\theta =\frac{9\pi}{20}$).}%
\label{InfinityPsi}%
\end{figure}

\clearpage

\section{\label{Con}Conclusions}

While it may be difficult to discern an overall pattern in the
computer-generated output from the cascades, the expectation is that there are
domains of starting points (functions or coefficients) which lead to ``nice''
limit functions, $L^{2}\left(  \mathbb{R}\right)  $ or continuous, while the
other extreme
ones lead
to fractal-like pictures. This is based on two analogies,
in addition to the existing results in
\cite{DaLa92},
\cite{CoHe92},
\cite{CoHe94},
\cite{Wan95},
\cite{Wan96}:

(i)~First, we think of (\ref{eq1}) as a version of an iterated function system
in the sense of Hutchinson \cite{Hut81}. Hutchinson considers Borel mappings
$f_{1},\dots,f_{n}$ on a
complete
metric space $X$, probabilities $p_{1},\dots,p_{n}$,
$p_{i}>0$, $\sum_{i}p_{i}=1$. Such a system defines a dynamical system
$x_{0}\rightarrow x_{1}\rightarrow\cdots$ in $X$ where%
\[
x_{i+1}=f_{\alpha_{i}}\left(  x_{i}\right)
\]
and where the indices $\alpha_{i}$ are chosen randomly for each $i$, with
probability $p_{\alpha_{i}}$. For Borel probability measures $\mu$ on $X$, set
$M_{i}\left(  \mu\right)  =\mu\circ f_{i}^{-1}$, i.e., $M_{i}\left(
\mu\right)  \left(  E\right)  =\mu\left(  f_{i}^{-1}\left(  E\right)  \right)
$, where $E\subset X$ is a Borel set, and $f_{i}^{-1}\left(  E\right)
=\left\{  x\in X\mid f_{i}\left(  x\right)  \in E\right\}  $. A fixed point
for the system is a measure $\mu$ such that%
\begin{equation}
\mu=\sum_{i}p_{i}M_{i}\left(  \mu\right)  , \label{eqCon.1}%
\end{equation}
and we think of (\ref{eq1}) as a version of this, but of course, in
(\ref{eq1}) we do not
necessarily
have the coefficients $a_{i}$ positive, and we do not
impose the same normalization. Nonetheless, we may take $X=\mathbb{R}$,
$f_{i}\left(  x\right)  :=\frac{1}{2}\left(  x+i\right)  $, and $p_{i}%
:=\frac{1}{\sqrt{2}}a_{i}$, $i=0,1,\dots,N$.
(Figure \ref{ckplot} shows examples with four
coefficients, where generically three are positive
and one is negative.)
Hutchinson's theorem for the
general version of (\ref{eqCon.1}) yields existence and uniqueness of $\mu$
provided the probabilistic assumptions hold, and the mappings $f_{i}$ have
contractive Lipschitz constants. Even in the wavelet setting, one may ask for
a
signed
measure $\mu$ solving (\ref{eqCon.1}) and make the distinction between
solutions $\mu$ which are absolutely continuous with respect to Lebesgue
measure, vs.\ the measures $\mu$ which are singular. In the first case, there
is a Radon-Nikodym derivative $\varphi$ which may be viewed as a solution to
the original problem (\ref{eq1}).

(ii)~The second analogy is to a problem studied by P.~Erd\"os \cite{Erd40} for
an iterated function system on
\begin{equation}
\begin{aligned}
X &=\left[  0,1\right]  ,\qquad 0<\lambda<1, \\
f_{0}\left(  x\right)  &=\lambda x, \\
f_{1}\left(  x\right)  &=\lambda x+1-\lambda, \text{\qquad and} \\
p_{0} &=p_{1}=\frac{1}{2}. 
\end{aligned}%
\label{eqCon.2}%
\end{equation}
The corresponding problem
(\ref{eqCon.1}) 
and
(\ref{eqCon.2}) 
leads to a probability measure $\mu_{\lambda}$, and Erd\"os
showed that $\mu_{\lambda}$ is either absolutely continuous, or else totally
singular. Total singularity means that there is a subset $E\subset X$
($=\left[  0,1\right]  $) such that $\mu_{\lambda}\left(  E\right)  =1$ and
$E$ is of zero Lebesgue measure. This shows up in computer output as
fractal-like appearance. It is known that if $\lambda<\frac{1}{2}$, then the
support of $\mu_{\lambda}$ is a Cantor set of Hausdorff dimension $-\frac
{\ln2}{\ln\lambda}$, while the case $\lambda>\frac{1}{2}$ is not fully
understood, the expectation being that $\mu_{\lambda}$ will be more regular
(or ``less fractal'') if $\lambda$ is closer to $1$.

The case $\lambda=\frac12$ in the Erd\"os
construction is clear, of course,
and yields $d\mu_{\frac12}\left(x\right)=\varphi\left(x\right)\,dx$
where $\varphi =\psi^{(0)}$ is precisely the
Haar function from (\ref{eq13}) above, and $dx$ is Lebesgue measure,
i.e., $\varphi$ satisfies
\[
\varphi\left(\frac x2\right) =\varphi \left(x\right) +\varphi \left(x-1\right),
\]
which is also (\ref{eq1}) for the Haar
scaling function.

It should also be mentioned that
Solomyak \cite{Sol95} proved that for
almost all $\lambda$ in $\left( \frac{1}{2},1\right) $, the
measure $\mu_{\lambda}$ does have an $L^{2}\left( \mathbb{R}\right) $
density. An example of $\lambda >\frac{1}{2}$ when $\mu_{\lambda}$
is known to be singular is
$\lambda^{-1}=\left( 1+\sqrt{5}\right) /2$, the golden ratio.

\section*{Appendix\protect\frogleg\mdseries\upshape\protect\frogleft
by Brian F. Treadway\protect\frogright}

\setcounter{equation}{0}\renewcommand{\theequation}{A.\arabic{equation}}
By the method of Daubechies \cite[Section 6.5, pp.\ 204--206]{Dau92}, the
iteration (\ref{eqScheme.starstar}) can be replaced by a ``local'' iteration
in which the two constant values of $\psi^{\left(  m\right)  }\left(
x\right)  $ on adjacent intervals at the left of the point $x=k\cdot
2^{-\left(  m-1\right)  }$,
\begin{align*}
\psi_{2k-2}^{\left(  m\right)  }  &  :=\psi^{\left(  m\right)  }\left(
x\right)  \text{\qquad for }x\in\left[  \left(  2k-2\right)  \cdot
2^{-m},\left(  2k-1\right)  \cdot2^{-m}\right)  ,\\
\psi_{2k-1}^{\left(  m\right)  }  &  :=\psi^{\left(  m\right)  }\left(
x\right)  \text{\qquad for }x\in\left[  \left(  2k-1\right)  \cdot
2^{-m},2k\cdot2^{-m}\right)  ,
\end{align*}
are determined as a linear combination of two constant values of
$\psi^{\left(  m-1\right)  }\left(  x\right)  $ on adjacent intervals at the
left of the same $x$ from the previous iteration, viz.\
\begin{align*}
\psi_{k-2}^{\left(  m-1\right)  }  &  :=\psi^{\left(  m\right)  }\left(
x\right)  \text{\qquad for }x\in\left[  \left(  k-2\right)  \cdot2^{-\left(
m-1\right)  },\left(  k-1\right)  \cdot2^{-\left(  m-1\right)  }\right)  ,\\
\psi_{k-1}^{\left(  m-1\right)  }  &  :=\psi^{\left(  m\right)  }\left(
x\right)  \text{\qquad for }x\in\left[  \left(  k-1\right)  \cdot2^{-\left(
m-1\right)  },k\cdot2^{-\left(  m-1\right)  }\right)  .
\end{align*}
The relation is%
\begin{equation}%
\begin{pmatrix}
\psi_{2k-2_{\mathstrut}}^{\left(  m\right)  }\\
\psi_{2k-1}^{\left(  m\right)  }%
\end{pmatrix}
=A%
\begin{pmatrix}
\psi_{k-2_{\mathstrut}}^{\left(  m-1\right)  }\\
\psi_{k-1}^{\left(  m-1\right)  }%
\end{pmatrix}
\label{eqA.1}
\end{equation}
with%
\[
A=%
\begin{pmatrix}
\sqrt{2}a_{2} & \sqrt{2}a_{0}\\
\sqrt{2}a_{3} & \sqrt{2}a_{1}%
\end{pmatrix}%
,
\]
where the $a_{i}$ are as in (\ref{eqSome.2}).
With the definition
\begin{equation}
\xi_{m-1+i}:=
\begin{pmatrix}
\psi_{\left(  k\cdot2^{i}\right)  -2_{\mathstrut}}^{\left(  m-1+i\right)  }\\
\psi_{\left(  k\cdot2^{i}\right)  -1}^{\left(  m-1+i\right)  }%
\end{pmatrix}%
,
\label{eqA.2}
\end{equation}
the relation (\ref{eqA.1}) becomes
\[
\xi_{m}=A\xi_{m-1}.
\]

The matrix relation (\ref{eqA.1}) is in fact
a transcription of the equations (6.5.9)
and (6.5.10) in \cite{Dau92}, and the
reader should be warned that the
relations are not merely an algebraic
consequence of the cascade iteration
(\ref{eq9}), but depend heavily on the orthogonality
of the translates of the
successive iterates $\psi^{\left( m\right) }$.
See Remark \ref{RemOrth} for comments on this dependence
on orthogonality.
The relation
between Daubechies's iterates and ours
is that
\[
\psi_{+}^{\left( m\right) }\left( n\right) 
=\left( M^{m}\chi^{}_{\left[ 0,1\right\rangle}\right) 
\left( n\cdot 2^{-N}\right) 
=\psi_{n\cdot 2^{-N+m}}^{\left( m\right) },
\]
at least when $m\geq M$. Thus, when
Daubechies's scheme is iterated more than
$N$ times, the two methods give the
same values on $\mathbb{Z}2^{-N}$. In order
to get the values of $\varphi$ at dyadic
rationals it is necessary to switch
from Daubechies's algorithm, which
doubles the number of points at
each iteration, to the algorithm based
on (\ref{eq9}), where the number of points
is constant. The extreme cusps in
the movie reel Figure \ref{Movie}(b)--(j) are therefore
computed by using (\ref{eq9}) with large $n$'s,
or by computing the limit $n\rightarrow \infty $ as in
(\ref{eqA.4}) below.
The connection between
the two versions of the
cascade algorithm is
spelled out in \cite[p.~204]{Dau92},
where equation (6.5.3) shows
that the value of
$\varphi$ at the dyadic rational $n\cdot 2^{-N}$ may be
approximated by
$\psi^{\left( j\right) }_{2^{j-N}}$ as $j\rightarrow\infty$ where
$\psi^{\left( j\right) }_{m}=2^{j/2}\ip{\varphi }{\varphi_{-j,m}}$ and
$\varphi_{j,m}\left( x\right) 
=2^{-j/2}\varphi\left( 2^{-j}x-m\right) $.
When the approximations at a given stage
are written out in terms of products of
the coefficients $a_{i}$ in the order they
arise in the iteration, it is found
that the two algorithms supply these factors
in precisely the opposite order.

If $\sin\theta\neq-1$, the eigenvectors of $A=\left(
\begin{smallmatrix}
\sqrt{2}a_{2} & \sqrt{2}a_{0}\\
\sqrt{2}a_{3} & \sqrt{2}a_{1}%
\end{smallmatrix}
\right)  $ are: for eigenvalue $\lambda_{1}=1$, $e_{1}=\left(
\begin{smallmatrix}
1\\
1
\end{smallmatrix}
\right)  $; for eigenvalue $\lambda_{2}=-\sin\theta$, $e_{2}=\left(
\begin{smallmatrix}
\sqrt{2}a_{1}\\
\sqrt{2}a_{2}%
\end{smallmatrix}
\right)  $. A starting vector $\xi_{m-1}=\left(
\begin{smallmatrix}
\psi_{k-2}^{\left(  m-1\right)  }\\
\psi_{k-1}^{\left(  m-1\right)  }%
\end{smallmatrix}
\right)  $ composed of ordinates in the $\left(  m-1\right)  $'st cascade
approximant gives, by $n$ applications of $A$,
\[
A^{n}\xi_{m-1}
=\xi_{m-1+n},
\]
where $\xi_{m-1+n}$, given by (\ref{eqA.2}) above with $i=n$,
is
composed of ordinates in the $\left(  m-1+n\right)
$'th cascade approximant on two intervals at the left of the same point
$x=k\cdot2^{-\left(  m-1\right)  }$. Expanding the starting vector $\xi_{m-1}$
in the two eigenvectors, we get%
\begin{align*}
\xi_{m-1} &=\alpha_{1}e_{1}+\alpha_{2}e_{2}  \\
&=\frac{\sqrt{2}a_{2}\psi
_{k-2}^{\left(  m-1\right)  }-\sqrt{2}a_{1}\psi_{k-1}^{\left(  m-1\right)  }%
}{\sqrt{2}\left(  a_{2}-a_{1}\right)  }e_{1}  \\
&\qquad +\frac{-\psi_{k-2}^{\left(
m-1\right)  }+\psi_{k-1}^{\left(  m-1\right)  }}{\sqrt{2}\left(  a_{2}%
-a_{1}\right)  }e_{2}.
\end{align*}
Thus
\begin{align}
\xi_{m-1+n} &  =A^{n}\xi_{m-1}  \label{eqA.3}\\
&=\lambda_{1}^{n}\alpha_{1}e_{1}+\lambda_{2}%
^{n}\alpha_{2}e_{2} \notag \\
&  =\frac{\sqrt{2}a_{2}\psi_{k-2}^{\left(  m-1\right)  }-\sqrt{2}a_{1}%
\psi_{k-1}^{\left(  m-1\right)  }}{\sqrt{2}\left(  a_{2}-a_{1}\right)  }%
e_{1}  \notag \\
&\qquad +\left(  -\sin\theta\right)  ^{n}\frac{-\psi_{k-2}^{\left(  m-1\right)
}+\psi_{k-1}^{\left(  m-1\right)  }}{\sqrt{2}\left(  a_{2}-a_{1}\right)
}e_{2}. \notag
\end{align}
This form allows us to read off the $n\rightarrow\infty$ limit directly except
in the two special cases $\theta=\pm\frac{\pi}{2}$, as the second term
on the right-hand side of (\ref{eqA.3})
vanishes in the limit
unless $\sin\theta=-1$ (already excluded) or $\sin\theta=1$:%
\begin{align}
\lim_{n\rightarrow\infty}\xi_{n}
&=\lim_{n\rightarrow\infty}\xi_{m-1+n}  \label{eqA.4}\\
&=\frac
{1}{\sqrt{2}\left(  a_{2}-a_{1}\right)  }%
\begin{pmatrix}
\sqrt{2}a_{2}\psi_{k-2}^{\left(  m-1\right)  }-\sqrt{2}a_{1}\psi
_{k-1}^{\left(  m-1\right)  }\\
\sqrt{2}a_{2}\psi_{k-2}^{\left(  m-1\right)  }-\sqrt{2}a_{1}\psi
_{k-1}^{\left(  m-1\right)  }%
\end{pmatrix}
.
\notag
\end{align}

Note that the two components of this limit
vector are equal, so the ``jump'' in $\psi$
disappears at the dyadic rational point $x$ in question.
This does not imply that $\psi$ converges to a
continuous function $\varphi$. While the values of
$\psi^{\left(  m-1+n\right)  }\left(  x\right)  $ on the two intervals in the
stage-$\left(  m-1+n\right)  $ partition at the immediate left of the $x$ in
question do move closer together as $n$ increases until their difference
becomes zero in the limit, values at higher-order dyadic rationals are not
so well behaved: for example, we can see, by extending the treatment above to a
$3\times3$ matrix $A^{\prime}=\left(
\begin{smallmatrix}
\sqrt{2}a_{3} & \sqrt{2}a_{1} & 0\\
0 & \sqrt{2}a_{2} & \sqrt{2}a_{0}\\
0 & \sqrt{2}a_{3} & \sqrt{2}a_{1}%
\end{smallmatrix}
\right)  $ iterating three intervals at the left of the point $x$ instead of
just two, that the value of $\psi^{\left(  m-1+n\right)  }\left(  x\right)  $
on the third interval
added at the left of the original two
generally grows without limit as $n\rightarrow\infty$
if the new eigenvalue $\sqrt{2}a_{3}$, corresponding to the eigenvector
$\left(
\begin{smallmatrix}
1\\
0\\
0
\end{smallmatrix}
\right)  $, is greater than $1$, which is the case for $0<\theta<\frac{\pi}%
{2}$.

As an example, the negative ``peak'' at $x=1$ (see Figure \ref{Movie}(b)--(j))
can be obtained by using the starting vector $\xi_{0}=\left(
\begin{smallmatrix}
\psi_{-1}^{\left(  0\right)  }\\
\psi_{0}^{\left(  0\right)  }%
\end{smallmatrix}
\right)  =\left(
\begin{smallmatrix}
0\\
1
\end{smallmatrix}
\right)  $ (following the convention that $\psi^{\left(  m\right)  }\left(
x\right)  $ for $x<0$ for all $m$, and thus $\psi_{-1}^{\left(  0\right)  }%
=0$; see Figure \ref{Res045_00}): then the two components of $\lim_{n\rightarrow\infty}\xi_{n}$ are both
equal to%
\begin{align*}
-\frac{\sqrt{2}a_{1}}{\sqrt{2}\left(  a_{2}-a_{1}\right)  }&=-\frac
{1-\cos\theta-\sin\theta}{2\cos\theta}  \\
&=\frac{1}{2}\left(  1+\tan\frac
{\theta-\frac{\pi}{2}}{2}\right)  .
\end{align*}
Similarly, the peak at $x=2$ in the same plots can be obtained by using the
starting vector $\xi_{0}=\left(
\begin{smallmatrix}
\psi_{0}^{\left(  0\right)  }\\
\psi_{1}^{\left(  0\right)  }%
\end{smallmatrix}
\right)  =\left(
\begin{smallmatrix}
1\\
0
\end{smallmatrix}
\right)  $ (see Figure \ref{Res045_00}); it is
\begin{align*}
\frac{\sqrt{2}a_{2}}{\sqrt{2}\left(  a_{2}-a_{1}\right)  }&=\frac{1+\cos
\theta-\sin\theta}{2\cos\theta}  \\
&=\frac{1}{2}\left(  1-\tan\frac{\theta
-\frac{\pi}{2}}{2}\right)  .
\end{align*}
Note that these two peaks sum to $1$ for all $\theta$ (excluding $\pm\frac
{\pi}{2}$). The limiting value at $x=3/2$ is also independent of $\theta$. For
this, use the starting vector $\xi_{1}=\left(
\begin{smallmatrix}
\psi_{1}^{\left(  1\right)  }\\
\psi_{2}^{\left(  1\right)  }%
\end{smallmatrix}
\right)  =\left(
\begin{smallmatrix}
\sqrt{2}a_{1}\\
\sqrt{2}a_{2}%
\end{smallmatrix}
\right)  $ (see Figure \ref{Res045_stages}(a)). This is just $e_{2}$, so
\[
\lim_{n\rightarrow\infty}\xi_{n}=%
\begin{pmatrix}
0\\
0
\end{pmatrix}
.
\]
The numerical values of $\lim_{n\rightarrow\infty}\psi^{\left(  n\right)
}\left(  x\right)  $, $x=1,\frac{3}{2},2$ for the ``movie reel'' plots in
Figure \ref{Movie} are given below.
\[%
\begin{tabular}
[c]{rc|rcr}%
Fig.\ \ref{Movie} & $\theta$ & $x=1$ & $x=\frac{3}{2_{\mathstrut}}$ & $x=2$\\\hline
(b) &  $-\frac{9\pi^{\mathstrut}}{20_{\mathstrut}}$ & $-5.8531$ & $0$ & $6.8531$\\
(c) &  $-\frac{2\pi}{5_{\mathstrut}}$ & $-2.6569$ & $0$ & $3.6569$\\
(d) &  $-\frac{7\pi}{20_{\mathstrut}}$ & $-1.5826$ & $0$ & $2.5826$\\
(e) &  $-\frac{3\pi}{10_{\mathstrut}}$ & $-1.0388$ & $0$ & $2.0388$\\
(f) &  $-\frac{\pi}{4_{\mathstrut}}$ & $-0.7071$ & $0$ & $1.7071$\\
(g) &  $-\frac{\pi}{5_{\mathstrut}}$ & $-0.4813$ & $0$ & $1.4813$\\
(h) &  $-\frac{3\pi}{20_{\mathstrut}}$ & $-0.3159$ & $0$ & $1.3159$\\
(i) &  $-\frac{\pi}{10_{\mathstrut}}$ & $-0.1882$ & $0$ & $1.1882$\\
(j) &  $-\frac{\pi}{20_{\mathstrut}}$ & $-0.0854$ & $0$ & $1.0854$\\
(k) &  $0\vphantom{\frac{\pi}{20_{\mathstrut}}}$ & $0.0000$ & $0$ & $1.0000$\\
(l) &  $\frac{\pi}{20_{\mathstrut}}$ & $0.0730$ & $0$ & $0.9270$\\
(m) &  $\frac{\pi}{10_{\mathstrut}}$ & $0.1367$ & $0$ & $0.8633$\\
(n) &  $\frac{3\pi}{20_{\mathstrut}}$ & $0.1936$ & $0$ & $0.8064$\\
(o) &  $\frac{\pi}{5_{\mathstrut}}$ & $0.2452$ & $0$ & $0.7548$\\
(p) &  $\frac{\pi}{4_{\mathstrut}}$ & $0.2929$ & $0$ & $0.7071$\\
(q) &  $\frac{3\pi}{10_{\mathstrut}}$ & $0.3375$ & $0$ & $0.6625$\\
(r) &  $\frac{7\pi}{20_{\mathstrut}}$ & $0.3800$ & $0$ & $0.6200$\\
(s) &  $\frac{2\pi}{5_{\mathstrut}}$ & $0.4208$ & $0$ & $0.5792$\\
(t) &  $\frac{9\pi}{20}$ & $0.4606$ & $0$ & $0.5394$%
\end{tabular}
\]
These values can be observed in the plots, though for $\theta$ near $
\frac{\pi}{2}$ the cascade level plotted ($m=8$) is not high enough to show a
close approach to the limit.

\begin{acknowledgements}
We are grateful to Brian Treadway for insightful observations on the cascade
implementations, and to Rune Kleveland for help with implementing the algorithm.
\end{acknowledgements}

Most of the work in the present paper was done while O.B. visited Iowa in
December of 1997, but it was interrupted at the end of the visit, when P.J.
landed in the hospital after an accident.

\bibliographystyle{bftalpha}
\bibliography{jorgen}
\end{document}